\documentclass{mcom-l}

\usepackage{amssymb}

\usepackage{amssymb}

\usepackage[utf8]{inputenc}
\usepackage{bbm}
\usepackage[bbgreekl]{mathbbol}
\usepackage{graphicx,subfig,tikz} 

\usepackage{comment}

\usepackage{tikz}
\usepackage{tikz-cd}
\usetikzlibrary{arrows,matrix}

\usepackage{mathtools}
\mathtoolsset{showonlyrefs=true}

\usepackage{multirow}

\newtheorem{theorem}{Theorem}[section]

\theoremstyle{definition}
\newtheorem{definition}[theorem]{Definition}

\newtheorem{assumption}[theorem]{Assumption}
\newtheorem{prop}[theorem]{Proposition}

\theoremstyle{remark}
\newtheorem{remark}[theorem]{Remark}

\numberwithin{equation}{section}

\usepackage[T3,T1]{fontenc}
\DeclareSymbolFont{tipa}{T3}{cmr}{m}{n}
\DeclareMathAccent{\invbreve}{\mathalpha}{tipa}{16}

\def\bs{\boldsymbol}

\def\ib{{\boldsymbol{i}}}
\def\jb{{\boldsymbol{j}}}
\def\mb{{\boldsymbol{m}}}
\def\bzero{{\boldsymbol{0}}}

\def\ba{{\bs a}}

\def\bi{{\bs i}}
\def\bj{{\bs j}}

\def\bp{{\bs p}}

\def\bu{{\bs u}}
\def\bv{{\bs v}}

\def\bx{{\bs x}}

\def\bB{{\bs B}}

\def\bE{{\bs E}}

\def\bJ{{\bs J}}

\def\bP{{\bs P}}

\def\btau{{\bs \tau}}

\def\bLambda{{\bs \Lambda}}

\def\bgamma{{\bs \gamma}}

\def\hx{{\hat x}}

\def\hbx{{\hat {\bs x}}}

\def\cC{\mathcal{C}}

\def\cE{\mathcal{E}}
\def\cF{\mathcal{F}}

\def\cK{\mathcal{K}}

\def\cT{\mathcal{T}}

\def\cV{\mathcal{V}}

\def\AA{\mathbbm{A}}

\def\NN{\mathbbm{N}}

\def\QQ{\mathbbm{Q}}
\def\RR{\mathbbm{R}}

\def\VV{\mathbbm{V}}

\usepackage{ esint }

\def\pw{\ensuremath{\mathrm{pw}}} 

\def\beq{\begin{equation}}
\def\eeq{\end{equation}}
\def\beqq{\begin{equation*}}
\def\eeqq{\end{equation*}}

\def\tb{\hbox{$\|\kern -.09em |$}}
\def\bigtb{\hbox{$\big\|\kern -.09em \big|$}}
\def\Bigtb{\hbox{$\Big\|\kern -.09em \Big|$}}

\providecommand{\sprod}[2]{\langle#1,#2\rangle}

\providecommand{\openint}[2]{\mathopen{]}#1, #2\mathclose{[}}

\def\bbb{\phantom{\big|} }

\usepackage{amsopn}

\DeclareMathOperator{\Int}{int}

\DeclareMathOperator{\Div}{div}

\DeclareMathOperator{\Span}{Span}

\DeclareMathOperator{\curl}{curl}
\DeclareMathOperator{\bcurl}{{\mathbf{curl}}}

\DeclareMathOperator{\bgrad}{{\mathbf{grad}}}

\newcommand{\arr}[1]{{\bs{\mathsf{#1}}}}
\newcommand\arre{{\bs{\mathsf{e}}}}

\newcommand\arrb{{\bs{\mathsf{b}}}}
\newcommand\arrT{{\bs{\mathsf{T}}}}
\newcommand\arrC{{\bs{\mathsf{C}}}}
\newcommand\arrP{{\bs{\mathsf{P}}}}
\newcommand\arrG{{\bs{\mathsf{G}}}}

\newcommand\arrj{{\bs{\mathsf{j}}}}

\newcommand\arrM{{\bs{\mathsf{M}}}}

\newcommand\arrE{{\bs{\mathsf{E}}}}
\newcommand\arrR{{\bs{\mathsf{R}}}}

\newcommand\arrI{{\bs{\mathsf{I}}}}

\usepackage{mathrsfs}

\def\scrE{\mathscr{E}}
\def\scrR{\mathscr{R}}

\newcommand\ttC{{\mathtt{C}}}

\newcommand\vertex{\mathrm{v}}
\newcommand\vertices{\cV}

\newcommand{\iden}{{\mathrm{id}}}

\newcommand{\edge}{\mathrm{e}}
\newcommand{\edges}{\cE}

\newcommand{\polspace}{\QQ_{r-1, r-1}[\hat \bx]}
\newcommand{\scalarpolspace}{\QQ_{r-1}[\hat x]}

\newcommand{\gammaod}{\arr{\gamma}^{\text{1d}}}

\newcommand{\coeffs}{\cC}

\newcommand{\ndof}{N_{\mathrm{dof}}}
\newcommand{\nc}{N_{\mathrm{cells}}}

\def\bcI{{\bs {\mathcal I}}}

\usepgfmodule{nonlineartransformations} 
\makeatletter
\def\latticetilt{%
\pgf@xa=\pgf@x%
\pgf@ya=\pgf@y%
\pgfmathsetmacro{\myx}{\pgf@xa+\pgfkeysvalueof{/tikz/lattice/amplitude}*sin((\pgf@ya/\pgfkeysvalueof{/tikz/lattice/spacing})*360/\pgfkeysvalueof{/tikz/lattice/superlattice period})}%
\pgf@x=\myx pt%
\pgfmathsetmacro{\myy}{\pgf@ya+\pgfkeysvalueof{/tikz/lattice/amplitude}*sin((\pgf@xa/\pgfkeysvalueof{/tikz/lattice/spacing})*360/\pgfkeysvalueof{/tikz/lattice/superlattice period})}%
\pgf@y=\myy pt}
\makeatother

\usepackage{hyperref}

\usepgfmodule{nonlineartransformations} 
\makeatletter
\def\latticetilt{%
\pgf@xa=\pgf@x%
\pgf@ya=\pgf@y%
\pgfmathsetmacro{\myx}{\pgf@xa+\pgfkeysvalueof{/tikz/lattice/amplitude}*sin((\pgf@ya/\pgfkeysvalueof{/tikz/lattice/spacing})*360/\pgfkeysvalueof{/tikz/lattice/superlattice period})}%
\pgf@x=\myx pt%
\pgfmathsetmacro{\myy}{\pgf@ya+\pgfkeysvalueof{/tikz/lattice/amplitude}*sin((\pgf@xa/\pgfkeysvalueof{/tikz/lattice/spacing})*360/\pgfkeysvalueof{/tikz/lattice/superlattice period})}%
\pgf@y=\myy pt}

\begin{document}

\title[Broken-FEEC on locally refined multipatch domains]%
	{Broken-FEEC on multipatch domains with local refinements}


\author{Frederik Schnack}
\address{Max-Planck-Institut für Plasmaphysik, Boltzmannstr. 2, 85748 Garching, Germany}
\email{frederik.schnack@ipp.mpg.de}

\author{Martin Campos Pinto}
\address{Max-Planck-Institut für Plasmaphysik, Boltzmannstr. 2, 85748 Garching, Germany}
\email{martin.campos-pinto@ipp.mpg.de}

\date{August 06, 2026}

\begin{abstract}
This article introduces a novel approach for spline broken-FEEC (Finite Element Exterior Calculus) on multipatch domains, extending its application to refined spline patches with non-matching interfaces. 
Traditional broken-FEEC allows for discontinuous discretizations at patch interfaces, preserving the de Rham structure and offering computational benefits. However, local refinements may lead to numerical artifacts at non-matching interfaces.
Our solution involves developing moment-preserving discrete conforming projection operators. These operators are explicit, localized, and metric-independent, ensuring $H^1$ 
and $H(\curl)$ continuity across non-matching interfaces while preserving high-order polynomial moments. This results in broken-FEEC de Rham sequences with accurate strong and weak derivatives, leading to energy-preserving Maxwell solvers that are explicit and virtually free of spurious modes. Numerical simulations confirm the efficacy of our method in eliminating spurious waves.
\end{abstract} 



\maketitle

\setcounter{tocdepth}{1}
\tableofcontents
\newpage 

\section{Introduction}
\label{sec:introduction}
Finite element methods that preserve the de Rham structure underlying partial differential equations (PDEs) in electromagnetics and fluid dynamics have been extensively studied in the last decades \cite{Raviart.Thomas.1977.lnm,Nedelec.1980.numa,Girault.Raviart.1986.sv,bochev_principles_2006}. 
Following the pioneering works of Bossavit \cite{Bossavit.1988.IEE-A,Bossavit.1998.ap} and Hiptmair \cite{Hiptmair.1999.mcomp,Hiptmair.2002.anum}, the Hilbert complex analysis of Arnold, Falk, and Winther~\cite{Arnold.Falk.Winther.2006.anum,Arnold.Falk.Winther.2010.bams} 
has identified finite element exterior calculus (FEEC) as a unifying framework in which 
the existence of bounded cochain projections drives the study of the intrinsic properties of these methods, such as discrete stability, spectral accuracy and structure preservation previously established 
by other authors; see e.g., \cite{Compatible.2006.IMA,Boffi.2006.csd,boffi2011discrete}.

A powerful extension of these works has been developed in the context of isogeometric analysis methods \cite{Hughes.Cottrell.Bazilevs.2005.cmame}, utilizing structure-preserving spline finite element spaces as proposed in \cite{Buffa_2011}. These discretizations employ compatible sequences of tensor-product spline spaces, initially defined on a Cartesian parametric domain and then mapped to subdomains (referred to as patches) through pullback and pushforward operators. The parametric tensor-product structure is particularly advantageous, as it facilitates the implementation of efficient numerical algorithms while enabling the discretization of 
general domains with complex geometries. 

In \cite{guclu_broken_2023,campos_pinto_broken-feec_2025}, a {\em broken-FEEC} 
framework was proposed for multipatch spline spaces, following 
previous work on nonconforming versions of N\'edel\'ec's elements for Maxwell's equations 
\cite{Campos-Pinto.Sonnendrucker.2016.mcomp}.
In this framework, PDEs are discretized using finite element spaces 
which are fully discontinuous at the patch interfaces:
local de Rham sequences are constructed on individual patches and strong differential operators are obtained by combining the single-patch (broken) derivatives with conforming projection operators which enforce the appropriate continuity conditions at the patch interfaces.
An important point is that the underlying continuity is the minimal one to guarantee 
a Hilbert de Rham subcomplex property. In particular, no high-order smoothness is required 
across patch interfaces, which allows the discrete conforming projection operators to be of geometric nature and metric independent.
In turn, broken-FEEC discretizations have the attractive feature of preserving
the de Rham structure and being equipped with local $L^2$-projection operators. 
Allowing discontinuities at the patch interfaces indeed leads to mass matrices which 
are block-diagonal, with blocks corresponding to individual patches.
Inverse mass matrices are then also block-diagonal, which results in $L^2$-projections, 
Galerkin Hodge and discrete coderivatives being all local operators
in the sense that their values on a given patch only depend on function values from that patch and its adjacent neighbors.

Although the discrete spaces are discontinuous across patch interfaces, broken-FEEC
differs from standard discontinuous Galerkin (DG) methods
\cite{Buffa.Perugia.2006.sinum,Buffa.Perugia.Warburton.2009.jsc}.
DG methods couple broken fields through numerical fluxes and, typically, jump penalties
in the variational formulation. In broken-FEEC, the discrete differentials are instead
obtained by composing patch-wise differential operators with discrete projections onto a
conforming de Rham subcomplex; this construction allows to preserve the essence of the 
discrete de Rham structure at the level of the discontinuous spaces, as well 
as commuting projections, both for the primal (strong) and dual (weak) broken-FEEC sequences.
Stabilization is used in stationary formulations to control the interface discontinuities corresponding to the kernels
of the conforming projections, whereas it is not needed for the time-dependent Maxwell system
\cite{campos_pinto_broken-feec_2025,guclu_broken_2023}.
Thus, broken-FEEC is a nonconforming structure-preserving approach in which nonconformity is
controlled by an explicit conforming subcomplex. In contrast to mortar methods
\cite{Buffa_Mortar}, it introduces no independent interface multiplier unknowns.

\begin{figure}[!ht]
    \centering
    \includegraphics[width=0.35\columnwidth]{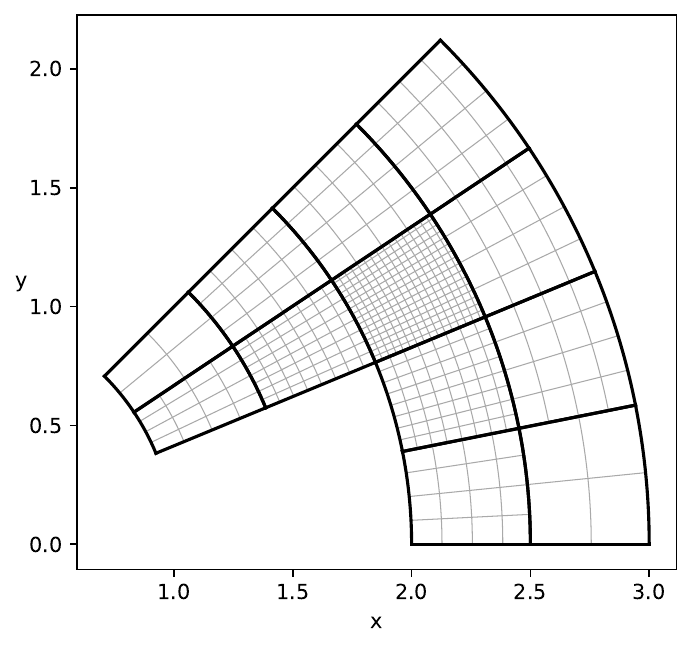}
    \hspace{1cm}
    \includegraphics[width=0.35\columnwidth]{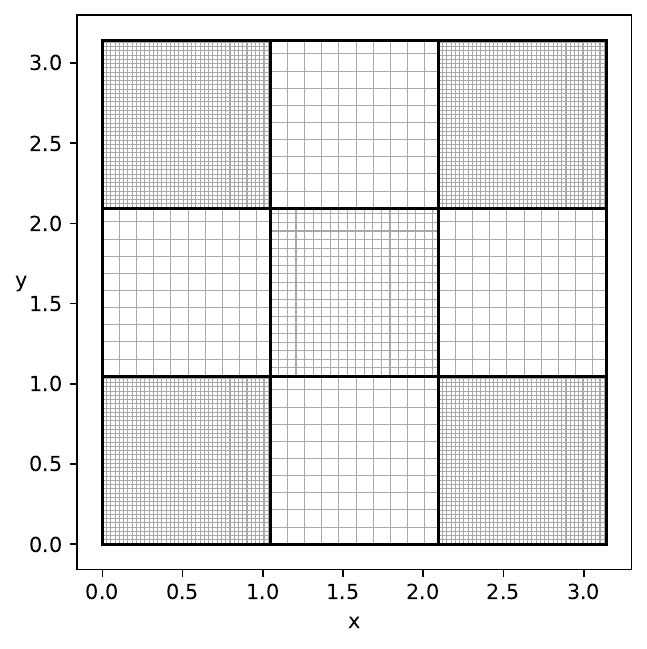}
    \caption{Multipatch geometries with patch-local refinements.}
    \label{fig:intro}
\end{figure}

In this article, we consider the extension of the broken-FEEC constructions to refined patches, as illustrated in Figure \ref{fig:intro}. A common scenario arises when adjacent patches are discretized using spline spaces 
of different resolutions, e.g., different knot vectors, 
resulting in non-matching interfaces as described in \cite{buffa_approximation_2015}. Such configurations allow for flexible handling of nonuniform smoothness and singularities in complex geometries, but they also pose several challenges.
First, preserving the de Rham structure with locally refined splines remains an active area of research. Notable contributions include the construction of discrete de Rham sequences for T-splines and locally refined B-splines \cite{buffa_isogeometric_2014,johannessen_divergence-conforming_2015}, the study of exact sequences of hierarchical spline spaces \cite{evans_hierarchical_2020,shepherd_toshniwal_2024}, and the design of de Rham sequences for splines with varying degrees on domains with polar singularities \cite{toshniwal_isogeometric_2021,patrizi_isogeometric_2021}.
In \cite{BCP}, we constructed local $L^p$-bounded commuting projection operators for
two-dimensional multipatch de Rham sequences with non-matching interfaces, under
additional connectivity and grading assumptions.
A second challenge associated with local refinements is the appearance of numerical artifacts at non-matching patch interfaces. This issue is well known in the context of wave propagation problems where spurious reflections have been observed on finite-difference grids 
with local refinements, and several interesting approaches have been devised to mitigate 
them. In \cite{Frank2004OnSR}, the authors study the good dispersion properties 
of an implicit box scheme, and its ability to avoid spurious modes on nonuniform grids.
In \cite{vay_2004}, a numerical method is proposed that uses anisotropic perfectly matched layers to absorb spurious reflections around non-matching patch interfaces, 
and in \cite{chilton_phd,chilton-collela} the authors devise a fourth-order solver 
for Maxwell's equations on locally refined Cartesian grids, that involves artificial damping and transition layers to mitigate the reflections.

To address these challenges, we develop a structure-preserving alternative to these methods, which
builds on three earlier developments. The abstract analysis of
\cite{campos_pinto_broken-feec_2025} identifies stable moment-preserving conforming
projection operators as an ingredient for high-order broken-FEEC approximations.
The mapped multipatch formulation of \cite{guclu_broken_2023} develops the framework
and its electromagnetic applications for spline discretizations with matching interfaces. The commuting
projections of \cite{BCP} provide a priori approximation and stability results for non-matching
multipatch complexes. These works do not provide the explicit moment-preserving
discrete conforming projections required for non-matching interfaces constructed here.
Our approach is to equip broken-FEEC spaces with discrete conforming projection 
operators that preserve polynomial moments of arbitrary orders.
In particular, we detail the construction of explicit projection matrices 
which enforce full and tangential continuity (corresponding to 
$H^1$ and $H(\curl)$ conformity, respectively) of the discrete fields across 
non-matching interfaces, while preserving integrals against high-order polynomials 
on each logical patch. Our projection matrices are localized close to the patch interfaces
and, like their basic low-order counterparts, they are metric independent.
As a result, our broken-FEEC de Rham sequences have strong and weak (adjoint) 
derivatives which are both local and high-order accurate, leading to
energy-preserving Maxwell solvers which are both explicit and asymptotically free of 
spurious modes. 

The outline is as follows: in Section~\ref{sec:broken-feec} we recall the main lines of the broken-FEEC framework. 
Special attention is devoted to the description of weak differential operators associated with several types of boundary conditions, 
and to the importance of moment preservation in the underlying conforming projection operators. 
The construction of moment-preserving conforming projection operators on non-matching 
patch interfaces is then described in Section~\ref{sec:conf-proj} and Section~\ref{sec:conf-proj-v1}, for $H^1$ and $H(\curl)$ spaces respectively.
Finally, numerical simulations are presented in Section~\ref{sec:num_ex}. These experiments 
assess the properties of our numerical approach on the different discrete models previously described, 
and in particular they allow us to verify that spurious waves virtually disappear 
when moment-preserving projection matrices are used.

\section{Broken-FEEC}
\label{sec:broken-feec}
  
In this article, we consider spline finite element discretizations of the two-dimensional grad-curl de Rham sequence
\begin{equation} \label{dR}
	V^{0} = H^1(\Omega)
	~ \xrightarrow{ \bgrad }  ~
	V^{1} = H(\curl;\Omega) 
	~ \xrightarrow{ \curl}  ~
	V^{2} = L^2(\Omega) 
\end{equation}
defined on an open, bounded domain $\Omega \subset \RR^2$ with Lipschitz boundary,
for which the adjoint sequence \cite{Arnold.Falk.Winther.2006.anum} is
\begin{equation} \label{dR_*}
	L^2(\Omega)
	~ \xleftarrow{ -\Div }  ~
	H_0(\Div;\Omega) 
	~ \xleftarrow{ \bcurl}  ~
	H_0(\bcurl;\Omega). 
\end{equation}
We will also address the sequence with homogeneous boundary conditions
\begin{equation} \label{dR0} V^{0} = H^1_0(\Omega)
	~ \xrightarrow{ \bgrad }  ~
	V^{1} = H_0(\curl;\Omega) 
	~ \xrightarrow{ \curl}  ~
	V^{2} = L^2(\Omega)
\end{equation}
for which the adjoint sequence is the inhomogeneous counterpart of \eqref{dR_*}.

We refer to \cite{Arnold.Falk.Winther.2006.anum,Arnold.Falk.Winther.2010.bams} 
for the analysis of these Hilbert complexes and their structure-preserving
discretizations with finite element exterior calculus (FEEC) spaces.
In this section we describe broken-FEEC spline spaces whose general properties
have been studied in \cite{guclu_broken_2023,campos_pinto_broken-feec_2025},
and we introduce the notation that will allow us to consider multipatch discretizations
with nonuniform resolution.
Note that a strong discretization of \eqref{dR} is naturally associated with 
a weak discretization of the adjoint sequence \eqref{dR_*}.
By using a standard rotation argument, see e.g., \cite[Rem.~2.3.2]{Boffi.Brezzi.Fortin.2013.scm}, one may also apply the constructions presented in this article 
to strong discretizations of the two-dimensional curl-div sequences, namely 
\begin{equation} \label{dR_curl_div}
	H(\bcurl;\Omega)
	~ \xrightarrow{ \bcurl }  ~
	H(\Div;\Omega) 
	~ \xrightarrow{ \Div}  ~
	L^2(\Omega) 
\end{equation}
or its homogeneous counterpart \eqref{dR_*}.

An extension to the three-dimensional de Rham sequence can also be derived by following the same approach,
but it is outside the scope of this article and will be addressed in a forthcoming work.

\begin{remark}[Terminology]
	In this article we follow a recent trend \cite{kraus2016gempic,BCP} to refer to structure-preserving spline finite elements on single or multipatch domains \cite{Buffa_2011} as a FEEC method. Although they differ from the standard (full or trimmed) polynomial elements defined on simplicial complexes \cite{Arnold.Falk.Winther.2006.anum}, they share the fundamental properties that lead to structure-preservation and stability, in that they form discrete de Rham sequences and admit commuting projection operators.
	The term broken-FEEC is then used to refer to the discontinuous discretization of the de Rham sequence on multipatch domains, which is obtained by relaxing the continuity constraints at patch interfaces. 
	This approach was first introduced in \cite{Campos-Pinto.Sonnendrucker.2016.mcomp} in the context of simplicial polynomial FEEC (N\'edel\'ec) spaces, and later extended to multipatch spline spaces in \cite{guclu_broken_2023}.
\end{remark}

\subsection{Multipatch spline finite element spaces}

The discretizations are based on spline finite element spaces associated 
with a multipatch domain, that is a domain 
of the form
\begin{equation} \label{Om}
	\Omega = \Int\Big(\bigcup_{k \in \cK} \overline{\Omega_k} \Big) 
	\quad \text{ with } \quad \Omega_k = F_k(\hat \Omega), 
\end{equation}
where the patches $\Omega_k$ have pairwise disjoint interiors and are geometrically conforming.
More precisely, for any two distinct patches $k_1,k_2 \in\mathcal K$, the
intersection $\overline{\Omega}_{k_1}\cap\overline{\Omega}_{k_2}$ is either empty,
a common vertex, or one complete common patch edge. 
The patches are parametrized by $C^1$-diffeomorphisms $F_k$ 
with patch indices $k \in \cK = \{1, \dots, K\}$,
and a reference Cartesian domain $\hat \Omega = \openint{0}{1}^2$.

Following \cite{Buffa_2011,buffa_isogeometric_2010},
on the mapped subdomains we consider mapped spline spaces
obtained by pushing forward logical tensor-product splines.
For each patch $k$ we then first consider 
two univariate spline spaces on the interval $[0,1]$
forming a one-dimensional de Rham sequence,
\begin{equation} \label{VV01}
	\VV^{0}_{k} \coloneqq \mathbbm{S}^p[\Xi_k]
	~ \xrightarrow{ \partial_{\hx} }  ~
	\VV^{1}_{k} \coloneqq \mathbbm{S}^{p-1} [\Xi_k'].
\end{equation}
These spaces consist of clamped splines of degree $p$ and $p-1$, respectively, 
defined by two knot vectors
$\Xi_k = (\xi^k_i)_{i=0}^{n_k+1+p}$
and $(\Xi_k)' = (\xi^k_i)_{i=1}^{n_k+p}$ with interior-knot multiplicity $1 \le m \le p$,
associated with knots of the form
\begin{equation} \label{knots}
	0 = \xi_0^k = \dots = \xi_p^k < \xi_{p+1}^k = \dots = \xi_{p+m}^k < \dots  
	< \xi_{n_k+1}^k = \dots = \xi_{n_k+1+p}^k = 1.
\end{equation}
Here $n_k = \dim(\VV^1_k) = \dim(\VV^0_k)-1$, and 
we refer to, e.g., \cite{Schumaker_2007} for a detailed description of
univariate spline spaces and B-spline bases.
In this work we will consider that the distinct knots are uniformly spaced in each patch:
this spacing defines the resolution in the patch $k$,
which is allowed to vary from that of the other patches
under a nestedness assumption that will be specified below.

A de Rham sequence on the logical 
domain $\hat \Omega$ is then obtained by considering the tensor-product  
spaces 
\begin{equation*}
\hat V^{0}_{k}  \coloneqq \VV^{0}_k \otimes \VV^{0}_k
~~ \xrightarrow{ \bgrad }  ~~
\hat V^{1}_{k}  \coloneqq \begin{pmatrix}
  \VV^{1}_{k} \otimes \VV^{0}_k
  \\ \noalign{\smallskip}
  \VV^{0}_k\otimes \VV^{1}_{k}
\end{pmatrix}
~~ \xrightarrow{ \curl}  ~~
\hat V^{2}_{k}  \coloneqq \VV^{1}_{k} \otimes \VV^{1}_{k},
\end{equation*}
and on the mapped patch $\Omega_k = F_k(\hat \Omega)$
we consider the spaces
\begin{equation} \label{dR_k}
	V^{0}_{k} \coloneqq \cF^0_k(\hat V^{0}_k)
	~ \xrightarrow{ \bgrad }  ~
	V^{1}_{k} \coloneqq \cF^1_k(\hat V^{1}_k)
	~ \xrightarrow{ \curl}  ~
	V^{2}_{k} \coloneqq \cF^2_k(\hat V^{2}_k)
\end{equation}
associated with the pushforward operators associated with 
$F_k$. In two dimensions, they read
\begin{equation} \label{pf_gc}
  \left\{
  \begin{aligned}
  &\cF^0_k : \hat \phi \mapsto \phi \coloneqq \hat \phi \circ F_k^{-1}
  \\
  &\cF^1_k : \hat \bu \mapsto \bu \coloneqq  \big(DF_k^{-T} \hat \bu \big)\circ F_k^{-1}
  \\
  &\cF^2_k : \hat f \mapsto f \coloneqq  \big(J_{F_k}^{-1} \hat f \big)\circ F_k^{-1}
  \end{aligned}
  \right.
\end{equation}
where $DF_k(\hat \bx) = \big(\partial_b (F_k)_a(\hat \bx)\big)_{1 \le a,b \le 2}$
is the Jacobian matrix of $F_k$, and $J_{F_k}$ 
its metric determinant, 
see e.g., \cite{Hiptmair.2002.anum,kreeft2011mimetic,Buffa_Doelz_Kurz_Schoeps_Vazquez_Wolf_2019,%
Holderied_Possanner_Wang_2021}.

The broken-FEEC discretization involves two sequences of spaces
on the global domain $\Omega$:
a sequence of broken spaces which are fully discontinuous at patch interfaces and therefore smooth only within each patch,
\begin{equation} \label{Vpw}
	V^{\ell}_\pw \coloneqq 
		\{ v \in L^2(\Omega) : v|_{\Omega_k} \in V^{\ell}_k  \text{ for } k \in \cK\}, \quad \ell \in \{0, 1, 2 \}
\end{equation}
and a sequence of conforming subspaces 
\begin{equation} \label{Vh}
	V^{\ell}_h \coloneqq V^{\ell}_\pw \cap V^{\ell}  
		= 
		\{ v \in V^\ell : v|_{\Omega_k} \in V^{\ell}_k  \text{ for } k \in \cK\}.
\end{equation}
As a result of the global de Rham sequence structure 
\eqref{dR} or \eqref{dR0}, and of the local ones 
\eqref{dR_k}, 
which follow from the commutation properties
of the pushforward operators \cite{Hiptmair.2002.anum},
the conforming spaces also form a de Rham sequence,
\begin{equation} \label{dRh}
	V^{0}_{h} 
	~ \xrightarrow{ \bgrad }  ~
	V^{1}_{h} 
	~ \xrightarrow{ \curl}  ~
	V^{2}_{h}.
\end{equation}
Its stability and structure-preserving properties have been analyzed in \cite{BCP},
where $L^p$-stable commuting projection operators were constructed for shape-regular patches under two additional topological and grading hypotheses: every
interior vertex is shared by exactly four patches, and the nested resolutions
of neighboring patches do not form a local checkerboard pattern. These two assumptions are
  specific to the stability analysis;
  neither is required for the
  construction nor the algebraic properties of
  the moment-preserving
  conforming projections developed in the
  present work.
  Under these
assumptions one obtains the commuting diagram

\begin{equation} \label{cd}
	\begin{tikzcd}[row sep=large, column sep=large]
		V^0 \arrow{r}{\bgrad} \arrow{d}{\Pi^0} & V^1 \arrow{r}{\curl} \arrow{d}{\Pi^1} & V^2 \arrow{d}{\Pi^2} \\
		V^0_h \arrow{r}{\bgrad} & V^1_h \arrow{r}{\curl} & V^2_h.
	\end{tikzcd}
\end{equation}

Broken-FEEC schemes rely on {\em conforming projection operators}
\begin{equation} \label{intro:confproj}
    P^\ell: V^{\ell}_\pw \rightarrow V^{\ell}_h \subset V^\ell_\pw
\end{equation}
whose primary purpose is to enforce the relevant continuity conditions at the patch interfaces.
They are regarded as endomorphisms of the broken space
$V^\ell_{\mathrm{pw}}$, with range $V^\ell_h\subset
V^\ell_{\mathrm{pw}}$. Their adjoints are consequently taken with respect to the $L^2$-inner product on the ambient broken space.
These discrete projections endow the broken spaces \eqref{Vpw} with a 
discrete de Rham structure
\begin{equation} \label{dRpw}
	V^{0}_{\pw} 
	~ \xrightarrow{ \bgrad P^0 }  ~
	V^{1}_{\pw} 
	~ \xrightarrow{ \curl P^1 }  ~
	V^{2}_{\pw}
\end{equation}
and an adjoint discrete sequence may be defined by $L^2$-duality,
\begin{equation} \label{dRpw_*}
	V^{0}_{\pw} 
	~ \xleftarrow{ (\bgrad P^0)^* }  ~
	V^{1}_{\pw} 
	~ \xleftarrow{ (\curl P^1)^* }  ~
	V^{2}_{\pw}.
\end{equation}
A key observation is that \eqref{dRpw_*} can only be a high-order discretization of the 
continuous sequence \eqref{dR_*} if the adjoint conforming projection operators
$(P^\ell)^*$ are high-order accurate. That is, if the projections \eqref{intro:confproj}
preserve high-order polynomial moments. 
Their construction on non-matching multipatch spline spaces will be the main goal of this work. 

\subsection{Tensor-product and mapped spline bases} 
\label{tensorprodspaces}

The univariate spaces in \eqref{VV01} admit B-spline bases,
\begin{equation} \label{eq:basis-hatV-1D}
\VV^{0}_{k} = \Span(\{\lambda^{0, k}_{i} : i = 0, \dots, n_k\}),
\quad
\VV^{1}_{k} = \Span(\{\lambda^{1, k}_{i} : i = 0, \dots, n_k-1\}).
\end{equation}
Since we will mostly use the former ones to define our projection operators, we drop the space index for conciseness and write 
\begin{equation}
\lambda^{k}_{i} \coloneqq \lambda^{0, k}_{i}, \qquad 0 \le i \le n_k.
\end{equation}
From the form of the open knot vector \eqref{knots}, it follows that
these functions satisfy an interpolation property at the endpoints, namely
\begin{equation} \label{ipe}
	\lambda^{k}_{i}(0) = \delta_{i,0}
	\qquad \text{ and } \qquad
	\lambda^{k}_{i}(1) = \delta_{i,n_k}.
\end{equation}
For later purposes we assume a minimal resolution for the patches.
\begin{assumption} \label{ass:nbc}
	The spline space dimension satisfies $n_k > p + 1$ for all $k \in \cK$.
\end{assumption}

The bivariate basis functions on the logical domain 
$\hat \Omega$
are then defined as tensor products:
for $\hat V^0_k$ we set
\begin{equation} \label{eq:basis-hatV0}
	\hat \Lambda^{k}_{\bi}(\hbx) \coloneqq \lambda^{k}_{i_1}(\hx_1) \lambda^{k}_{i_2}(\hx_2)
\qquad \text{ for } ~
\bi = (i_1, i_2) \in \bcI^k \coloneqq \{0, \dots, n_k\}^2,
\end{equation}
for the space $\hat V^1_k$ we have 
\begin{equation*}
\left\{\begin{aligned}
	\hat \bLambda^{1, k}_{1, \bi}(\hbx) \coloneqq \begin{pmatrix}
		\lambda^{1, k}_{i_1}(\hx_1) \lambda^{k}_{i_2}(\hx_2) \\ 0
	\end{pmatrix} \qquad \text{ for } ~
	\bi \in \bcI^{1,k}_{1} \coloneqq \{0, \dots, n_k-1\}\times\{0, \dots, n_k\}, \\
	\hat \bLambda^{1, k}_{2, \bi}(\hbx) \coloneqq \begin{pmatrix}
		0 \\ \lambda^{k}_{i_1}(\hx_1) \lambda^{1, k}_{i_2}(\hx_2) 
	\end{pmatrix} \qquad \text{ for } ~
\bi \in \bcI^{1, k}_2 \coloneqq \{0, \dots, n_k\}\times\{0, \dots, n_k-1\},
\end{aligned}\right.
\end{equation*}
and set $\bcI^{1, k} \coloneqq \left\{ (d, \bi) : d \in \{1,2\}, \ \ib \in  \bcI^{1, k}_d  \right\}$.
Finally, for $\hat V^2_k$ we set 
\begin{equation} \label{eq:basis-hatV2}
		\hat \Lambda^{2, k}_{\bi}(\hbx) \coloneqq \lambda^{1, k}_{i_1}(\hx_1) \lambda^{1, k}_{i_2}(\hx_2)
	\qquad \text{ for } ~
	\bi \in \bcI^{2, k} \coloneqq \{0, \dots, n_k-1\}^2.
\end{equation}

Finally, the basis functions for the broken spaces 
\eqref{Vpw} are
obtained by pushing forward the reference basis functions
on the mapped patches $\Omega_k$, $k \in \cK$,
and extending them by zero outside their patch.
The basis of $V^{0}_\pw$ thus reads
\begin{equation} \label{basis_pw}
\Lambda^{k}_{\bi}(\bx) \coloneqq \cF^0_k(\hat \Lambda^{k}_{\bi})(\bx)
	=
	\begin{cases}
		\hat \Lambda^{k}_{\bi}(\hbx^k)  &\text{ for } \bx \in \Omega_k
		\\
		0 &\text{ elsewhere}
	\end{cases}
	\quad \text{ for } k \in \cK, ~ \bi \in \bcI^k ,
\end{equation} 
where $\hbx^k \coloneqq (F_k)^{-1}(\bx)$. 
Similarly, the basis of $V^1_\pw$ reads
\begin{equation} \label{eq:V1-basis}
	\bLambda^{1, k}_{d, \bi}(\bx) \coloneqq \cF^1_k(\hat \bLambda^{1, k}_{d, \bi})(\bx)
		=
		\begin{cases}
			DF_k^{-T} \hat \bLambda^{1,k}_{d,\bi}(\hbx^k)  &\text{for } \bx \in \Omega_k
			\\
			0 &\text{ elsewhere}
		\end{cases}
		\quad \text{for } k \in \cK,  (d,\bi) \in \bcI^{1,k},
\end{equation}
and the basis of $V^2_\pw$ reads
\begin{equation} \label{eq:V2-basis}
\Lambda^{2, k}_{\bi}(\bx) \coloneqq \cF^2_k(\hat \Lambda^{2, k}_{\bi})(\bx)
	=
	\begin{cases}
		J_k^{-1} \hat \Lambda^{2,k}_{\bi}(\hbx^k)  &\text{ for } \bx \in \Omega_k
		\\
		0 &\text{ elsewhere}
	\end{cases}
	\quad \text{ for } k \in \cK, \bi \in \bcI^{2,k}.
\end{equation}

\subsection{Conforming constraints on patch interfaces} \label{sec:conf-constr}

Since the local spaces $V^\ell_k$ consist of continuous functions,
the conforming subspaces $V^\ell_h \subset V^\ell_\pw$ may be
characterized by continuity constraints on the patch interfaces,
namely on the intersections of patch boundaries
$\partial \Omega_k \cap \partial \Omega_l$, $k \neq l$.
Because we assume that patches are geometrically conforming,
these interfaces are either an edge or a vertex of both patches,
and it follows from the Lipschitz assumption on $\partial \Omega$ 
that if a vertex is shared by two patches, then it belongs to an edge also 
shared by two patches (an interior edge).
In particular, the continuity constraints characterizing 
the conforming spaces \eqref{Vh} may be expressed on interior edges.

For the first space $V^0_h = V^0_\pw \cap V^0$, with $V^0 = H^1(\Omega)$ 
as in \eqref{dR}, the continuity condition reads:
$\phi \in V^{0}_\pw$ belongs to $V^0$ if and only if
\begin{equation} \label{V0_confcond}
	\phi|_{\Omega_{k^-}} = \phi|_{\Omega_{k^+}}
	\qquad \text{ on every interior edge } \edge = \partial \Omega_{k^-} \cap \partial \Omega_{k^+}.
\end{equation}
In the case of homogeneous boundary conditions \eqref{dR0} where $V^0 = H^1_0(\Omega)$,
an additional condition is imposed on the boundary edges, namely
\begin{equation} \label{V0_confcond_hom}
	\phi|_{\edge} = 0
	\qquad \text{ on every boundary edge } \edge \subset \partial \Omega.
\end{equation}
Similarly, a function $\bu \in V^{1}_\pw$ belongs to $V^1 = H(\curl;\Omega)$ if and only if
\begin{equation} \label{V1_confcond}
	\btau \cdot \bu|_{\Omega_{k^-}} = \btau \cdot \bu|_{\Omega_{k^+}} 
	\qquad \text{ on every interior edge } \edge = \partial \Omega_{k^-} \cap \partial \Omega_{k^+}
\end{equation}
where $\btau = \btau(\bx)$ denotes an arbitrary vector tangent to the edge. For the sequence \eqref{dR}
this characterizes the vector fields in $V^1_h$.
In the case of homogeneous boundary conditions \eqref{dR0} where $V^1 = H_0(\curl;\Omega)$,
an additional condition is imposed on the boundary edges, namely
\begin{equation} \label{V1_confcond_hom}
	\btau \cdot \bu|_{\edge} = 0
	\qquad \text{ on every boundary edge } \edge \subset \partial \Omega.
\end{equation}
For the last space of the sequence there are no constraints: since $V^2 = L^2(\Omega)$,
the conforming space is simply $V^2_h = V^2_\pw$.

For later purposes we denote by $\edges$ and $\vertices$ 
the respective sets of patch edges and vertices. 
Given $\edge \in \edges$, we gather the indices of adjacent patches in
\begin{equation*}
\cK(\edge) = \{k \in \cK: \edge \subset \partial \Omega_k\},
\end{equation*}
and similarly for a given $\vertex \in \vertices$, we gather the 
indices of the adjacent patches in the set
\begin{equation*}
\cK(\vertex) = \{k \in \cK: \vertex \in \partial \Omega_k\}. 
\end{equation*}
We also define patch neighborhoods of edges and vertices as
\begin{equation} \label{Omega_ev}
	\Omega(\edge) \coloneqq \Int\big(\cup_{k \in \cK(\edge)} \overline{\Omega_k}\big),
	\qquad
	\Omega(\vertex) \coloneqq \Int\big(\cup_{k \in \cK(\vertex)} \overline{\Omega_k}\big).
\end{equation}
Clearly, $\cK(\edge)$ contains two patches if $\edge$ is an interior edge, 
and only one patch if it is a boundary edge.
We make the following assumption on shared interfaces:
\begin{assumption} \label{as:e_nested}
	Interior edges $\edge$ are shared by two patches with nested resolution.
	That is, for $\cK(\edge) = \{k^-(\edge), k^+(\edge)\}$, the knot vector $\Xi_+$ of the fine patch $k^+(\edge)$ results from a uniform refinement of the knot vector $\Xi_-$ of the coarse patch $k^-(\edge)$. More precisely, every interval between two consecutive distinct knots of $\Xi_-$ is
  subdivided into the same number of equal subintervals by the distinct
  knots of $\Xi_+$.
	In addition, the spline degree $p$ and the interior-knot multiplicity $m$ are constant
  	across all patches.
\end{assumption}

\begin{remark}
	In practice, this assumption can be lifted by just requiring  $\VV^{0}_{k^-(\edge)} \subset \VV^{0}_{k^+(\edge)}$ and allowing for different degrees across patches. The analysis of the resulting projection operators would be more involved, but the main ideas would remain the same.
\end{remark}

\subsection{Weak differential operators in conforming FEEC sequences} 
\label{sec:wdiff_conf}

In this section and the following one, we specify in which sense discrete adjoint sequences such as 
\eqref{dRpw_*} approximate the continuous ones. Our main observation will be that 
broken-FEEC spaces allow more flexibility than conforming ones 
in the choice of the underlying weak boundary conditions. We will also
specify a notion of moment preservation that will allow us to reach 
high-order accuracy for the weak differential operators.

Let us begin by reminding that if one considers the grad-curl sequence 
{\em without} boundary conditions \eqref{dR} 
as a primal de Rham sequence, then the associated adjoint sequence, defined by $L^2$-duality,
is the curl-div sequence {\em with} homogeneous boundary conditions \eqref{dR_*}.
Conversely, if the grad-curl sequence with boundary conditions \eqref{dR0} 
is chosen as a primal de Rham sequence, then the corresponding adjoint sequence 
is the curl-div sequence without boundary conditions.

At the discrete level a similar behavior can be observed. 
For instance, in the case of conforming spaces $V^\ell_h$ approximating the inhomogeneous sequence \eqref{dR}, one may
expect that the discrete adjoint operators approximate 
the continuous differential operators with homogeneous boundary conditions. 
This can be verified by considering the discrete adjoint operators
$\widetilde{\Div}_{h,0} := -(\bgrad)^* : V^1_{h} \rightarrow V^0_{h}$
and $\widetilde{\bcurl}_{h,0} := (\curl)^* : V^2_{h} \rightarrow V^1_{h}$
defined by the relations 
\begin{equation} \label{div_curl_h0}
	\left\{\begin{aligned}	
	\sprod{\phi_h}{\widetilde{\Div}_{h,0} \bu_h} &= - \sprod{\bgrad \phi_h}{\bu_h}, 
	\quad &&\forall \bu_h \in V^1_{h}, ~ \phi_h \in V^0_{h},
	\\
	\sprod{\bv_h}{\widetilde{\bcurl}_{h,0} w_h} &= \sprod{\curl \bv_h}{w_h}, 
	\quad &&\forall w_h \in V^2_{h}, ~ \bv_h \in V^1_{h},
	\end{aligned}\right.
\end{equation}
and letting $\tilde \Pi_h^\ell$ be the $L^2$-projection onto $V^\ell_h$.
Here, we have used a tilde to indicate a weak discretization (involving test functions) and a subscript $0$ to indicate that the discrete operators approximate the continuous ones on spaces with homogeneous boundary conditions.
Indeed the following relations hold:
\begin{equation} \label{dual-comm_0}
	\left\{\begin{aligned}	
	\widetilde{\Div}_{h, 0} \tilde \Pi_{h}^1 \bu &= \tilde \Pi_{h}^0 \Div \bu,
	\quad &&\forall \bu \in H_0(\Div;\Omega),
	\\
	\widetilde{\bcurl}_{h, 0} \tilde \Pi_{h}^2 w &= \tilde \Pi_{h}^1 \bcurl w,
	\quad &&\forall w \in H_0(\bcurl;\Omega).
	\end{aligned}\right.
\end{equation}
To verify the first one, we observe that for all 
$\phi_h \in V^0_h$ and $\bu \in H_0(\Div;\Omega)$ we have 
\begin{align*}
	\sprod{\phi_h}{\widetilde{\Div}_{h, 0} \tilde \Pi_{h}^1 \bu}
	&= -\sprod{\bgrad \phi_h}{\tilde \Pi_{h}^1 \bu} \\
	&= -\sprod{\bgrad \phi_h}{\bu}
	= \sprod{\phi_h}{\Div \bu}
	= \sprod{\phi_h}{\tilde \Pi^0_h \Div \bu},	
\end{align*}
and the relation involving $\widetilde{\bcurl}_{h, 0}$ is verified in the same way.

Conversely, if the conforming spaces (denoted here $V^\ell_{h,0}$ for clarity)
approximate the sequence \eqref{dR0} {\em with} homogeneous boundary conditions, 
then the discrete adjoint operators approximate the continuous differential operators {\em without} boundary conditions. 
To specify this, consider the operators
$\widetilde{\Div}_{h} := -(\bgrad)^* : V^1_{h,0} \rightarrow V^0_{h,0}$
and $\widetilde{\bcurl}_{h} := (\curl)^* : V^2_{h,0} \rightarrow V^1_{h,0}$
defined by the relations 
\begin{equation} \label{div_curl_h}
	\left\{\begin{aligned}	
	\sprod{\phi_h}{\widetilde{\Div}_{h} \bu_h} &= - \sprod{\bgrad \phi_h}{\bu_h}, 
	\quad &&\forall \bu_h \in V^1_{h,0}, ~ \phi_h \in V^0_{h,0},
	\\
	\sprod{\bv_h}{\widetilde{\bcurl}_{h} w_h} &= \sprod{\curl \bv_h}{w_h}, 
	\quad &&\forall w_h \in V^2_{h,0}, ~ \bv_h \in V^1_{h,0},
	\end{aligned}\right.
\end{equation}
and let $\tilde \Pi^\ell_{h,0}$ be the $L^2$-projection onto $V^\ell_{h,0}$. 
For these weak discrete operators, the absence of the subscript $0$ in the notation indicates that they approximate the continuous ones on spaces without boundary conditions. Indeed the following relations hold:
\begin{equation} \label{dual-comm}
	\left\{\begin{aligned}	
	\widetilde{\Div}_{h} \tilde \Pi_{h, 0}^1 \bu &= \tilde \Pi_{h, 0}^0 \Div \bu,
	\quad &&\forall \bu \in H(\Div;\Omega),
	\\
	\widetilde{\bcurl}_{h} \tilde \Pi_{h, 0}^2 w &= \tilde \Pi_{h, 0}^1 \bcurl w,
	\quad &&\forall w \in H(\bcurl;\Omega).
	\end{aligned}\right.
\end{equation}
The same argument applies: for all 
$\phi_h \in V^0_{h,0}$ and $\bu \in H(\Div;\Omega)$, compute
\begin{align*}
	\sprod{\phi_h}{\widetilde{\Div}_{h} \tilde \Pi_{h,0}^1 \bu}
	&= -\sprod{\bgrad \phi_h}{\tilde \Pi_{h,0}^1 \bu} \\
	&= -\sprod{\bgrad \phi_h}{\bu}
	= \sprod{\phi_h}{\Div \bu}
	= \sprod{\phi_h}{\tilde \Pi^0_{h,0} \Div \bu},	
\end{align*}
and the relation involving $\widetilde{\bcurl}_{h}$ is verified in the same way.
We may summarize these observations as follows.

\begin{theorem} \label{thm:weak-diff_conf}
	Let $\tilde \Pi^\ell_h$ be the $L^2$-projection onto 
	the conforming space $V^\ell_h$ without boundary conditions, and let 
	\begin{equation}
		\widetilde{\Div}_{h,0} := - (\bgrad)^* : V^1_h \rightarrow V^0_h, 
		\qquad 
		\widetilde{\bcurl}_{h, 0} := (\curl)^*: V^2_h \rightarrow V^1_h
	\end{equation}
	be the discrete differential operators corresponding to \eqref{div_curl_h0}.
	Then in the following diagram the operators denoted by solid arrows commute:
	\begin{equation} \label{cd-dual}
		\begin{tikzcd}[row sep=large, column sep=large]
			V^0_{h} \arrow[shift left, dashed]{r}{\bgrad} & V^1_{h} \arrow[shift left, dashed]{r}{\curl} \arrow[shift left]{l}{\widetilde{\Div}_{h, 0}} & \arrow[shift left]{l}{\widetilde{\bcurl}_{h, 0}} V^2_{h} \\
			\arrow{u}{\tilde\Pi^0_{h}} L^2 &  \arrow{l}{\Div} \arrow{u}{\tilde\Pi^1_{h}} H_0(\Div) & \arrow{l}{\bcurl} \arrow{u}{\tilde\Pi^2_{h}} H_0(\bcurl)
		\end{tikzcd}
	\end{equation}
	(where the dashed arrows represent the strong differential operators defining
	the discrete adjoints).
	Conversely, let $\tilde \Pi^\ell_{h,0}$ be the $L^2$-projection onto 
	the conforming space with homogeneous boundary conditions $V^\ell_{h,0}$, 
	and let 
	\begin{equation}
		\widetilde{\Div}_{h} := - (\bgrad)^* : V^1_{h,0} \rightarrow V^0_{h,0}, 
		\qquad 
		\widetilde{\bcurl}_{h} := (\curl)^*: V^2_{h,0} \rightarrow V^1_{h,0}
	\end{equation}
	be the discrete differential operators corresponding to 
	\eqref{div_curl_h}.
	Then in the following diagram the operators denoted by solid arrows commute:
	\begin{equation} \label{cd-dual-0}
	\begin{tikzcd}[row sep=large, column sep=large]
		V^0_{h,0} \arrow[shift left, dashed]{r}{\bgrad} & V^1_{h,0} \arrow[shift left, dashed]{r}{\curl} \arrow[shift left]{l}{\widetilde{\Div}_h} & \arrow[shift left]{l}{\widetilde{\bcurl}_h} V^2_{h,0} \\
		\arrow{u}{\tilde\Pi^0_{h,0}} L^2 &  \arrow{l}{\Div} \arrow{u}{\tilde\Pi^1_{h,0}} H(\Div; \Omega) & \arrow{l}{\bcurl} \arrow{u}{\tilde\Pi^2_{h,0}} H(\bcurl; \Omega)
	\end{tikzcd}
	\end{equation}
	(where again the dashed arrows represent the strong differential operators defining
	the discrete adjoints).
	In particular, the boundary conditions associated with the weak
	differential operators are determined by the underlying conforming finite element spaces.
\end{theorem}

\subsection{Weak broken differential operators and moment preservation} 
\label{sec:mom-pres}

In the broken-FEEC framework, boundary conditions are not imposed 
on the finite element spaces themselves 
(which are defined patch-wise without taking into account 
neighboring patches or global constraints) but through the conforming projection operators.
As a result, one can expect that the above duality principle
applies in a flexible way, allowing one to choose the weak boundary
conditions by selecting a suitable conforming projection.
To verify this flexible duality principle, let us assume that 
\begin{equation} \label{Pell_Pell0}
	P^\ell: V^\ell_\pw \rightarrow V^\ell_{h} \subset V^\ell_\pw
	\quad \text{ and } \quad
	P^\ell_0: V^\ell_\pw \rightarrow V^\ell_{h,0} \subset V^\ell_\pw
\end{equation}
are conforming projection operators onto the discrete spaces without and with homogeneous boundary conditions, respectively.
Examples of such projections will be described later in this article.
We may then consider two kinds of weak differential operators, as above.
First, the adjoints of the strong broken-FEEC differential operators associated with {\em no} boundary conditions, 
which should in principle approximate the adjoint differential operators {\em with} homogeneous boundary conditions:
$\widetilde{\Div}_{\pw,0} := -(\bgrad P^0)^* : V^1_\pw \rightarrow V^0_\pw$
and $\widetilde{\bcurl}_{\pw,0} := (\curl P^1)^* : V^2_\pw \rightarrow V^1_\pw$,
defined by the relations 
\begin{equation} \label{div_curl_pw0}
	\left\{\begin{aligned}	
	\sprod{\phi_h}{\widetilde{\Div}_{\pw,0} \bu_h} &= - \sprod{\bgrad P^0 \phi_h}{\bu_h}, 
	\quad &&\forall \bu_h \in V^1_\pw, ~ \phi_h \in V^0_\pw,
	\\
	\sprod{\bv_h}{\widetilde{\bcurl}_{\pw,0} w_h} &= \sprod{\curl P^1 \bv_h}{w_h}, 
	\quad &&\forall w_h \in V^2_{\pw}, ~ \bv_h \in V^1_{\pw}.
	\end{aligned}\right.
\end{equation}
Second, the adjoints of the broken-FEEC differential operators associated {\em with} homogeneous boundary conditions, 
which should approximate the adjoint differential operators {\em without} boundary conditions:
$\widetilde{\Div}_{\pw} := -(\bgrad P^0_0)^* : V^1_\pw \rightarrow V^0_\pw$
and $\widetilde{\bcurl}_{\pw} := (\curl P^1_0)^* : V^2_\pw \rightarrow V^1_\pw$,
defined by the relations 
\begin{equation} \label{div_curl_pw}
	\left\{\begin{aligned}	
	\sprod{\phi_h}{\widetilde{\Div}_{\pw} \bu_h} &= - \sprod{\bgrad P^0_0 \phi_h}{\bu_h}, 
	\quad &&\forall \bu_h \in V^1_\pw, ~ \phi_h \in V^0_\pw,
	\\
	\sprod{\bv_h}{\widetilde{\bcurl}_{\pw} w_h} &= \sprod{\curl P^1_0 \bv_h}{w_h}, 
	\quad &&\forall w_h \in V^2_{\pw}, ~ \bv_h \in V^1_{\pw}.
	\end{aligned}\right.
\end{equation}
We finally consider two kinds of patch-wise approximation operators, namely
\begin{equation} \label{tPipw}
	\tilde \Pi^\ell_\pw := (P^\ell)^* Q^\ell_\pw : L^2(\Omega) \to V^\ell_\pw
	\quad \text{ and } \quad
	\tilde \Pi^\ell_{\pw,0} := (P^\ell_0)^* Q^\ell_\pw : L^2(\Omega) \to V^\ell_\pw	
\end{equation}
where $Q^\ell_\pw$ is the $L^2$-projection onto $V^\ell_\pw$.
Here we have used again a tilde to indicate a weak discretization (involving test functions) and a subscript $0$ to indicate which discrete operators approximate the continuous ones on spaces with homogeneous boundary conditions.
Indeed, the following relations hold:
\begin{equation} \label{dual-comm_pw0}
	\left\{\begin{aligned}	
	\widetilde{\Div}_{\pw, 0} \tilde \Pi_{\pw}^1 \bu &= \tilde \Pi_{\pw}^0 \Div \bu,
	\quad &&\forall \bu \in H_0(\Div;\Omega),
	\\
	\widetilde{\bcurl}_{\pw, 0} \tilde \Pi_{\pw}^2 w &= \tilde \Pi_{\pw}^1 \bcurl w,
	\quad &&\forall w \in H_0(\bcurl;\Omega),
	\end{aligned}\right.
\end{equation}
and 
\begin{equation} \label{dual-comm_pw}
	\left\{\begin{aligned}	
	\widetilde{\Div}_{\pw} \tilde \Pi_{\pw,0}^1 \bu &= \tilde \Pi_{\pw,0}^0 \Div \bu,
	\quad &&\forall \bu \in H(\Div;\Omega),
	\\
	\widetilde{\bcurl}_{\pw} \tilde \Pi_{\pw,0}^2 w &= \tilde \Pi_{\pw,0}^1 \bcurl w,
	\quad &&\forall w \in H(\bcurl;\Omega).
	\end{aligned}\right.
\end{equation}
These relations are easily verified by applying the definitions and integrating by parts as above: for instance, for the first relation in \eqref{dual-comm_pw0}, we have for all $\phi_h \in V^0_\pw$ and $\bu \in H_0(\Div;\Omega)$
\begin{align*}
	\sprod{\phi_h}{\widetilde{\Div}_{\pw, 0} \tilde \Pi_{\pw}^1 \bu}
	&= -\sprod{\bgrad P^0 \phi_h}{\tilde \Pi_{\pw}^1 \bu} \\
	&= -\sprod{P^1 \bgrad P^0 \phi_h}{Q^1_\pw \bu} \\
	&= -\sprod{\bgrad P^0 \phi_h}{Q^1_\pw \bu} \\
	&= -\sprod{\bgrad P^0 \phi_h}{\bu} \\
	&= \sprod{P^0 \phi_h}{\Div \bu}
	= \sprod{P^0 \phi_h}{Q^0_{\pw} \Div \bu}
	= \sprod{\phi_h}{\tilde \Pi^0_{\pw} \Div \bu}.
\end{align*}
The other relations in \eqref{dual-comm_pw0} and \eqref{dual-comm_pw} are proven in the same way.

We summarize the preceding results in a formal statement.

\begin{theorem} \label{thm:weak-diff_pw}
	Let $P^\ell: V^\ell_\pw \to V^\ell_h$ 
	and $P^\ell_0: V^\ell_\pw \to V^\ell_{h,0}$ be discrete projection operators onto 
	the conforming spaces without and with homogeneous boundary conditions, respectively.
	Then the discrete differential operators corresponding to \eqref{div_curl_pw0}, i.e.
	\begin{equation}
		\widetilde{\Div}_{\pw,0} := - (\bgrad P^0)^* : V^1_\pw \rightarrow V^0_\pw, 
		\qquad 
		\widetilde{\bcurl}_{\pw, 0} := (\curl P^1)^*: V^2_\pw \rightarrow V^1_\pw,
	\end{equation}
	yield a commuting diagram (for the solid arrows as above, the dashed arrows representing the strong differential operators defining the discrete adjoints):
	\begin{equation} \label{cd-dual_pw}
		\begin{tikzcd}[row sep=large, column sep=large]
			V^0_{\pw} \arrow[shift left, dashed]{r}{\bgrad P^0} & V^1_{\pw} \arrow[shift left, dashed]{r}{\curl P^1} \arrow[shift left]{l}{\widetilde{\Div}_{\pw, 0}} & \arrow[shift left]{l}{\widetilde{\bcurl}_{\pw, 0}} V^2_{\pw} \\
			\arrow{u}{\tilde\Pi^0_{\pw}} L^2 &  \arrow{l}{\Div} \arrow{u}{\tilde\Pi^1_{\pw}} H_0(\Div) & \arrow{l}{\bcurl} \arrow{u}{\tilde\Pi^2_{\pw}} H_0(\bcurl)
		\end{tikzcd}
	\end{equation}
	where $\tilde \Pi^\ell_{\pw} := (P^\ell)^* Q^\ell_\pw$ is defined as in \eqref{tPipw}.
	Similarly, the discrete differential operators corresponding to \eqref{div_curl_pw}, i.e.
	\begin{equation}
		\widetilde{\Div}_{\pw} := - (\bgrad P^0_0)^* : V^1_\pw \rightarrow V^0_\pw, 
		\qquad 
		\widetilde{\bcurl}_{\pw} := (\curl P^1_0)^*: V^2_\pw \rightarrow V^1_\pw,
	\end{equation}
	yield a commuting diagram (for the solid arrows):
	\begin{equation} \label{cd-dual_pw0}
		\begin{tikzcd}[row sep=large, column sep=large]
			V^0_{\pw} \arrow[shift left, dashed]{r}{\bgrad P^0_0} & V^1_{\pw} \arrow[shift left, dashed]{r}{\curl P^1_0} \arrow[shift left]{l}{\widetilde{\Div}_{\pw}} & \arrow[shift left]{l}{\widetilde{\bcurl}_{\pw}} V^2_{\pw} \\
			\arrow{u}{\tilde\Pi^0_{\pw,0}} L^2 &  \arrow{l}{\Div} \arrow{u}{\tilde\Pi^1_{\pw,0}} H(\Div) & \arrow{l}{\bcurl} \arrow{u}{\tilde\Pi^2_{\pw,0}} H(\bcurl)
		\end{tikzcd}
	\end{equation}
	where $\tilde \Pi^\ell_{\pw,0} := (P^\ell_0)^* Q^\ell_\pw$ is defined as in \eqref{tPipw}.
\end{theorem}

\begin{remark}
	Let us emphasize that the main difference with the conforming case summarized by Theorem~\ref{thm:weak-diff_conf} is that boundary conditions in the broken-FEEC framework are not imposed by the finite element spaces themselves but through the conforming projection operators. Since projecting onto a space with homogeneous boundary conditions simply amounts to setting boundary degrees of freedom to zero, this framework leads to a greater flexibility in the selection 
	of the weak boundary conditions.
\end{remark}

The observations above also suggest that the approximation properties of the weak differential operators depend on the approximation properties of the 
operator $\tilde \Pi^\ell_\pw = (P^\ell)^* Q^\ell_\pw: L^2(\Omega) \to V^\ell_\pw$. 
This operator is a projection onto a subspace of $V^\ell_\pw$: namely, onto $(P^\ell)^* V^\ell_\pw$.
For this projection to be a high-order approximation operator, we see that $(P^\ell)^*$ 
must preserve the high-order polynomials in $V^\ell_\pw$.
In other terms, $P^\ell$ must preserve polynomial moments (in logical coordinates).
\begin{definition}\label{mom_pres} 
An operator $P: V^\ell_\pw \rightarrow V^\ell_\pw$ is called moment-preserving of order 
$r \ge 0$
if for any $u \in V^\ell_\pw$ and any patch $k \in \cK$, the pullbacks 
\begin{equation*}
\hat v_k = (\cF^\ell_k)^{-1}((Pu) |_{\Omega_k}) \quad \text{ and } \quad 
\hat u_k = (\cF^\ell_k)^{-1}(u |_{\Omega_k})
\end{equation*}
satisfy the equality
\begin{equation} \label{eq:moment-preservation}
	\sprod{\hat v_k}{q} 
	= \sprod{\hat u_k}{q} 
\end{equation}
for any scalar (or vector-valued) polynomial $q \in \polspace$ defined on the logical domain $\hat \Omega$. 
\end{definition}

\begin{remark}
	To check or enforce \eqref{eq:moment-preservation}, one may choose a basis of $\polspace$ consisting of polynomials of the form $q(\hat \bx) = q_1(\hat x_1) q_2(\hat x_2)$, where $q_1$ and $q_2$ are univariate polynomials of degree at most $r-1$. This decomposes \eqref{eq:moment-preservation} into univariate conditions. 
\end{remark}

\begin{remark}
	In the case $r=0$, we have $\QQ_{-1, -1}[\hat \bx] = \{ 0 \}$ and the preceding definition is empty.
	By convention, this case denotes the absence of an imposed moment-
preservation condition.
\end{remark}
\begin{remark} \label{rem:mom-pres}
	The $L^2$-projection operator $Q^\ell_\pw: L^2(\Omega) \rightarrow V^\ell_\pw$ defined by
	\begin{equation*}
		\langle Q^\ell_\pw v , u\rangle = \langle v, u \rangle,  \quad \forall v \in L^2(\Omega), \quad \forall u \in V^\ell_\pw,
	\end{equation*}
	is moment-preserving of order $p+1$ for $\ell = 0$, and of order $p$ for $\ell \ge 1$, where $p$ is the polynomial degree of the discrete space $V^0$.
	Accordingly, the form of $\tilde \Pi^\ell_\pw$ above motivates conforming projection operators $P^\ell$ that preserve moments of order $r = r(P^\ell)$ with $r(P^0) \le p+1$ and 
	$r(P^1) \le p$.
\end{remark}

  \section{Conforming projection operators on \texorpdfstring{$V^0_\pw$}{V0 patch-wise}}
  \label{sec:conf-proj}
 On the first space of \eqref{dR}, where the $H^1$ conformity amounts to continuity of functions across all patches
as described in Section~\ref{sec:conf-constr}, 
the discrete conforming space amounts to $V^0_h = V^0_\pw \cap H^1(\Omega)$.
We construct a conforming projection 
\begin{equation}
    P^0: V^0_\pw \to V^0_\pw    
\end{equation}
onto the space $V^0_h$
by averaging degrees of freedom of interface basis functions.
In our approach, we first ensure continuity at domain vertex basis functions
by applying a vertex-based conforming projection $P_\vertices$.
Subsequently, we apply an edge-based conforming projection $P_\edges$ 
for continuity along edge-interior basis functions, while preserving continuity at vertices.
The first operation is straightforward, but the second one involves 
extension and restriction operators, mapping to the shared interface-spaces, introduced in Section~\ref{sec:ext-res}, 
since the patch resolutions may be different across the edges.
In addition to the averaging procedure that enforces the continuity at the patch
interfaces, we complete our
projection operators with local correction terms in order to preserve polynomial moments.
A benefit of this approach is to allow for a local and fully decoupled treatment of vertices and edges.
For the clarity of the presentation, 
we first list some properties of the vertex- and edge-based projection operators that will be discussed in the following Section~\ref{sec:p_vertex} and  Section~\ref{sec:p_edge}.
The following results are stated under Assumptions~\ref{ass:nbc} and
\ref{as:e_nested}.
 
\begin{prop} \label{prop:vertex}
    Let $\phi \in V^0_\pw$, $\vertex \in \vertices$, 
    and the order of moment preservation be $r = r(P^0)$ with $0 \le r \le p + 1$.
    Then, the local vertex-based operator 
    $P_\vertex: V^0_\pw \to V^0_\pw$ from Definition~\ref{eq:Pv} below satisfies the following properties
    \begin{enumerate}
        \item $P_\vertex \phi$ is continuous at $\vertex$.
        \item $P_\vertex \phi = \phi$ if and only if $\phi$ is continuous at $\vertex$.
        \item $P_\vertex$ preserves polynomial moments of order $r$ in the sense of Definition~\ref{mom_pres}.
        \item $P_\vertex P_{\vertex'} \phi = P_{\vertex'} P_\vertex \phi $ for any $\vertex' \in \vertices$.
    \end{enumerate}
\end{prop}

\begin{prop}\label{prop:edge}
    Let $\phi \in V^0_\pw$, such that it is continuous on vertices, $\edge \in \edges$, and
    the order of moment preservation be $r = r(P^0)$ with $0 \le r \le p + 1$.
    Then, the local edge-based operator $P_\edge: V^0_\pw \to V^0_\pw$ from Definition~\ref{eq:Pe} below
    satisfies the following properties
    \begin{enumerate}
        \item $P_\edge \phi$ is continuous on $\edge$.
        \item $P_\edge \phi = \phi$ if and only if $\phi$ is continuous on $\edge$.
        \item $P_\edge$ preserves polynomial moments of order $r$ in the sense of Definition~\ref{mom_pres}.
        \item $P_\edge P_{\edge'} \phi = P_{\edge'} P_\edge \phi $ for any $\edge' \in \edges$.
    \end{enumerate}
\end{prop}

 Combining both operators yields the full $V^0$ conforming projection:

\begin{theorem}
    Denote the vertex- and edge-based operators by
    \begin{equation} \label{eq:PvPe}
        P_\vertices = \prod_{\vertex \in \vertices} P_\vertex, \quad P_\edges = \prod_{\edge \in \edges} P_\edge,
    \end{equation}
consisting of the local operators from Definition~\ref{eq:Pv} and Definition~\ref{eq:Pe}, which respectively satisfy Proposition~\ref{prop:vertex} and Proposition~\ref{prop:edge}.
Then, the conforming projection operator
\begin{equation} \label{eq:P0}
    P^0 = P_\edges P_\vertices: V^0_\pw \rightarrow V^0_\pw
\end{equation}
is a projection onto the conforming space $V^0_h$ which is moment-preserving of order $0 \le r \le p + 1$.
\end{theorem}

\begin{proof} 
    Note that the products in \eqref{eq:PvPe} are well-defined, as their order does not matter 
    thanks to $(4)$ in Proposition~\ref{prop:vertex} and Proposition~\ref{prop:edge}.
    For the range property of $P^0$, let $\phi \in V^0_\pw$ and write
    $\tilde \phi  = P_\vertices \phi$, which is continuous at any vertex 
    by $(1)$ in Proposition~\ref{prop:vertex}. 
    This continuity is preserved by $P_\edges$,  and $P_\edges \tilde \phi$ is additionally
    continuous on any edge by $(1)$ in Proposition~\ref{prop:edge}. Thus, $P^0$ maps into the continuous space $V^0_h$.
    The projection property of $P^0$ directly follows from $(2)$ in Proposition~\ref{prop:vertex} and Proposition~\ref{prop:edge}.
By $(3)$ in Proposition~\ref{prop:vertex} and Proposition~\ref{prop:edge}, 
since every operator in the composition preserves polynomial moments,
the product preserves polynomial moments. 
\end{proof}

  \subsection{Vertex-based projections}
  \label{sec:p_vertex}
  In this section, we detail the construction of the vertex-conforming projection
of Proposition~\ref{prop:vertex} resulting in Definition~\ref{eq:Pv}. Let the order of moment preservation be given by
 $0 \le r \le p+1$ 
and fix a vertex $\vertex \in \vertices$. For ease of notation,
assume that the ordering of basis functions of $V^0_k$, cf. \eqref{basis_pw}, in every adjacent patch $k \in \cK(\vertex)$
starts with the one basis function on the vertex $\vertex$, i.e. $\Lambda^k_{\bzero}$ is associated to $\vertex$ for every $k\in \cK(\vertex)$.

\subsubsection*{Local vertex-conforming projection} 
The local vertex-conforming projection $P_\vertex$ enforces continuity at the vertex $\vertex$.
In absence of moment preservation ($r = 0$), it maps the vertex basis function from each
adjacent patch to their common average:
\begin{equation} \label{eq:Pv_r0}
    P_\vertex = \Lambda^k_\bj \mapsto \left\{ \begin{aligned}
                            &\bar\Lambda^\vertex &&\text{ if } \bj = \bzero \text{ and } k \in \cK(\vertex), \\
                            &\Lambda^k_\bj &&\text{ else,}
    \end{aligned}\right.
\end{equation}
where $k\in \cK$ and the vertex-conforming function $\bar\Lambda^\vertex$ is defined as the average of the patch-wise basis functions at the vertex: 
\begin{equation} \label{eq:vertex-basis}
    \bar\Lambda^\vertex = \frac{1}{\#\cK(\vertex)} \sum_{k' \in \mathcal{K}(\vertex)} \Lambda^{k'}_{\bzero},
\end{equation}
where $\#\cK(\vertex)$ denotes the number of patches incident to $\vertex$.

\subsubsection*{Preservation of polynomial moments} 
\label{subsec:vertex-pol-mom}
For $r \ge 1$, we need to add correction terms to the definition of $P_\vertex$ in \eqref{eq:Pv_r0}.
Following Definition~\ref{mom_pres}, we choose
correction terms that involve additional interior degrees of
freedom, and we use those that are close to the
vertex and that preserve the continuity.
So for any patch $k \in \cK(\vertex)$, we consider the ansatz 
\begin{align} \label{eq:ansatz-pv}
P_\vertex \Lambda^{k}_{\bzero} = \bar\Lambda^\vertex + \sum_{m_1, m_2 = 1}^r  \ba_{\mb}  \Lambda^{k}_{\mb} 
    + \sum_{\mathclap{\substack{k' \in \mathcal{K}(\vertex)\\ k' \not = k}}} \ \  \sum_{ m_1, m_2 = 1}^r  \ \tilde \ba_{\mb} \Lambda^{k'}_{  \mb} ,
\end{align}
with coefficients $\ba_\mb = a_{1, m_1} a_{2, m_2}$ -- similarly $\tilde \ba_\mb$ -- to be determined. 
As shown in Remark~\ref{rem:patch-indep},
these correction terms can be determined independently of the patch, which is an important property of our construction
and permits us to omit patch-dependent notation here. 
In the case $r=0$, this ansatz is consistent with the prior definition.

The pullbacks on the patch $\Omega_k$ read
    \begin{align*}
        (\cF^0_k)^{-1}( \Lambda^{k}_{\bi}|_{\Omega_k}) = \hat \Lambda^{k}_{\bi}, \quad 
        (\cF^0_k)^{-1}( \Lambda^{k'}_{\bi}|_{\Omega_k}) = 0,  \\
        (\cF^0_k)^{-1}(P_\vertex \Lambda^{k}_{\bzero} |_{\Omega_k}) = \frac{1}{\#\cK(\vertex)} \hat\Lambda^{k}_{\bzero} + \sum_{m_1, m_2 = 1}^r \ba_\mb \hat\Lambda^{k}_{\mb}, 
    \end{align*}
which gives \eqref{eq:moment-preservation}, for any $q(\hbx) = q_1(\hx_1) q_2(\hx_2) \in \polspace$, as
\begin{align}
    & &&\int_{\hat \Omega}\left( \frac{1}{\#\cK(\vertex)} \hat \Lambda^{k}_{\bzero} + \sum_{m_1, m_2 = 1}^r  \ba_\mb \hat \Lambda^{k}_{\mb} \right) q = \int_{\hat \Omega} \hat \Lambda^{k}_{\bzero} q \\
    &\Leftrightarrow &&  \sum_{m_1, m_2 = 1}^r \ba_\mb \int_{\hat \Omega} \hat \Lambda^{k}_{\mb} q = \left( 1 - \frac{1}{\#\cK(\vertex)}\right) \int_{\hat \Omega}  \hat \Lambda^{k}_{\bzero} q \\
    &\Leftarrow &&\begin{aligned}
        \sum_{m = 1}^r a_{t, m} \int_0^1  \lambda^{k}_{m} q_{t} = \sqrt{ 1 - \frac{1}{\#\cK(\vertex)}} \int_0^1  \lambda^{k}_{0} q_{t}, \quad \text{ for } t = 1, 2,
    \end{aligned}
\end{align}
by the tensor-product structure of the basis functions \eqref{basis_pw}.
On the other hand, on any patch $\Omega_{k'}$ with $k' \in \cK(\vertex) \setminus \{k\}$, a similar derivation yields 
\begin{align}
        \begin{aligned}
            \sum_{m = 1}^r  \tilde a_{t, m} \int_0^1  \lambda^{k'}_{m} q_{t} &= (-1)^t \sqrt\frac{1}{\#\cK(\vertex)} \int_0^1  \lambda^{k'}_{0} q_{t} , \quad \text{ for } t = 1, 2,
        \end{aligned}
\end{align}
where the sign $(-1)^t$ is an arbitrary choice to distribute the sign.
It remains to solve the linear system for the correction coefficients $\ba$ and $\tilde \ba$. 
Let $(q_j)_{1 \le j \le r}$ be the basis of the univariate polynomial space $\scalarpolspace$
and define 
the basis-polynomial duality-matrices and vertex-basis-polynomial vectors  as 
\begin{align} \label{eq:pol_mat}
    &\arrM^{k} = \left(\int_0^1 \lambda^{k}_{m} q_{j} \right)_{1 \le j \le r, 1 \le m \le r},
      &&\arrb^{k} = \left( \int_0^1  \lambda^{k}_{0} q_{j} \right)_{1 \le j \le r}.
\end{align}
We can thus write the linear systems as
\begin{align} \label{eq:lin-sys}
    \arrM^k  \ba_t = \sqrt{ 1 - \frac{1}{\#\cK(\vertex)}} \arrb^{k}, 
    \quad   \arrM^{k'}  \tilde \ba_t &=  (-1)^t \sqrt\frac{1}{\#\cK(\vertex)} \arrb^{k'}.
\end{align}

\begin{remark} \label{rem:patch-indep}
    By Assumption~\ref{as:e_nested}, we can find a uniform refinement factor $z \in \NN_{>0}$ 
    between the resolutions of two patches $k$ and $k'$, such that we have a scaling relation between 
    one-dimensional basis functions on the logical tensor-product domain:  
       $ \lambda^k_i(x) = \lambda^{k'}_{i'}(z x).$
    Performing a change of basis in \eqref{eq:pol_mat} by $q_j(x) = \tilde q_j(z x)$, i.e.
    \begin{equation} \label{eq:scaling}
        \int_0^1  \lambda^k_i q_j = \int_0^1 \frac{ \lambda^{k'}_{i'} \tilde q_j}{z},
    \end{equation}
    this scaling factor cancels out in \eqref{eq:lin-sys}; hence the result is independent of the patch $k$ and direction $t$.
\end{remark}

We define the one-dimensional correction coefficient 
\begin{equation} \label{eq:gamma1d}
    \gammaod = (\arrM^k)^{-1}\arrb^{k} \in \RR^r,
\end{equation} 
and the tensor-product correction coefficient 
\begin{equation} \label{eq:gamma}
    \bgamma =  \gammaod \otimes \gammaod.
\end{equation}
These coefficients are independent of the patch $k\in \cK(\vertex)$ and of the map $F_k$, which is an important property of our construction.
Thus, the solution of the linear system \eqref{eq:lin-sys} yields 
\begin{equation}
    \ba = \left(1 - \frac{1}{\#\cK(\vertex)} \right) \bgamma, \quad \tilde \ba = -\frac{1}{\#\cK(\vertex)} \bgamma.
\end{equation}

The full operator for a given vertex $\vertex$ then reads:
\begin{definition}\label{eq:Pv}
The local vertex-conforming projection $P_\vertex: V^0_\pw \to V^0_\pw$ is defined as
\begin{equation} 
    P_\vertex: \Lambda^k_\bj \mapsto \left\{ \begin{aligned} 
                            &\bar\Lambda^\vertex 
                            + \sum_{m_1, m_2 = 1}^r  \bgamma_{\mb}   \left( \Lambda^{k}_{\mb} -  \frac{1}{\#\cK(\vertex)}  \sum_{k' \in \mathcal{K}(\vertex)} \Lambda^{k'}_{ \mb}\right)
                               &&\text{ if } \bj = \bzero \text{ and } k \in \cK(\vertex), \\
                            &\Lambda^k_\bj &&\text{ else,}
    \end{aligned}\right.
\end{equation}
for all $\bj \in \bcI^k$.
Here, $\Lambda^k_\bj$ are the patch-wise basis functions defined in \eqref{basis_pw}, 
the order of moment preservation is $0 \le r \le p+1$,
$\bgamma$ contains the correction coefficients in \eqref{eq:gamma}, and $\bar\Lambda^\vertex$ is the vertex-conforming basis function in \eqref{eq:vertex-basis}.
\end{definition}

\begin{remark} \label{rem:nbc}
    Assumption~\ref{ass:nbc} gives 
    \begin{equation}
        r \le p+1 < n_k,
    \end{equation}
    hence the correction terms 
    do not involve other vertex basis functions. For $r=0$, the sum is empty, and the definition reduces to \eqref{eq:Pv_r0}.
\end{remark}

We now prove Proposition~\ref{prop:vertex}.
    \begin{proof}[Proof of Proposition~\ref{prop:vertex}]  $ $\par\nobreak\ignorespaces
        Let $\phi = \sum_k \sum_\ib \phi^k_\ib \Lambda_\ib^k \in V^0_\pw$ and $\vertex \in \vertices$.
     \begin{enumerate}
        \item  To show continuity of $P_\vertex \phi$ at $\vertex$, we write its projection on a patch $\Omega_{k'}, k' \in \cK(\vertex)$, as
        \begin{align}
        (P_\vertex \phi)|_{ \Omega_{k'}}(\vertex)
            = \sum_{k \in \cK(\vertex)}   \phi^{k}_{\bzero}  \ \bar \Lambda^\vertex |_{ \Omega_{k'}}(\vertex) 
            = \sum_{k \in \cK(\vertex)}   \frac{\phi^{k}_{\bzero}}{\#\cK(\vertex)}  \underbrace{\Lambda^{k'}_{\bzero}(\vertex)}_{= 1} ,
        \end{align}
        where all interior basis functions in Definition~\ref{eq:Pv} vanish at the vertex $\vertex$.
        Since the restriction to any patch $k' \in \cK(\vertex)$ yields the same value, it proves continuity at $\vertex$.

        \item  It remains to show that if $\phi$ is continuous at $\vertex$, then $P_\vertex \phi = \phi$ as $(1)$ implies the other direction. 
         If $\phi$ is continuous at $\vertex$, all its vertex coefficients have a common value $\phi_\bzero$. Summing the images of the corresponding vertex basis functions yields 
            \begin{align}
                \sum_{k \in \cK(\vertex)} P_\vertex \Lambda^{k}_{\bzero} = \#\cK(\vertex)   \bar\Lambda^\vertex 
                +  \sum_{m_1, m_2 = 1}^r  \bgamma_{\mb}   \left( \sum_{k \in \cK(\vertex)} \Lambda^{k}_{\mb} -   \sum_{k' \in \mathcal{K}(\vertex)} \Lambda^{k'}_{ \mb}\right)  
                = \sum_{k \in \cK(\vertex)} \Lambda^{k}_{\bzero}.
            \end{align}
        All other basis functions are unchanged, hence $P_\vertex \phi = \phi$.

        \item Moment preservation holds by construction in Section~\ref{subsec:vertex-pol-mom}.
        
        \item Remark~\ref{rem:nbc} shows that the vertex-based operators commute for different vertices, as the projection  $P_\vertex$ does not change coefficients of other vertex basis functions.

        \end{enumerate}

     \end{proof}

  \subsection{Edge-based projections}
  \label{sec:p_edge}
  For this section, we again let the moment preservation order be given by $0 \le r \le p + 1$, fix an interface $\edge \in \edges$ and denote by $k^- \in \cK(\edge)$
the coarse and $k^+\in \cK(\edge)$ the fine patch adjacent to the edge $\edge$, where we drop the $k$ in superscripts and just write $\pm$. Since the ordering of the patch-wise basis functions depends on the orientation of the edge, we have to distinguish between vertical 
and horizontal interfaces, but without loss of generality we can assume for the following definitions that
the first index is the edge-parallel coordinate and the second index the edge-perpendicular coordinate, i.e.
\begin{equation}
\bj = (j_1, j_2), \quad j_1 = j_{\parallel(\hat \edge)}, \quad j_2 = j_{\perp(\hat \edge)}.
\end{equation}
Furthermore, we assume that for both patches the basis functions in the edge-perpendicular direction 
are ordered starting with the edge-basis functions associated with the edge $\edge$.
At the same time, we assume that the basis functions in the edge-parallel 
direction are ordered in the same direction on both adjacent patches. 

For notational convenience, we consider an edge-local parametrization of the patches $k^-$ and $k^+$
consisting of a single mapping $F_\edge: \hat \Omega(\edge) \rightarrow \Omega(\edge)$,
such that the edge-local logical domain $\hat \Omega(\edge) = [0, 1] \times [-1, 1]$,
is the juxtaposition of the tensor-product patches $\hat\Omega_\pm$,
and the edge-local physical domain $\Omega(\edge)$,
  is the juxtaposition of the physical patches $\Omega_\pm$. The mapping $F_\edge$ is a
  homeomorphism and a $C^1$-diffeomorphism on each logical patch.
  We assume that the restrictions of $F_\edge$ to $\hat\Omega_\pm$ coincide with the original patch parametrizations $F_{\pm}$, up to affine coordinate transformations accounting for the orientation of the common edge.
 An illustration is provided in Figure~\ref{fig:edge_domain}.

\begin{figure}[!ht]
    \centering
    \begin{tikzpicture}[lattice/.cd,spacing/.initial=5,superlattice
	period/.initial=30,amplitude/.initial=3]
   
	\begin{scope}[yshift=2cm, xshift=4cm]
		\draw[ step = 1.0] (0,0) grid (4,4);
	\end{scope}

  \begin{scope}[yshift=-2cm, xshift=4cm]
	  \draw[ step = 2.0] (0,0) grid (4,4);

	  \draw[->] (-0.5,-0.5)--(0.5,-0.5) node[right]{$\hat x_\parallel$};
	  \draw[->] (-0.5,-0.5)--(-0.5,0.5) node[above]{$\hat x _\perp$};

	  \draw[thick, color=red] (0, 4) -- (4, 4);
  \end{scope}
  
  \begin{scope}[xshift=-4cm]
	\pgftransformnonlinear{\latticetilt}
	  
	  \draw[ step = 2.0, yshift=-2cm] (0,0) grid (4,4);
	  
	  
	  
	  \draw[ step = 1.0, yshift=2cm] (0,0) grid (4,4);
	  
	  
	  \draw[thick, color=red] (0, 2) -- (4, 2);
  \end{scope} 
		  
  	\begin{scope}
		\node[xshift=-4cm] at (2,-3) {$\Omega(e)$};
		\node[xshift=4cm] at (2, -3) {$\hat \Omega(e)$};

		\node[xshift=-4cm] at (0.5, 2.5) {$\Omega_+$};
		\node[xshift=-4cm] at (0.5,-1.5) {$\Omega_-$};

		\node[xshift=4cm] at (3.5, 2.5) {$\hat \Omega_+$};
		\node[xshift=4cm] at (3.5,-1.5) {$\hat \Omega_-$};

		\node[right, xshift=-4cm] at (4.25, 2) {$e$} ;
		\node[left, xshift=4cm] at (-0.25, 2) {$ \hat e$} ;
		
		\draw[->] (4.25, 1.75) -- (4.75, 1.75) node[right, scale=0.75]{$j^-_1$};
		\draw[->] (4.25, 1.75) -- (4.25, 1.25) node[below, scale=0.75]{$j^-_2$};

		\draw[->] (4.25, 2.2) -- (4.75, 2.2) node[right, scale=0.7]{$j^+_1$};
		\draw[->] (4.25, 2.2) -- (4.25, 2.7) node[above, scale=0.7]{$j^+_2$};

		\draw[-{Latex[length=3mm]}] (3.75, 2.5) to [out=150,in=30]  node[midway, above]{$F_\edge$} (0.25, 2.5);

	\end{scope}
		  

   
		  
  \end{tikzpicture}
    \caption{The edge-local logical domain $\hat \Omega(\edge)$ (right) and the edge-local physical domain $\Omega(\edge)$ (left) 
    for a horizontal edge $\edge$, different resolutions and the mapping $F_\edge$. 
    Additionally, we show the orientation of the logical variables and indicate the direction of increasing 
    basis index ordering.}
    \label{fig:edge_domain}
\end{figure}

This extends the patch-wise logical basis functions in edge-perpendicular direction to the whole 
edge-local logical domain $\hat \Omega(\edge)$:
\begin{equation} \label{eq:e-log-basis}
    \lambda_{\perp, i}^-(\hat x_\perp) = \left\{ \begin{aligned}
        &\lambda^-_i(-\hat x_\perp) && \text{if $-1 \le \hat x_\perp \le 0$}\\
        & 0 && \text{if $0 < \hat x_\perp \le 1$}
    \end{aligned}\right., \quad 
    \lambda_{\perp, i}^+(\hat x_\perp) = \left\{ \begin{aligned}
        & 0 && \text{if $-1 \le \hat x_\perp \le 0$} \\
        &\lambda^+_i(\hat x_\perp) && \text{if $0 < \hat x_\perp \le 1$}
    \end{aligned}
    \right.
\end{equation} 
for all $0 \le i \le n_\pm$. These functions span the shifted univariate function spaces perpendicular to the edge $\hat V^{0}_{\perp, \pm} = \hat \VV^0_{\perp, \pm}$ in the edge-local logical domain $\hat \Omega(\edge)$.
On the other hand, the nested univariate function spaces along the edge $\hat \edge$, 
i.e.  $\hat V^{0}_{\parallel, \pm} = \VV^0_\pm$,
are naturally defined on the edge-local logical domain $\hat \Omega(\edge)$, but
have to be corrected for continuity by the conforming projection as they in general have different resolutions.

\subsubsection*{Extension-Restriction operators} 
\label{sec:ext-res}
In order to combine spline spaces of different refinement levels on shared interfaces,
we introduce projections between the univariate logical edge-parallel function spaces $\hat V^{0}_{\parallel, \pm} = \VV^0_\pm$. 
We refer to them as extension operators $\scrE$ (inclusion/change of basis) when mapping from the coarse spline space to the fine spline space,
and restriction operators $\scrR$ ($L^2$-projection) for the opposite direction. \\

We define the extension operator $\scrE$, mapping coarse one-dimensional spline functions to the fine spline space, as a map
\begin{equation} \label{eq:extension}
    \scrE: \VV^0_- \subset \VV^0_+ \longrightarrow \VV^0_+, 
\end{equation}

that acts as 
the change of basis describing the (coarse) basis functions of $\VV^{0}_-$ in the (fine) basis of $\VV^0_+$, 
which is possible due to Assumption~\ref{as:e_nested}.
Its matrix form $\arrE$, the change of basis matrix, is defined in Appendix~\ref{sec:impl_erv0}. 

The restriction operator $\scrR$ is a projection from the fine to the coarse spline space:
\begin{equation} \label{eq:restr}
    \scrR:\VV^0_+ \longrightarrow \VV^0_- \subset\VV^0_+,
\end{equation}
naturally this will be a left inverse of the extension operator, i.e. $\scrR \scrE = \iden$.
Since the edge-local conforming projection should not change the vertex degrees of freedom, 
as we mentioned in the beginning, we have additional constraints:
\begin{equation} \label{eq:vertex-constr}
    \scrR \lambda^+_{0}(0) = \scrR \lambda_0^-(0) = \lambda_0^-(0), \quad \scrR \lambda^+_{n_+}(1) = \scrR \lambda_{n_-}^-(1) = \lambda_{n_-}^-(1).
\end{equation}

First, we introduce the interior spaces $\VV_{\pm, 0}^{0} = \Span \left( \{\lambda^\pm_i \ | \  0 < i < n_\pm \} \right)$ 
and write the restriction operator as
\begin{equation} \label{eq:def-R-0}
    \scrR \hat \phi = \left( \iden - \cT\right)\hat \phi + \scrR_0 \cT \hat \phi,
\end{equation}
depending on the interior restriction operator $\scrR_0$ constructed below, and
where the truncation $\cT :\VV^0_+ \longrightarrow \VV^0_{+, 0}$ is defined as
\begin{align}
    \label{eq:trunc-1}
   \cT \hat \phi &= \hat \phi - \hat \phi_0 \lambda_0^- - \hat \phi_{n_+} \lambda_{n_-}^-, \\
\end{align}
note that $(\cT \hat \phi)(0) = (\cT \hat \phi)(1) = 0$
as the basis functions are interpolatory at the boundaries.
The interior restriction operator $\scrR_0: \VV^{0}_{+,0} \longrightarrow \VV^0_{-, 0}$ 
is defined on the interior fine space and thus does not involve vertex degrees of freedom.
This interior restriction operator $\scrR_0$ is defined using the full $L^2$-projection $\tilde \scrR$
from the interior fine space $\VV^0_{+, 0}$ onto the coarse space $\VV^0_{-}$:
\begin{equation} \label{eq:l2-r0}
    \int_0^1 \lambda^-_j \tilde \scrR \lambda^+_i = \int_0^1\lambda^-_j \lambda^+_i, \quad \forall j=0,\dots, n_-,i = 1, \dots, n_+-1,
\end{equation} 
which is moment-preserving since $\VV^0_{\pm}$ contains all polynomials of degree $p \ge r-1$ as mentioned in Remark~\ref{rem:mom-pres}.
To remove the vertex degrees of freedom from the range of the projection in a moment-preserving way, we 
introduce another truncation $\tilde \cT: \VV^0_- \longrightarrow \VV^0_{-, 0}$, defined as
\begin{equation}
    \label{eq:trunc-2}
    \tilde \cT   \hat \phi =  \hat \phi - \hat \phi_0 \lambda^-_0 - \hat \phi_{n_-} \lambda^-_{n_-} + \hat \phi_0 \sum_{i = 1}^r \gammaod_i  \lambda^-_i + \hat \phi_{n_-} \sum_{i = n_--r}^{n_--1} \gammaod_{n_- - i}  \lambda^-_i,
\end{equation}
where $\hat \phi \in \VV^0_{-}$. 
This correction accounts for the polynomial moments: for $\hat \phi^- \in \VV^0_{-}$ and $ q \in \scalarpolspace$,
\begin{equation}
    \int_0^1 \tilde\cT \hat \phi q = \int_0^1 \hat \phi q,
\end{equation}
by the definition of the correction coefficients $\gammaod$.
Hence, we set 
\begin{equation} \label{eq:int_rest}
    \scrR_0 = \tilde \cT \tilde \scrR.
\end{equation}

The moment preservation 
of the full restriction operator $\scrR$ in \eqref{eq:def-R-0} follows from the identities
\begin{align}
    \int_0^1 (\scrR \hat \phi) q &= \int_0^1 \left( (\iden - \cT)\hat \phi + \scrR_0 \cT \hat \phi \right) q \\
    &= \int_0^1 \left( (\iden - \cT)\hat \phi  +  \tilde \cT \tilde\scrR \cT \hat \phi \right) q \\   
    &= \int_0^1 \left( (\iden - \cT)\hat \phi  +  \tilde\scrR \cT \hat \phi \right) q \\ 
    &= \int_0^1 \left( (\iden - \cT)\hat \phi  + \cT \hat \phi \right) q \\
            &= \int_0^1 \hat \phi q,
\end{align}
holds for all $q \in \scalarpolspace$,
where we used the moment preservation of $\tilde \cT$ and $\tilde \scrR$.
Equation \eqref{eq:vertex-constr} then also holds by construction.

Note that this indeed defines a projection: Let $\hat \phi \in \VV^0_{-} \subset\VV^0_+$, then
\begin{equation}
    \cT \hat \phi = \hat \phi - \hat \phi_0 \lambda_0^- - \hat \phi_{n_+} \lambda_{n_-}^- \in \VV^{0}_{-, 0},
\end{equation}
thus $\scrR_0 \cT \hat \phi = \cT \hat \phi$ by definition of the $L^2$-projection and $\tilde \cT$, hence
\begin{equation}
    \scrR \hat \phi = \hat \phi_0 \lambda^-_0 + \hat \phi_{n_+} \lambda^-_{n_-} + \cT \hat \phi = \hat \phi.
\end{equation}

To summarize, we constructed the restriction operator 
\begin{equation}
    \label{eq:restr-def}
    \scrR:\VV^0_+ \longrightarrow \VV^0_-,
    \hat \phi \mapsto \scrR \hat \phi = \hat \phi_0 \lambda^-_0 + \hat \phi_{n_+} \lambda^-_{n_-} + \tilde \cT \tilde \scrR \cT \hat \phi,
\end{equation}
consisting of the two truncation operators $\cT$ defined in equation \eqref{eq:trunc-1} and $\tilde \cT$ defined in equation \eqref{eq:trunc-2}
and the $L^2$-projection $\tilde \scrR$ from equation \eqref{eq:l2-r0}.
 This operator is a projection onto the univariate logical edge-parallel space $\VV_-^{0}$ that
 preserves the vertex degrees of freedom and preserves polynomial moments up to order $r$.

The description of its matrix $\arrR$ and its detailed construction are given in the Appendix~\ref{sec:impl_erv0}.

\subsubsection*{Local edge-interior conforming projection} 
We introduce the averaged basis function perpendicular to the edge $\edge$:
\begin{equation} \label{eq:bar-lambda}
    \bar\lambda_{\perp}^\edge = \frac{1}{2}  \left(\lambda^{-}_{\perp, 0}  + \lambda^{+}_{\perp, 0} \right),
\end{equation}
which uses the parametrization in \eqref{eq:e-log-basis}. 
For coarse basis functions, the preliminary
averaging operator modifies only functions supported on the interface but not
at its vertices

\begin{equation}
   \bar P_\edge :  \Lambda^{-}_\jb \mapsto \left\{ \begin{aligned}
                            &\cF^0_\edge \left( \lambda^{-}_{j_1} \otimes \bar\lambda_{\perp}^\edge \right)
                            &&\text{ if } 0 < j_1< n_-, j_2= 0, \\
                            &\Lambda^{-}_\jb &&\text{ else,}
    \end{aligned}\right.
\end{equation}
note that this is well-defined since $\hat V^{0}_{\parallel, -} \subset \hat V^{0}_{\parallel, +}$. The application of the extension operator is left implicit here.
Here and below, $\cF^0_\edge$ denotes the scalar pushforward associated with the
edge-local mapping $F_\edge$ from Figure~\ref{fig:edge_domain}, defined analogously to $\cF_k^0$ in \eqref{pf_gc}.

For fine basis functions, 
we need to additionally apply the restriction operator in the parallel direction: 
\begin{equation}
    \bar P_\edge :  \Lambda^{+}_\jb \mapsto \left\{ \begin{aligned}
        &\cF^0_\edge \left( \scrR \lambda^{+}_{j_1} \otimes \bar\lambda_{\perp}^\edge \right)
        &&\text{ if } 0<j_1< n_+, j_2= 0, \\
        &\Lambda^{+}_\jb &&\text{ else.}
\end{aligned}\right.
\end{equation}
The final operator $P_\edge$ augments these averages with moment corrections.

\subsubsection*{Preservation of polynomial moments}\label{sec:mom-pres-Pe}

Following Definition~\ref{mom_pres}, as in Section~\ref{subsec:vertex-pol-mom}, we add 
 (interior) correction terms in the perpendicular direction of the interface, which are independent of the patch as explained in Remark~\ref{rem:patch-indep}. 
As for the vertex-conforming projection, we have to map coarse and fine basis functions on the edge, namely $\Lambda^{-}_{(i, 0)}$ and $\Lambda^{+}_{(i, 0)}$,
pulling them back to the logical domain from  both of the tensor-product patches, namely $\Omega_{-} $ and $\Omega_{+}$, and checking for the moment preservation. 
Since we want to preserve vertex degrees of freedom, we only consider $0 < i < n_\pm$ for now. 

Starting with a coarse basis function, the ansatz reads
\begin{equation*}
    P_\edge \Lambda^{-}_{(i, 0)}  =  \cF^0_\edge \left( \lambda^{-}_{i} \otimes  \left( \bar\lambda_{\perp}^\edge     + \sum_{m = 1}^r a_m  \lambda^{-}_{\perp, m} +  \sum_{ m = 1}^r \tilde a_{ m} \lambda^{+}_{\perp, m} \right) \right).
\end{equation*}

Thus, for any $q \in \polspace$, the condition for the coarse patch $\Omega_{-}$ reads
\begin{align}
     \sum_{m = 1}^r a_m \int_{\hat \Omega_-}     \lambda^{-}_{i} \otimes \lambda^{-}_{\perp,m} q 
            &= \frac{1}{2}\int_{\hat \Omega_-} \lambda^{-}_{i} \otimes \lambda^{-}_{\perp,0} q.
\end{align}
%
%
On the fine patch $\Omega_{+}$, for any $q \in \polspace$, we have
\begin{align}
    \sum_{ m = 1}^r \tilde a_{ m} \int_{\hat \Omega_+}   \scrE \lambda^{-}_{i} \otimes \lambda^{+}_{\perp, m}   q 
    = - \frac{1}{2} \int_{\hat\Omega_+} \scrE \lambda^{-}_{i} \otimes  \lambda^{+}_{\perp, 0} q.
\end{align}

The parallel component of the integrals is equal on both sides, so no correction is needed in that direction. 
The correction coefficients along the perpendicular direction are given by the same formula as in the vertex case, i.e. \eqref{eq:gamma1d}.

Introducing the patch-wise correction term $\ttC^\edge$ as
\begin{equation}\label{eq:Ce}
    \ttC^\edge = \frac{1}{2} \sum_{m = 1}^r  \gammaod_{m} \left(  \lambda^{-}_{\perp, m} -  \lambda^{+}_{\perp, m} \right),    
\end{equation}
the moment-preserving edge-conforming operator $P_\edge$ for coarse edge-interior basis functions is defined as
\begin{equation}
    P_\edge \Lambda^-_{(i, 0)} = \cF^0_\edge \left( \lambda^{-}_{i} \otimes \left( \bar\lambda_{\perp}^\edge  + \ttC^\edge  \right)  \right), \quad 0 < i < n_-.
\end{equation}

Analogously, the ansatz for fine basis functions reads
\begin{equation*}
    P_\edge \Lambda^{+}_{(i, 0)} = \cF^0_\edge \left( \scrR \lambda^{+}_{i} \otimes \left( \bar\lambda_{\perp}^\edge     + \sum_{m = 1}^r a_m \lambda^{+}_{\perp, m} +  \sum_{ m = 1}^r \tilde a_{  m}  \lambda^{-}_{\perp, m} \right) \right).
\end{equation*}

For any $q \in \polspace$, on the coarse patch $\Omega_{-}$, we obtain
\begin{align}
    \sum_{ m = 1}^r \tilde a_{ m} \int_{\hat \Omega_-}  \scrR \lambda^{+}_{i} \otimes \lambda^{-}_{\perp, m} q = -\frac{1}{2} \int_{\hat \Omega_-} \scrR \lambda^{+}_{i} \otimes \lambda^{-}_{\perp,0}  q,
\end{align}
whereas before, the parallel component cancels, leading only to one-dimensional correction terms in the perpendicular direction. 
On the other hand, on the fine patch $\Omega_{+}$, we get
\begin{align}
    \sum_{m = 1}^r a_m \int_{\hat \Omega_+}  \scrR \lambda^{+}_{i} \otimes \lambda^{+}_{\perp,m}  q = \int_{\hat \Omega_+} \left( \lambda^{+}_{i} - \frac{1}{2} \scrR \lambda^{+}_{i} \right) \otimes \lambda^{+}_{\perp,0} q,
\end{align}
for any $q \in \polspace$.
Here, we require moment preservation of the restriction operator $\scrR$, as derived in Section~\ref{sec:ext-res}; this makes the parallel component cancel as before.
Using the same correction term $\ttC^\edge$ defined in \eqref{eq:Ce},
the moment-preserving edge-conforming operator for fine edge-interior basis functions is defined as
\begin{equation} \label{eq:Pe_lam_+}
    P_\edge \Lambda^+_{(i, 0)} = \cF^0_\edge \left( \scrR \lambda^{+}_{i} \otimes \left( \bar\lambda_{\perp}^\edge  - \ttC^\edge  \right) \right), \quad 0 < i < n_+.
\end{equation}

The preceding derivations did not include the vertex basis functions (on the edge) because their direct correction would alter the degree of freedom at the vertex. 
In order to also preserve moments across the edge for these basis functions, we have to add another correction term,
where the idea is to consider a linear combination of logical basis functions in parallel direction that are supported on the edge, but not on the vertex, 
such that it has the same moments as the vertex basis function:
Assume the vertex $\vertex$ is one of the two vertices on the edge $\edge$ at (basis function) index $\bzero$ on both patches
and 
define a fine spline $ \mu_\vertex^{+} \in \hat V^0_{\parallel, +}$ in parallel coordinate such that $\mu_\vertex^{+}(\vertex_\parallel)= 0$ and
\begin{equation}
    \int_0^1 \mu_\vertex^{+} q = \int_0^1  \lambda_{0}^{+
    } q, \qquad \forall q \in \scalarpolspace ,
\end{equation}
holds.
This, again, leads to defining  
\begin{equation} \label{eq:mu-plus}
    \mu_\vertex^{+} = \sum_{m = 1}^r \gammaod_{m} \lambda^{+}_m.
\end{equation}
In the following proposition, we derive the existence and explicit form of an analogous coarse spline $\mu_\vertex^{-}$
and characterize the edge-based projection for the vertex basis functions of both patches: 
 
\begin{prop} \label{prop:mu-minus}
    There exists a coarse spline $\mu_\vertex^{-} \in \VV^0_-$, such that defining the edge-local projection of the fine vertex basis functions as 
    \begin{align} \label{eq:p_e_lam_+_0}
        P_\edge \Lambda^+_{\bzero} = \cF^0_\edge\left( \scrR \lambda^+_{0} \otimes \bar\lambda_{\perp}^\edge  
        - \scrR \mu_\vertex^+ \otimes \ttC^\edge \right),
    \end{align}
    where $\mu_\vertex^{+}$ is defined in \eqref{eq:mu-plus}, and of the coarse vertex basis functions as
    \begin{align} \label{eq:p_e_lam_-_0}
        P_\edge \Lambda^-_{\bzero} = \cF^0_\edge \left( \lambda^-_{0} \otimes \bar\lambda_{\perp}^\edge  
         + \mu_\vertex^- \otimes \ttC^\edge \right),
    \end{align}
    where $\mu_\vertex^{-}$ is derived in \eqref{eq:mu-minus} and $\ttC^\edge$ is defined in \eqref{eq:Ce},
    preserves polynomial moments of order $r$ and preserves continuity at the vertex. 
\end{prop}

\begin{proof}
Consider the continuous function $\phi \in V^0_h$ with pullback
\begin{equation}
    (\cF^0_\edge )^{-1} (\phi |_{\Omega(\edge)}) = \lambda^-_{0} \otimes \lambda^-_{\perp, 0}  + \scrE  \lambda^-_{0} \otimes \lambda^+_{\perp,0} 
\end{equation}
 which does not vanish on the vertex $\vertex$ and should stay unchanged under the projection.
 Since we can decompose $\scrE \lambda^-_0 = \lambda^+_{0} + (\scrE \lambda^-_{0} - \lambda^+_{0})$, such that the second term is only supported in the edge interior, 
 we can use the definitions \eqref{eq:Pe_lam_+}, \eqref{eq:p_e_lam_+_0} and \eqref{eq:p_e_lam_-_0} to write the projection of $\phi$ on the edge as
\begin{align}
    (\cF^0_\edge)^{-1} (P_\edge \phi |_{\Omega(\edge)}) &= \lambda^-_{0} \otimes \bar\lambda_{\perp}^\edge  + \mu_\vertex^- \otimes \ \ttC^\edge\\
       &+ \scrR \lambda^+_{0} \otimes \bar\lambda_{\perp}^\edge  - \scrR \mu_\vertex^+ \otimes \ttC^\edge\\
       &+  \scrR\left( \scrE \lambda^{-}_{0} - \lambda^{+}_{0}\right) \otimes \left( \bar\lambda_{\perp}^\edge - \ttC^\edge\right) \\
        &= \lambda^-_{0} \otimes 2 \bar\lambda_{\perp}^\edge 
        + \left(\mu_\vertex^- - \scrR \mu_\vertex^+ -\lambda^{-}_{0} + \scrR \lambda^{+}_{0}  \right) \otimes \ttC^\edge,
    \end{align}
where we used that $\scrR \scrE = \iden$.
    Subtracting the original function $\phi$ then leads to the condition

\begin{equation}
    0 = (P_\edge \phi - \phi)|_{\Omega(\edge)} = \cF^0_\edge \left( \left(\mu_\vertex^- - \scrR \mu_\vertex^+ -\lambda^{-}_{0} + \scrR \lambda^{+}_{0}  \right) \otimes \ttC^\edge \right),
\end{equation}
where we conclude with defining $\mu_\vertex^-$ as
\begin{equation} \label{eq:mu-minus}
    \mu_\vertex^- =  \lambda^{-}_{0} + \scrR \left( \mu_\vertex^+ - \lambda^{+}_{0}  \right),
\end{equation}
as the edge-conforming projection $P_\edge$ should leave $\phi$ unchanged. 

Since $\mu^\pm_\vertex(\vertex) = 0$ and $\ttC^\edge(\vertex) = 0$, vertex-continuity is preserved. Indeed, if $\phi = \sum_k \sum_\bi \phi_\bi^k \Lambda^k_\bi \in V^0_\pw$ is continuous at the vertex $\vertex$, then 
\begin{align}
    (P_\edge \phi)|_{\Omega_\pm}(\vertex) = \phi_\bzero \cF^0_\edge\left( \scrR \lambda_0^+ \otimes \bar\lambda_{\perp}^\edge + \lambda^-_0 \otimes \bar\lambda_{\perp}^\edge\ \right)|_{\Omega_\pm}(\vertex) = \phi_\bzero.
\end{align}
The preservation of polynomial moments comes from the ansatz as before and from the fact that $\mu^+_\vertex$ is constructed to preserve the moments of $\lambda^+_0$, which directly transfers to $\mu^-_\vertex$ preserving the moments of $\lambda^-_0$.
\end{proof}

The construction for the second vertex on the edge $\edge$ is analogous. This gives the following definition.
\begin{definition} \label{eq:Pe}
    The local edge-conforming projection $P_\edge: V^0_\pw \rightarrow V^0_\pw$ is defined as
\begin{align}
    P_\edge : \left\{
        \begin{aligned}
            &\Lambda^{-}_\jb \mapsto \left\{ \begin{aligned}
                                    & \begin{aligned} & \cF^0_\edge \left( \lambda^{-}_{j_1} \otimes \left( \bar\lambda_{\perp}^\edge   +  \ttC^\edge \right) \right)
                                    \end{aligned}
                                                &&\text{ if } 0 < j_1< n_-,~ j_2= 0, \\
                                    & \begin{aligned} &\cF^0_\edge \left( \lambda^-_{j_1} \otimes \bar\lambda_{\perp}^\edge + \mu_\vertex^- \otimes \ttC^\edge \right)
                                    \end{aligned}
                                                &&\text{ if } j_1 \in \{ 0, n_-\},~ j_2 = 0, \text{ for } \vertex \text{ at } \jb, \\
                                    &\Lambda^{-}_\jb &&\text{ else,}
                                \end{aligned}\right. 
             \\
            &\Lambda^{+}_\jb \mapsto \left\{ \begin{aligned}
                                    &\begin{aligned} &\cF^0_\edge \left( \scrR \lambda^+_{j_1} \otimes \left( \bar\lambda_{\perp}^\edge - \ttC^\edge\right) \right)
                                    \end{aligned}&&\text{ if } 0 < j_1< n_+,~ j_2= 0, \\
                                    &\begin{aligned} & \cF^0_\edge \left( \scrR \lambda^+_{j_1} \otimes \bar\lambda_{\perp}^\edge - \scrR \mu_\vertex^+ \otimes \ttC^\edge \right)
                                    \end{aligned}&&\text{ if } j_1 \in  \{0, n_+\},~ j_2 = 0, \text{ for } \vertex \text{ at } \jb, \\
                                    &\Lambda^{+}_\jb &&\text{ else,}
            \end{aligned}\right.\\
            &\Lambda^{k}_\jb \mapsto \left. \begin{aligned}
                                    &&\Lambda^{k}_\jb &&\text{ if } k \not= k^-, k^+, \\
            \end{aligned}\right.\\
        \end{aligned} \right.  \\
\end{align}


where $\vertex$ denotes the vertex corresponding to the index $\jb$ in both intermediate cases, $\lambda^k_\bj$ are the scalar patch-wise basis functions defined in \eqref{eq:basis-hatV0}, $\bar\lambda^\edge_\perp$ the edge-conforming spline in \eqref{eq:bar-lambda},
$\ttC^\edge$ the patch-wise correction term in \eqref{eq:Ce}, $\mu_\vertex^-$ and $\mu_\vertex^+$ the vertex correction terms in \eqref{eq:mu-minus} and \eqref{eq:mu-plus}, 
and the restriction operator $\scrR$ in \eqref{eq:restr}. On other patches $k \not= k^-, k^+$, it acts as the identity. 
\end{definition}

Using this definition, we conclude the section with the proof of Proposition~\ref{prop:edge}.



    \begin{proof}[Proof of Proposition~\ref{prop:edge}] $ $\par\nobreak\ignorespaces
        Let $\phi = \sum_k \sum_\ib \phi^k_\ib \Lambda_\ib^k  \in V^0_\pw$ be continuous on vertices and fix an edge $\edge \in \edges$.
            Furthermore, let
            \begin{equation}
                \hat \phi^-_\parallel = \sum_{i = 0}^{n_-} \phi^-_{(i, 0)} \lambda^-_i, \qquad
                \hat \phi^+_\parallel = \sum_{i = 0}^{n_+} \phi^+_{(i, 0)} \lambda^+_i.
            \end{equation}
            
        \begin{enumerate}
        \item Since $\ttC^\edge$ vanishes on the edge, evaluating $P_\edge \phi$ on $\edge$ gives
            \begin{align}
                (\cF^0_\edge)^{-1} \left( P_\edge \phi |_{\Omega(\edge)} \right)|_{\hat \edge} &= \left. \left( \left( \hat \phi^-_\parallel + \scrR \hat \phi^+_\parallel \right) \otimes\bar\lambda_{\perp}^\edge \right)\right|_{\hat\edge}  \\
                &= \frac{1}{2}\left( \hat \phi^-_\parallel + \scrR \hat \phi^+_\parallel \right) \otimes  \left( 1|_{\Omega_-} + 1|_{\Omega_+} \right) ,
                \end{align}
                    which shows that $P_\edge \phi$ is continuous on $\edge$.

        \item It remains to show that if $\phi$ is continuous on $\edge$, then $P_\edge \phi = \phi$  as the other implication follows from $(1)$. Assume $\phi$ is continuous on the edge $\edge$, which is equivalent to
        \begin{align}
            \hat \phi^-_\parallel = \hat \phi^+_\parallel.
        \end{align}
        This relation is invariant under the restriction operator, i.e. 
        \begin{align} \label{eq:restr-iden}
            \scrR \hat \phi^+_\parallel = \underbrace{\scrR \scrE}_{= \mathrm{id}} \hat \phi^-_\parallel = \hat \phi^-_\parallel  = \hat \phi^+_\parallel,
        \end{align}
        and yields the edge decomposition
        \begin{equation} \label{eq:decomp-edge}
            \hat \phi^\pm_\parallel \otimes \bar\lambda^\edge_\perp
             = \frac{1}{2}  \hat \phi^-_\parallel \otimes \lambda^-_{\perp, 0}
             + \frac{1}{2} \hat \phi^+_\parallel \otimes \lambda^+_{\perp, 0}.
        \end{equation}
        Using Definition \ref{eq:Pe}, the projection of $\phi$ can be written as 
            \begin{align}
                \left(\cF^0_\edge\right)^{-1} \left(P_\edge \phi |_{\Omega(\edge)}\right) =& \sum_{k \in \pm} \sum_{\mathclap{\substack{\ib \in \bcI^k \\ i_2 \not = 0 }}} \phi^k_\ib  \lambda^k_{i_1} \otimes  \lambda^k_{\perp, i_2} \\
                &+ \hat \phi^-_\parallel \otimes \bar\lambda_{\perp}^\edge + \scrR \hat \phi^+_\parallel \otimes \bar\lambda_{\perp}^\edge \\
                &+ \hat \phi^-_\parallel \otimes \ttC^\edge - \scrR \hat \phi^+_\parallel \otimes \ttC^\edge \\
                &- \sum_{i = 0, n_-} \phi^-_{(i, 0)} \lambda^{-}_{i} \otimes \ttC^\edge + \sum_{i = 0, n_+} \phi^+_{(i, 0)} \scrR \lambda^+_{i} \otimes \ttC^\edge \\
                &+ \sum_{i = 0, n_-}  \phi^-_{(i, 0)}  \mu_{\vertex(i)}^- \otimes \ttC^\edge - \sum_{i = 0, n_+}   \phi^+_{(i, 0)} \scrR \mu_{\vertex(i)}^+ \otimes \ttC^\edge. \\
            \end{align}
            Using relation \eqref{eq:restr-iden} together with relation \eqref{eq:decomp-edge}, we recover $\phi$ with the first two lines. The third line vanishes by relation \eqref{eq:restr-iden}, and the last two lines cancel by the definition of $\mu_\vertex^-$ in \eqref{eq:mu-minus}, where we also make use of the continuity at vertices. Note that we wrote $\vertex(i)$, with $i = 0, n_\pm$, to denote the vertex on the edge $\edge$ with associated index $(i,0)$. In total, we showed $P_\edge \phi = \phi$.


        \item Preservation of polynomial moments holds by construction in Section~\ref{sec:mom-pres-Pe}.
        
        \item The edge-based operators for different edges commute since the interiors of all edges are disjoint, 
        the projection $P_\edge$ only acts on basis functions in the edge interior of $\edge$ by construction,
        leaving the rest of the basis functions unchanged. The assumed continuity at vertices yields result.

        \end{enumerate}
    
\end{proof}

  \section{Conforming projection operators on \texorpdfstring{$V^1_\pw$}{V1 patch-wise}}  
  \label{sec:conf-proj-v1}
  
In this section we construct a conforming projection 
\begin{equation}
    \bP^1 : V^1_\pw \to V^1_\pw,
\end{equation}
onto the space $V^1_h = V^1_\pw \cap H(\curl;\Omega)$.
Here the curl conformity amounts to a tangential continuity constraint.
We use the same assumptions on the geometry of the patches as in
  Section~\ref{sec:p_edge},
  and consider the logical domain $\hat \Omega(\edge)$ around a horizontal edge $\hat \edge$
  corresponding to the logical patches $\hat \Omega_\pm$ associated with the edge-local
  parametrization $F_\edge: \hat \Omega(\edge) \rightarrow \Omega(\edge)$,
  see Figure~\ref{fig:edge_domain}.
On every patch, the logical spaces are of the form
  \begin{equation}
      \hat V^1_\pm = \begin{pmatrix}
          \VV^1_\pm \otimes \VV^0_{\perp, \pm} \\
          \VV^0_\pm  \otimes \VV^1_{\perp, \pm}
      \end{pmatrix},
  \end{equation}
  therefore the univariate function spaces along the edge $ \hat \edge$,
  are given by $(\hat V^{1}_\pm \cdot \hat\btau_{ \edge})_{\parallel} = \VV^{1}_\pm$, associated with the edge-parallel direction of the first component of the logical space $\hat V^1_\pm$.
  The corresponding component of the projected fields must 
  then be continuous across the edge $\hat \edge$, that is in the edge-perpendicular direction.

  Hence, Section~\ref{sec:ext-restr-v1} extends the previously introduced 
    extension $\scrE$ and restriction operators $\scrR$ 
    (which are applied in the edge-parallel coordinate) to the space $\VV^1_\pm$. The averaging and moment corrections are still done in the edge-perpendicular direction, i.e. on $\VV^0_{\perp, \pm}$, and thus allow us to reuse 
    prior definitions like the correction terms $\gammaod$ in \eqref{eq:gamma1d}.

The edge-local logical space is defined as $\hat V^1(e) = \hat V^1_- \oplus \hat V^1_+$, with first component
\begin{equation}
    \hat V^1_1(e) = \left(\VV^1_+ \otimes \VV^0_{\perp,+}\right) \oplus \left(\VV^1_- \otimes \VV^0_{\perp,-}\right).
\end{equation}
In Definition~\ref{eq:Pe-v1} below, the logical scalar operator $\hat P^1_{\edge}: \hat V^1_1(e) \rightarrow \hat V^1_1(e)$ is defined only on this first, tangential component.
Its vector-valued counterpart $\hat{\bP}_\edge^1$ on the
edge-local space is defined as
\begin{equation} \label{eq:bPe-v1}
\hat{\bP}_\edge^1
    \hat \bu 
    =
    \begin{pmatrix}
        \hat P^1_{\edge} \hat u_1 \\ \hat u_2
    \end{pmatrix},
\end{equation}
for $\hat \bu \in \hat V^1(e)$.
For $\bu \in V^1_\pw$, the physical edge operator $\bP_\edge^1$ is first defined on the two-patch neighborhood by
\begin{equation} \label{eq:bPe-v1-physical}
    \left.\bP_\edge^1 \bu\right|_{\Omega(\edge)}
    =
    \cF^1_\edge \hat{\bP}_\edge^1
    \left((\cF^1_\edge)^{-1}
    \left(\left.\bu\right|_{\Omega(\edge)}\right)\right),
\end{equation}
where $\cF^1_\edge$ denotes the vector-valued pushforward associated with
$F_\edge$, defined patch-wise as in \eqref{pf_gc}. It is then extended to the whole domain by the explicit identity rule 
\begin{equation} \label{eq:bPe-v1-identity}
        \left.\bP_\edge^1 \bu\right|_{\Omega_k}
    = \left.\bu\right|_{\Omega_k}, \quad k \notin \cK(\edge).
\end{equation}
  The construction satisfies the following properties, which are stated under Assumptions~\ref{ass:nbc} and
\ref{as:e_nested}.
  
  \begin{prop}\label{prop:edge-v1}
      Given $\edge \in \edges$ with tangential vector $\btau_\edge$, 
      let $\bP^1_\edge: V^1_\pw \to V^1_\pw$ be 
      the edge-conforming projection \eqref{eq:bPe-v1-physical} associated with
      the logical operator $\hat P^1_{\edge}$ from Definition~\ref{eq:Pe-v1} below,
      with an order of moment preservation $r = r(\bP^1)$ 
      such that $0 \le r \le p$. 
      Then, the following properties hold
      for all $\bu \in V^1_\pw$:
      \begin{enumerate}
          \item $\bP^1_\edge \bu \cdot \btau_\edge$ is continuous across $\edge$.
          \item $\bP^1_\edge \bu = \bu$ if and only if $\bu \cdot \btau_\edge$ is continuous.
          \item $\bP^1_\edge$ preserves polynomial moments of order $r$.
          \item $\bP^1_\edge \bP^1_{\edge'} \bu = \bP^1_{\edge'} \bP^1_\edge \bu $ for any $\edge' \in \edges$.
      \end{enumerate}
  \end{prop}
 The proof of Proposition~\ref{prop:edge-v1} is omitted, as it is analogous to the proof of Proposition~\ref{prop:edge} in Section~\ref{sec:p_edge}.
    This leads to the definition of the conforming projection operator of $V^1_h$. 
  
  \begin{theorem} \label{thm:P1}
    The operator
      \begin{equation} \label{eq:P1}
          \bP^1 = \prod_{\edge \in \edges} \bP^1_\edge: V^1_\pw \to V^1_\pw,
      \end{equation}
      involving the $r$-th order moment-preserving edge-conforming projection as in Proposition~\ref{prop:edge-v1},
      is an $r$-th order moment-preserving projection onto the conforming subspace $V^1_h$.
  \end{theorem}

  \begin{proof}
      Note that \eqref{eq:P1} is well-defined, as the order of the products does not matter 
      thanks to $(4)$ in Proposition~\ref{prop:edge-v1}.
      For the range property of $\bP^1$, let $\bu \in V^1_\pw$:
      $\bP^1 \bu$ is curl-conforming along any edge by $(1)$ in Proposition~\ref{prop:edge-v1}. 
      Thus, $\bP^1$ maps into the curl-conforming space $V^1_h$.
      The projection property directly follows from $(2)$ in Proposition~\ref{prop:edge-v1}.
      The moment preservation follows from $(3)$ in Proposition~\ref{prop:edge-v1}.
  \end{proof}

  \subsection{Edge-based projections}
  \label{sec:p_edge-v1}
  Fix an edge $\edge \in \edges$ and adopt a similar notation to Section~\ref{sec:p_edge} for the derived spaces $\hat V^{1}_\pm$.
The continuity constraint concerns only the component along the tangential vector $\btau_\edge$, 
i.e. in the space $(\hat V^{1}_\pm \cdot \hat\btau_{ \edge})_{\parallel} = \VV^{1}_\pm$.

\subsubsection*{Extension-restriction operators} 
\label{sec:ext-restr-v1}
Since the edge-local derived spaces are still nested, the extension operator on $\VV^{1}_{\pm}$, i.e. 
\begin{equation} \label{eq:extension-v1}
    \scrE^1: \VV^1_{-} \subset \VV^1_{+} \longrightarrow \VV^1_{+}
\end{equation}
is defined analogously to Section~\ref{sec:ext-res}.

As there are no vertex conformity constraints on the restriction operator
\begin{equation} \label{eq:restriction-v1}
    \scrR^1: \VV^1_{+} \longrightarrow \VV^1_{-} \subset \VV^1_{+},
\end{equation} 
it is simply defined as the $L^2$-projection 
from the fine space $\VV^1_{+}$ onto the coarse space $\VV^1_{-}$, i.e. 
\begin{equation} \label{eq:l2-r1}
    \int \lambda^{1, -}_j \scrR^1 \lambda^{1, +}_i = \int\lambda^{1, -}_j \lambda^{1, +}_i, \quad \forall j=0,\dots, n_- -1,i = 0, \dots, n_+ -1.
\end{equation} 
This operator is well-defined and moment-preserving of order $p$, by the same argument that 
$\VV^1_{-}$ contains polynomials up to degree $p-1$. 

The construction of the matrix forms $\arrE^1$ and $\arrR^1$ is discussed in Appendix~\ref{sec:impl_erv0}.

\subsubsection*{Local edge-conforming projection} 
Using the same edge-averaging function, that is defined in the edge-perpendicular direction in the same logical space as in the $V^0$ case, 
namely \eqref{eq:bar-lambda},
 we can define the local edge-based projection without moment preservation as

\begin{align}
    \hat P^1_\edge : \left\{
         \begin{aligned}
             &\lambda^{1, -}_{j_1} \otimes \lambda^-_{\perp, j_2} \mapsto \left\{ \begin{aligned}
                                     & \begin{aligned} & \lambda^{1, -}_{j_1} \otimes \bar\lambda_{\perp}^\edge   \\
                                     \end{aligned}
                                                 &&\text{ if } j_2= 0, \\
                                     &\lambda^{1, -}_{j_1} \otimes \lambda^-_{\perp, j_2}  &&\text{ else.}
                                 \end{aligned}\right. 
              \\
             &\lambda^{1, +}_{j_1} \otimes \lambda^+_{\perp, j_2} \mapsto \left\{ \begin{aligned}
                                     &\begin{aligned} &  \scrR^1 \lambda^{1, +}_{j_1} \otimes \bar\lambda_{\perp}^\edge  \\
                                     \end{aligned}&&\text{ if }  j_2= 0, \\     
                                     &\lambda^{1, +}_{j_1} \otimes \lambda^+_{\perp, j_2} &&\text{ else.}
             \end{aligned}\right.
         \end{aligned} \right. 
 \end{align}

\subsubsection*{Preservation of polynomial moments}  
 
    As in Section~\ref{sec:mom-pres-Pe}, moment preservation for the edge-based projection
    $\hat P^1_\edge$ follows Definition~\ref{mom_pres}. 
    Since $\scrR^1$ is moment-preserving and the correction terms only act in the edge-perpendicular direction,
    which is the same logical space $(\hat V^{1}_\pm \cdot \hat\btau_{ \edge})_{\perp} = \VV^0_\pm$ as in the $V^0$ case,
    we obtain an analogous definition with similar correction coefficients.
    The edge nestedness assumption $\VV^1_- \subset \VV^1_+$ ensures that the corrected functions remain in the edge-local space
    $\hat V^1_{1}(\edge)$ defined in Section~\ref{sec:conf-proj-v1}.
    For instance,
    \begin{equation}
        \lambda^{1,-}_{j_1} \otimes \bar\lambda_{\perp}^\edge
        =
        \frac{1}{2}\left(
            \lambda^{1,-}_{j_1} \otimes \lambda^-_{\perp,0}
            +
            \lambda^{1,-}_{j_1} \otimes \lambda^+_{\perp,0}
        \right)
        \in \hat V^1_{1}(\edge).
    \end{equation}
    \begin{definition}
        \label{eq:Pe-v1}
    The logical scalar edge-conforming projection
    $\hat P^1_\edge: \hat V^1_{1}(\edge) \rightarrow \hat V^1_{1}(\edge)$
    is defined as
    \begin{align}  
        \hat P^1_\edge : \left\{
            \begin{aligned}
                &\lambda^{1, -}_{j_1} \otimes \lambda^-_{\perp, j_2} \mapsto \left\{ \begin{aligned}
                                        & \begin{aligned} & \lambda^{1,-}_{j_1} \otimes \left( \bar\lambda_{\perp}^\edge + \ttC^\edge \right)                  
                                        \end{aligned}
                                                    &&\text{ if }  j_2= 0, \\
                                        &\lambda^{1, -}_{j_1} \otimes \lambda^-_{\perp, j_2} &&\text{ else,}
                                    \end{aligned}\right. 
                 \\
                &\lambda^{1, +}_{j_1} \otimes \lambda^+_{\perp, j_2} \mapsto \left\{ \begin{aligned}
                                        &\begin{aligned} &   \scrR^1 \lambda^{1,+}_{j_1} \otimes  \left( \bar\lambda_{\perp}^\edge - \ttC^\edge \right)   
                                        \end{aligned}&&\text{ if } j_2= 0, \\     
                                        &\lambda^{1, +}_{j_1} \otimes \lambda^+_{\perp, j_2} &&\text{ else,}
                \end{aligned}\right.
            \end{aligned} \right. ,
    \end{align}
    where $\lambda^{1,k}_j$ are the derived scalar patch-wise basis functions from \eqref{eq:basis-hatV-1D}, 
    $\bar\lambda^\edge_\perp$ is the edge-conforming spline in \eqref{eq:bar-lambda},
    $\ttC^\edge$ is the patch-wise correction term in \eqref{eq:Ce}, 
    and $\scrR^1$ is the restriction operator in \eqref{eq:restriction-v1}.
    \end{definition}

\section{Numerical experiments}
\label{sec:num_ex}
This section assesses the accuracy and well-posedness of numerical schemes based on the constructed conforming projection operators and implemented in the IGA library \textsc{psydac} \cite{psydac}. Based on the numerical test cases in Appendix~\ref{sec:disc-examples}, we compare our non-matching multipatch discretization to the matching broken-FEEC framework in \cite[Section 5]{guclu_broken_2023}.

\subsection{Approximation of the weak $\bcurl$ operator}
\label{sec:weak-div-num-ex}
\begin{figure}[!htb]
    \centering
    \includegraphics[width=0.32\columnwidth]{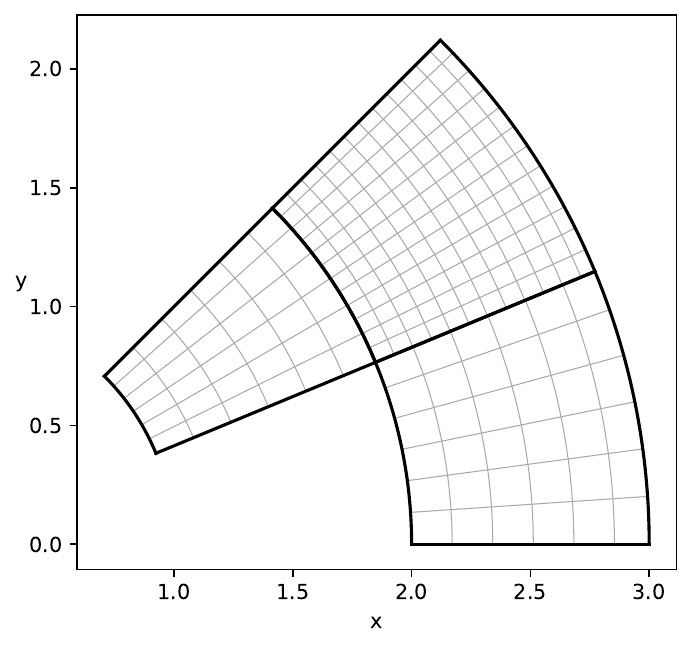}
    \includegraphics[width=0.32\columnwidth]{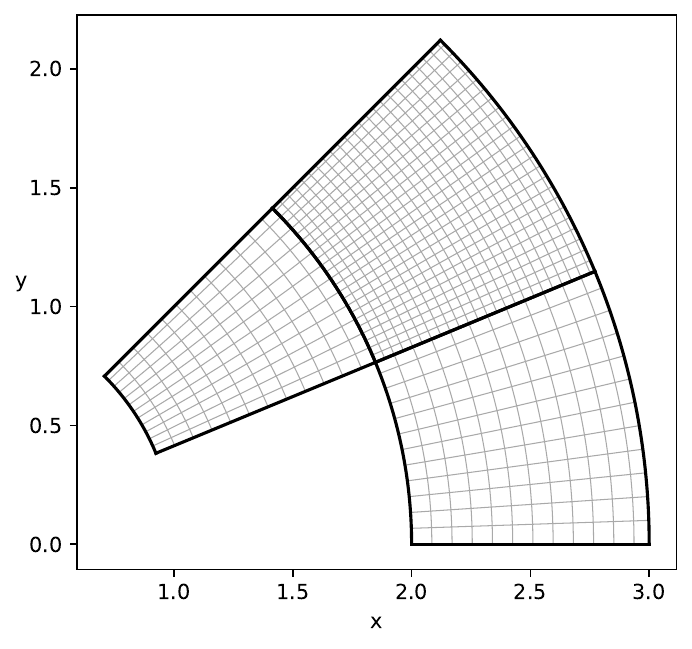}\\
    \includegraphics[width=0.32\columnwidth]{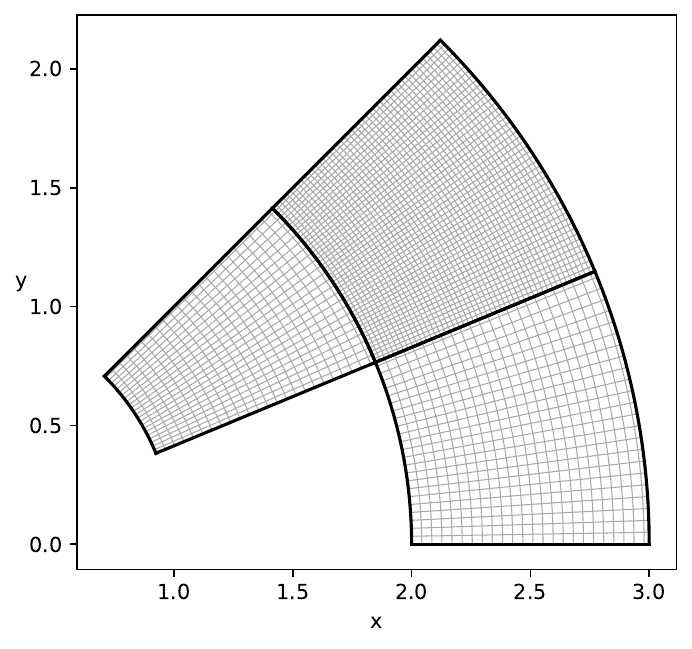}
    \includegraphics[width=0.32\columnwidth]{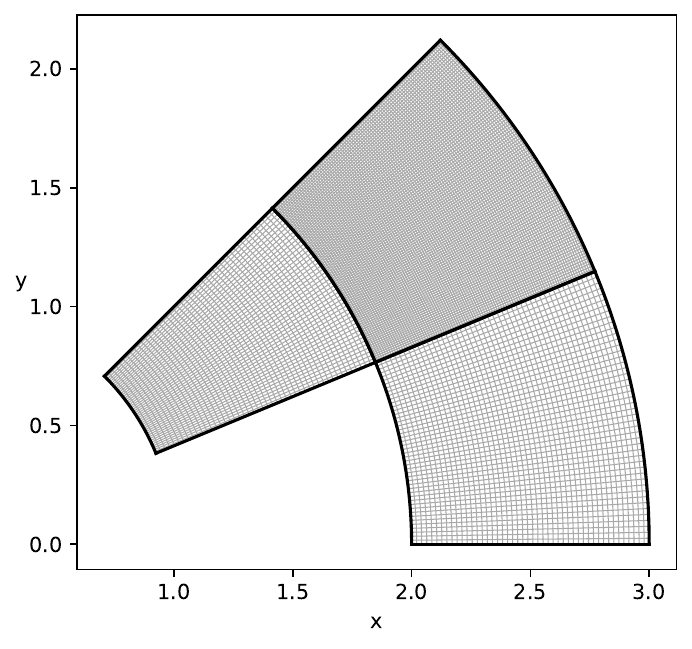}
    \caption{Non-matching refinement pattern of the curved L-shaped domain, where the top-right patch is refined by a factor of two.}
    \label{fig:weak_diff_grids}
\end{figure}
\begin{figure}[!htb]
    \centering
    \includegraphics[width=0.4\columnwidth]{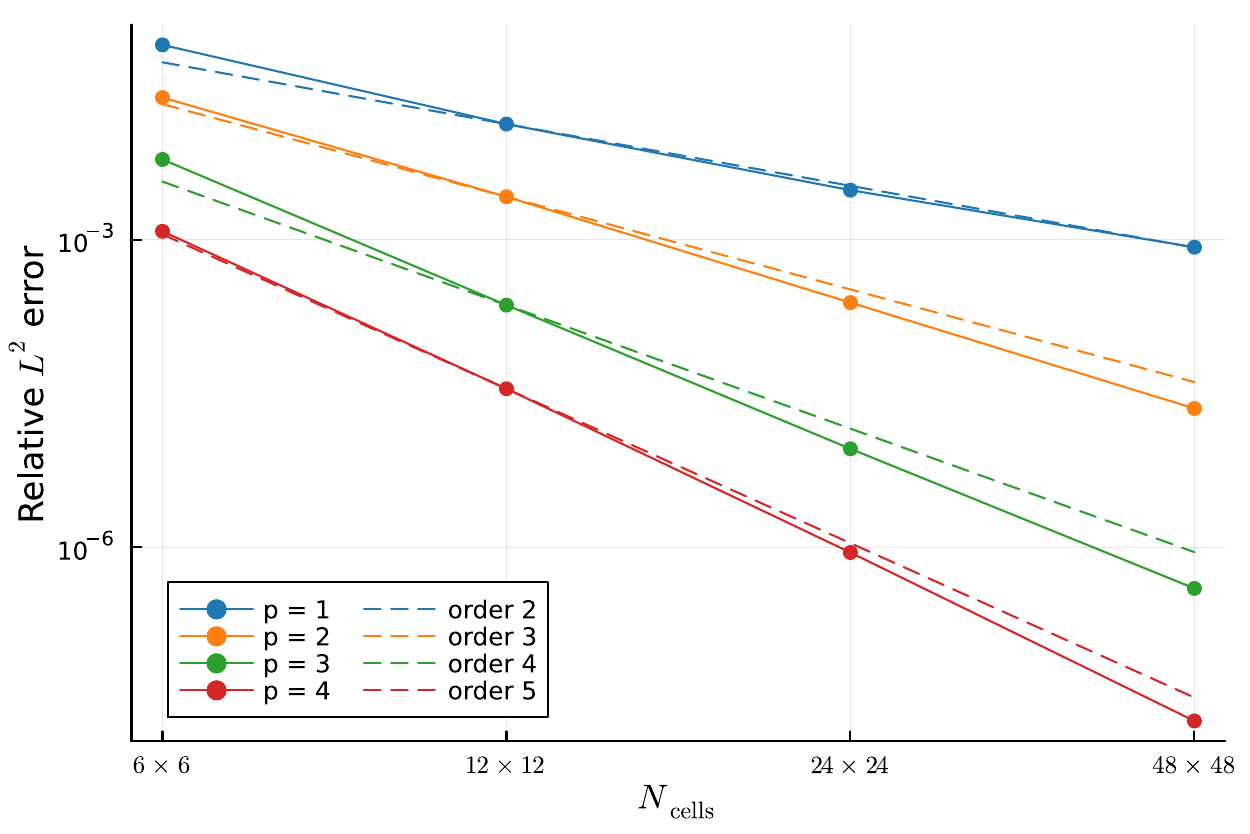} \hspace{1cm}
    \includegraphics[width=0.4\columnwidth]{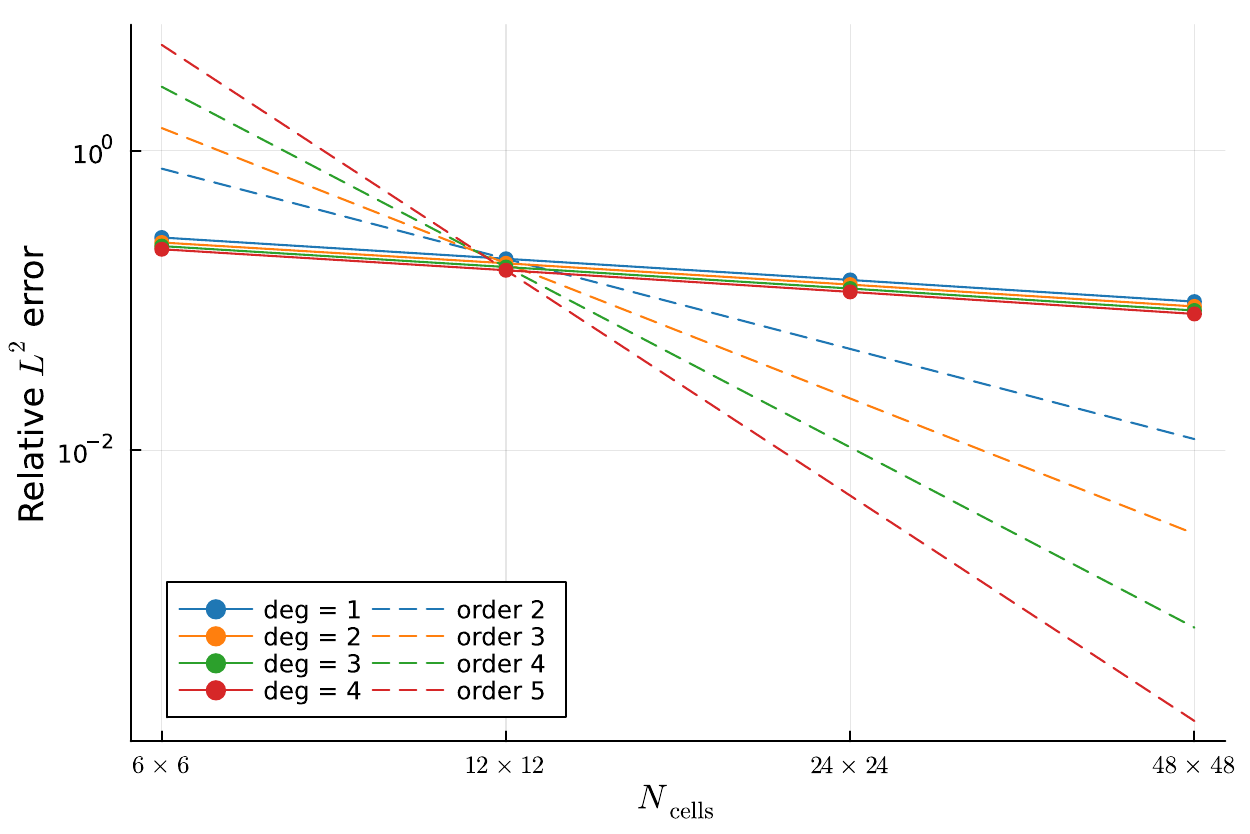}
    \caption{Convergence curves for the 
    weak $\bcurl$ test case with moment preservation of order $p$ (left) and without (right). 
    For different spline degrees $p$, we plot the relative $L^2$-error of the weak $\bcurl$ operator 
    on the refined curved L-shaped domains of Figure~\ref{fig:weak_diff_grids}, where $\nc$ denotes the number of cells on an unrefined patch.}
    \label{fig:weak_curl_conv}
\end{figure}

We apply the weak $\bcurl$ operator from Theorem~\ref{thm:weak-diff_pw}, with its matrix form described 
in Appendix~\ref{subsec:wdiff}, to a 
 smooth manufactured function without boundary conditions:
\begin{align}
    f(x, y) = \cos(\pi x) \cos(\pi y) \not\in H_0(\bcurl; \Omega),
\end{align}
and compute the relative $L^2$-error  $\| \widetilde \bcurl_\pw  \tilde \Pi_{\pw, 0}^2 f - \bcurl f\|_{L^2} / \| \bcurl f\|_{L^2}$.
In Figure~\ref{fig:weak_curl_conv}, we show the convergence curves for the weak $\bcurl$ operator.
We observe the expected convergence order when the weak differential operator preserves the required polynomial moments. Without moment preservation, the convergence deteriorates, demonstrating the importance of this property as is discussed in Section~\ref{sec:mom-pres}.


\subsection{Poisson problem} 
\begin{figure}[!htb]
        \centering
        \includegraphics[width=0.29\columnwidth]{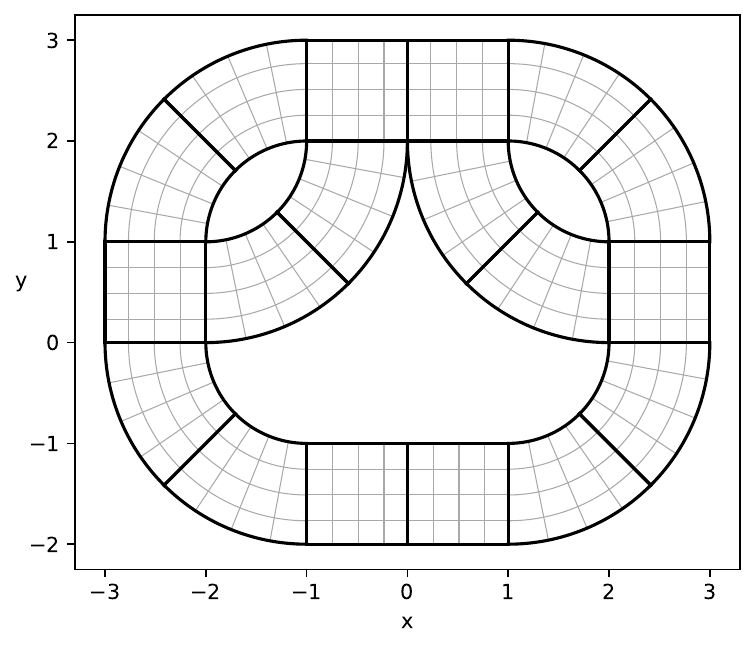}
        \includegraphics[width=0.29\columnwidth]{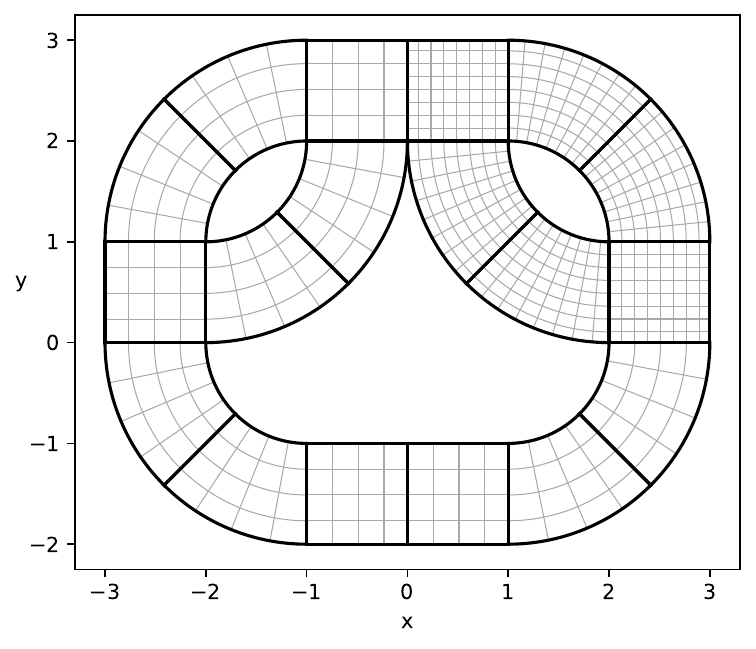}
        \includegraphics[width=0.29\columnwidth]{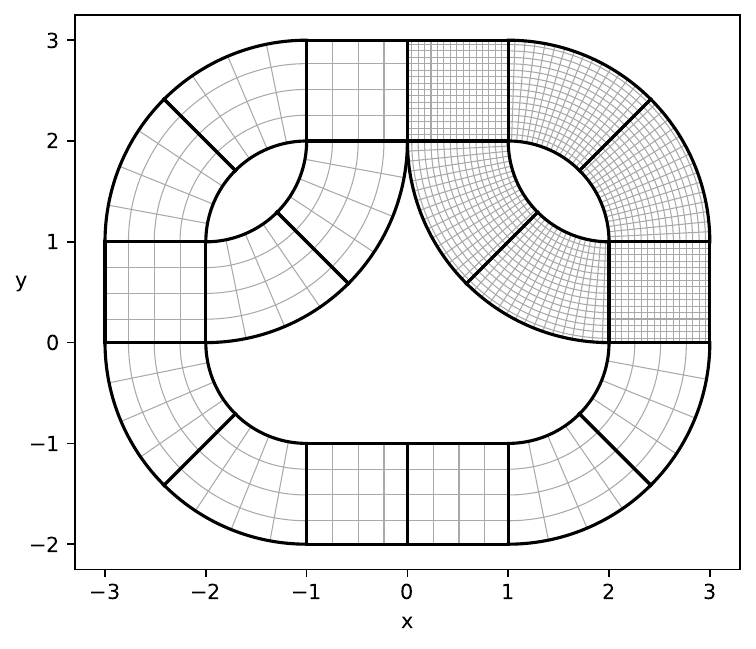}
        \caption{The non-matching refinement strategy of the pretzel domain.}
        \label{fig:poisson_refinement}
    \end{figure}

    For the Poisson problem of Appendix~\ref{subsec:poisson}, the experiments use the
manufactured solution and source from \cite[Section 5.1]{guclu_broken_2023}.
Both are localized near the upper-right handle of the pretzel-shaped domain.
We therefore compare uniform refinement of the full domain with the
non-matching refinement strategy of Figure~\ref{fig:poisson_refinement},
which refines only this handle.

\begin{figure}[!ht]
        \centering
        \includegraphics[width=0.3\columnwidth, trim={0 4cm 0 0}]{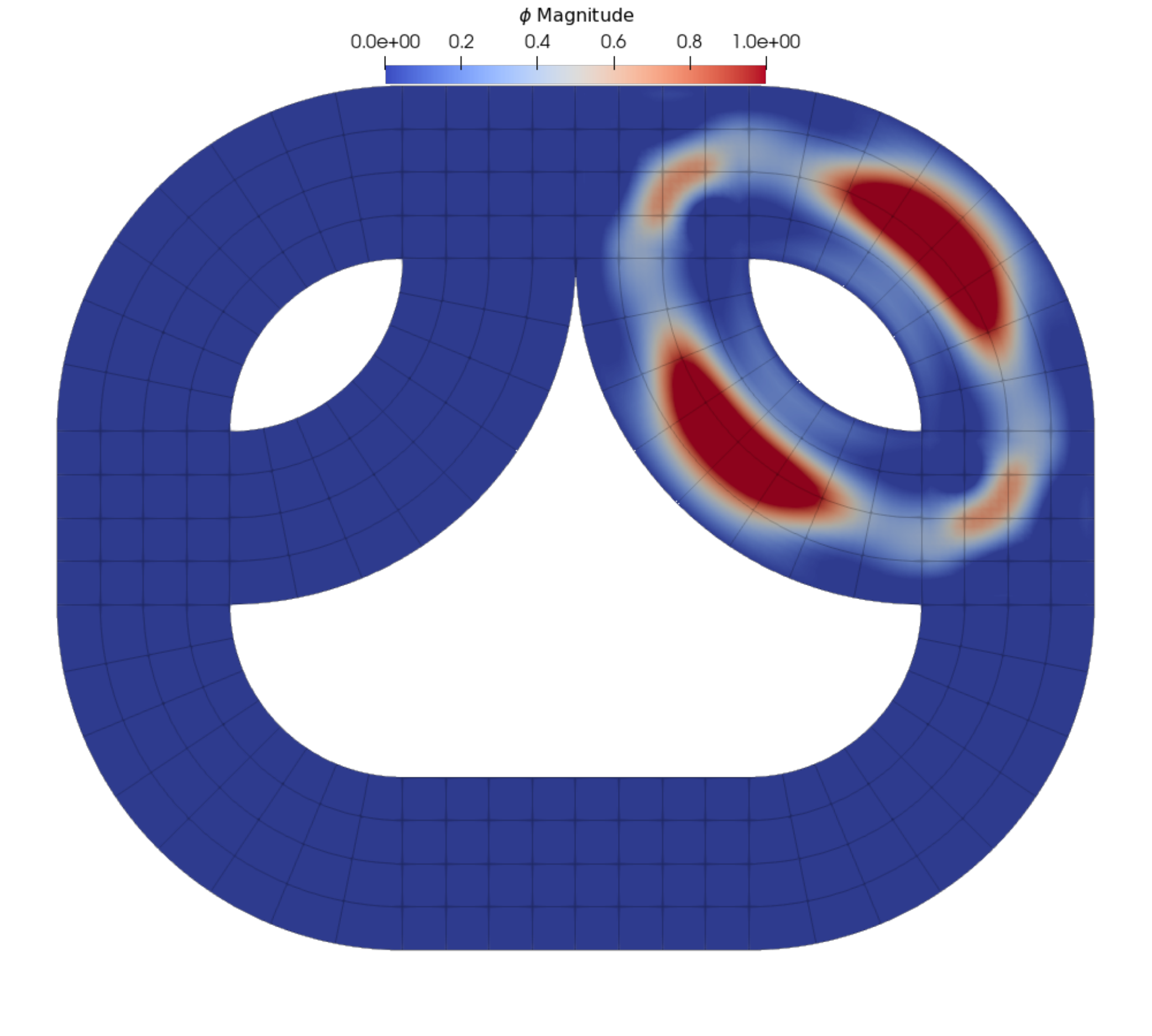}
        \includegraphics[width=0.3\columnwidth, trim={0 4cm 0 0}]{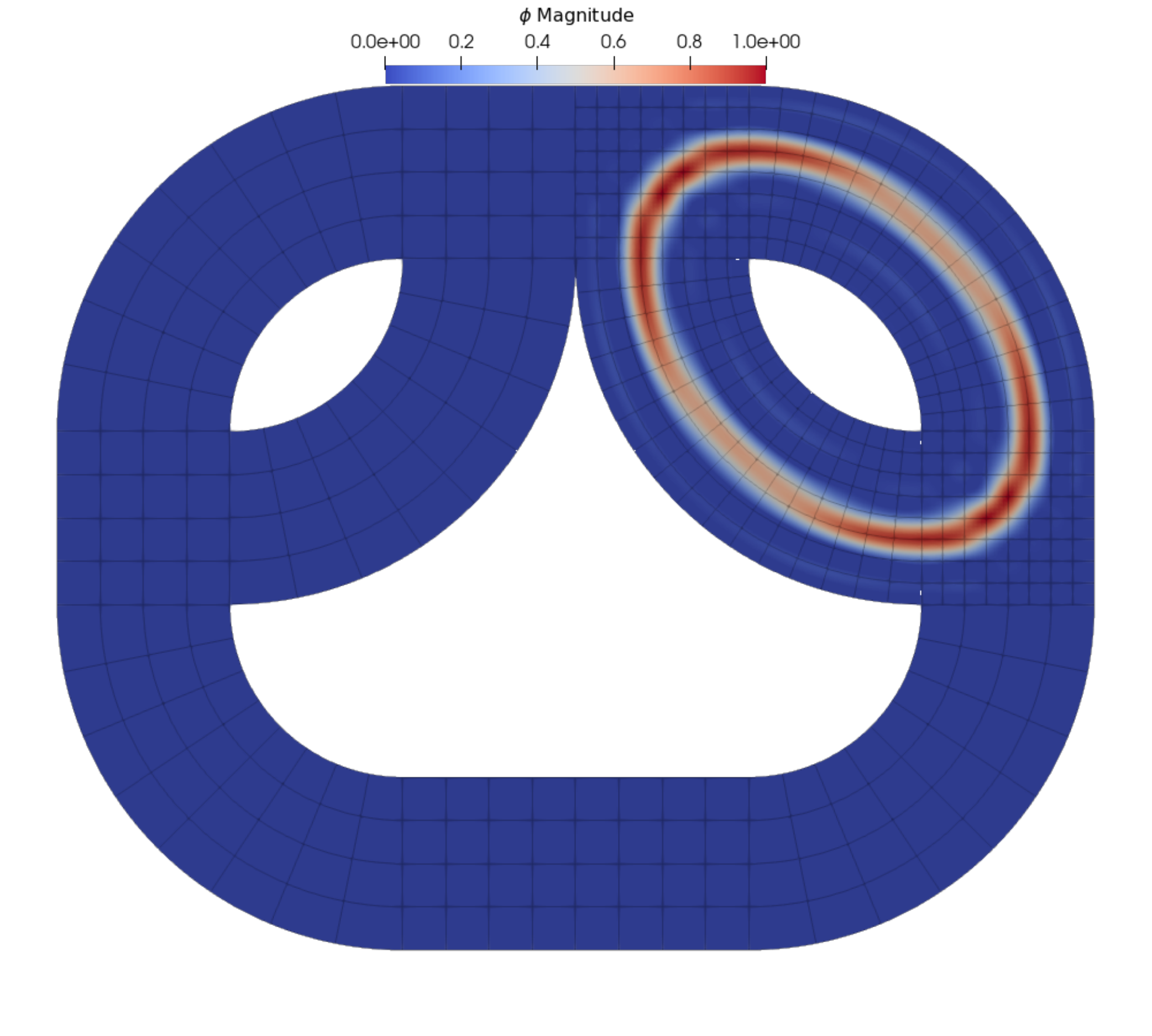}
        \includegraphics[width=0.3\columnwidth, trim={0 4cm 0 0}]{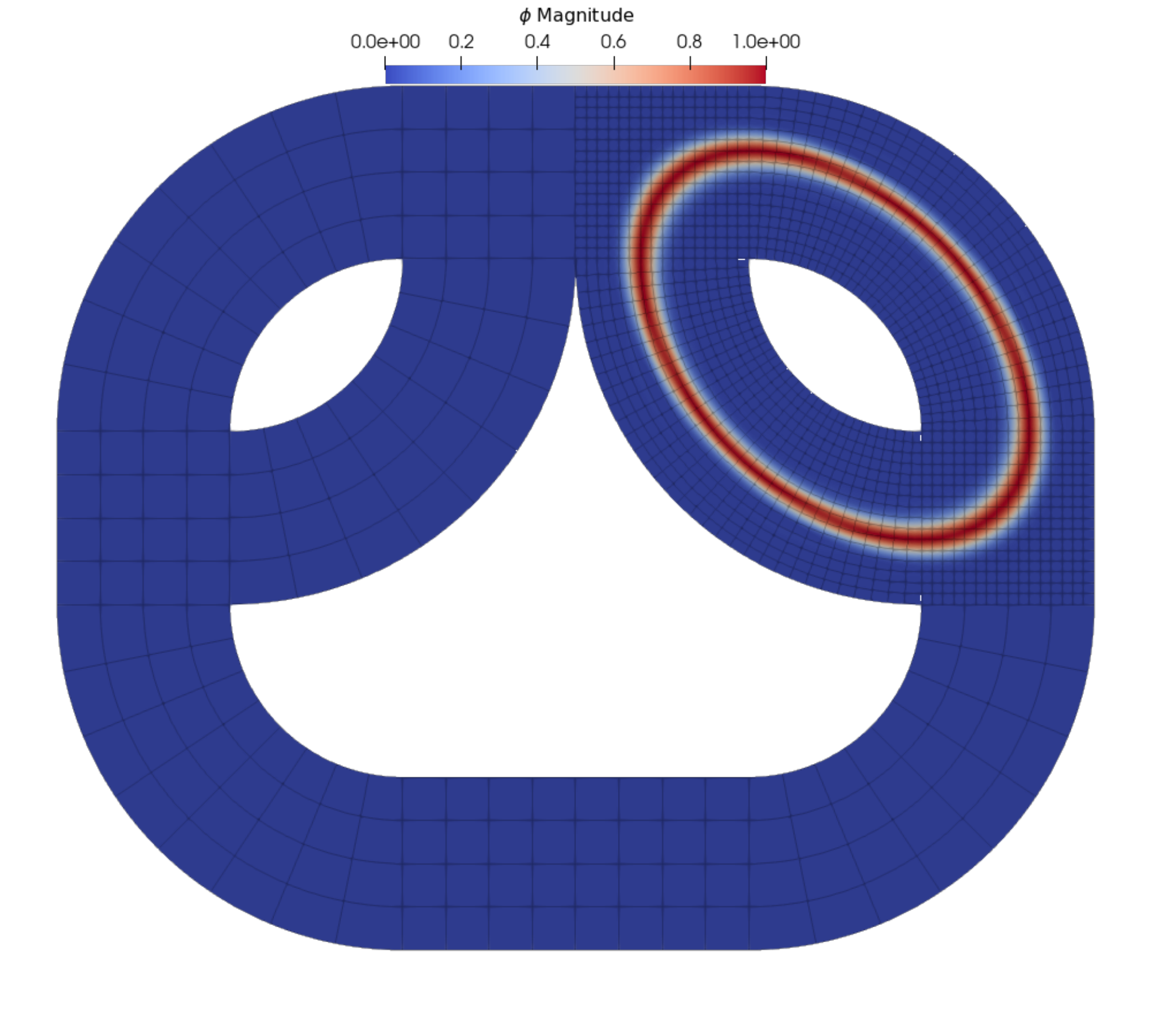}
        \caption{The numerical solution of the Poisson problem for spline degree $3$ on the refined domains of Figure~\ref{fig:poisson_refinement}.}
        \label{fig:poisson_sol_refinement}
    \end{figure}

Figure~\ref{fig:poisson_sol_refinement} shows the numerical solution for
\(p=3\) at the three refinement levels. The solution is negligible away from
the refined handle, consistent with the localization of the manufactured
solution.

Table~\ref{tab:poisson_errors} compares the relative \(L^2\)-errors and the
total number of degrees of freedom for both refinement strategies. At each
refinement level, the two discretizations attain essentially identical
errors. Once the refinement is localized, the refined-handle
discretization achieves this accuracy with substantially fewer degrees of
freedom. This demonstrates the effectiveness of the non-matching refinement strategy, which is enabled by the discretization developed in this work.

 \begin{table}[ht]
  \centering
    \begin{tabular}{|cc|cc|cc|}
    \hline
    & & \multicolumn{2}{c|}{uniform refinement} & \multicolumn{2}{c|}{refined handle} \\
    \cline{3-6}
    $\nc$  & $\bp$ & $\ndof$ & relative error & $\ndof$ & relative $L^2$-error \\
    \hline
     \multirow{3}{*}{$4 \times 4$} & $3 \times 3$ &  882 & 0.8652 & 882 & 0.8652 \\
                        & $4 \times 4$ & 1152 & 0.3446 & 1152 & 0.3446 \\
                        & $5 \times 5$ & 1458 & 0.1831 & 1458 & 0.1831 \\
    \hline
     \multirow{3}{*}{$8 \times 8$} & $3 \times 3$ &  2178 & 0.0460 & 1314 & 0.0460 \\
                        & $4 \times 4$ & 2592 & 0.0867 & 1632 & 0.0867 \\
                        & $5 \times 5$ & 3042 & 0.0556 & 1986 & 0.0557 \\
    \hline
    \multirow{3}{*}{$16 \times 16$} & $3 \times 3$ & 6498 & 0.0048 & 2754 & 0.0048 \\
                        & $4 \times 4$ & 7200 & 0.0057 & 3168 & 0.0057 \\
                        & $5 \times 5$ & 7938 & 0.0036 & 3618 & 0.0036 \\
    \hline
  \end{tabular}

  \caption{Relative $L^2$-errors for the Poisson problem. In the first column $\nc$ denotes
  the number of cells per parametric direction on a refined patch, where $4 \times 4$ cells is the resolution of the unrefined patches, $\bp$ is the spline degree per parametric direction, and
  $\ndof$ is the total number of degrees of freedom.
  The uniform case refines all patches; the refined handle case
  refines only the top-right handle as shown in
  Figure~\ref{fig:poisson_refinement}.}
  \label{tab:poisson_errors}

\end{table}

\subsection{Time-harmonic Maxwell problem} 
\begin{figure}[!ht]
        \centering
        \includegraphics[width=0.26\columnwidth]{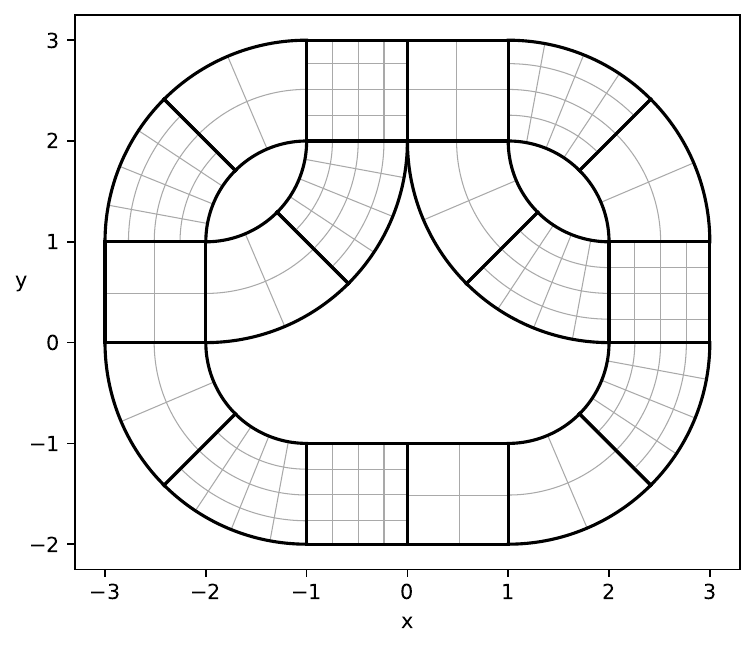}
        \includegraphics[width=0.26\columnwidth]{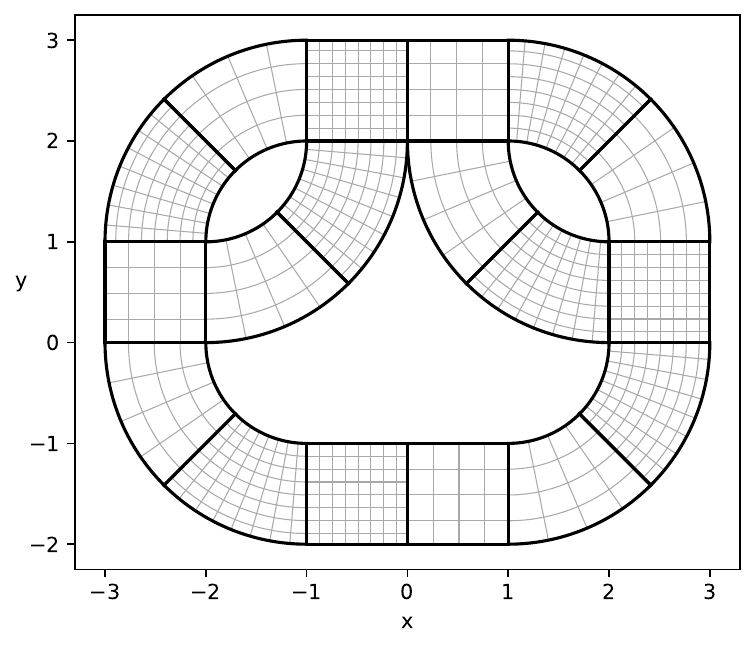}\\
        \includegraphics[width=0.26\columnwidth]{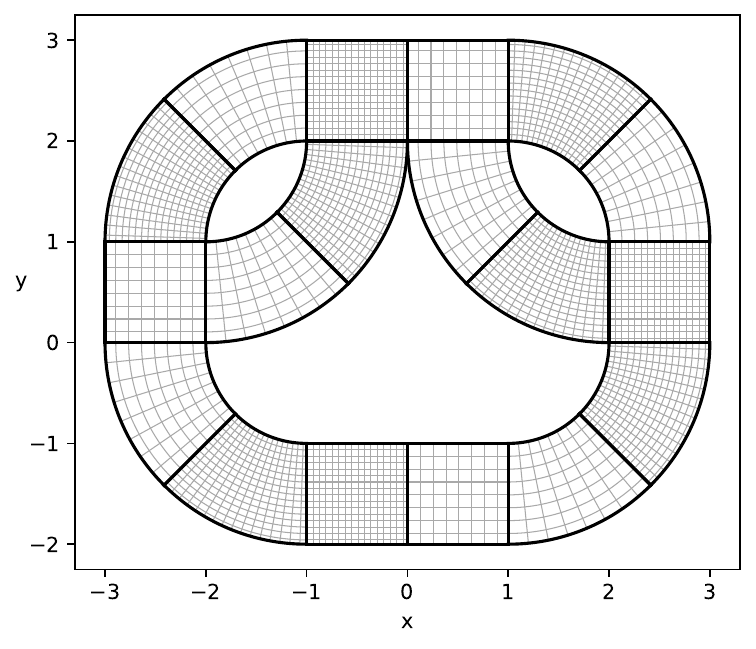}
        \includegraphics[width=0.26\columnwidth]{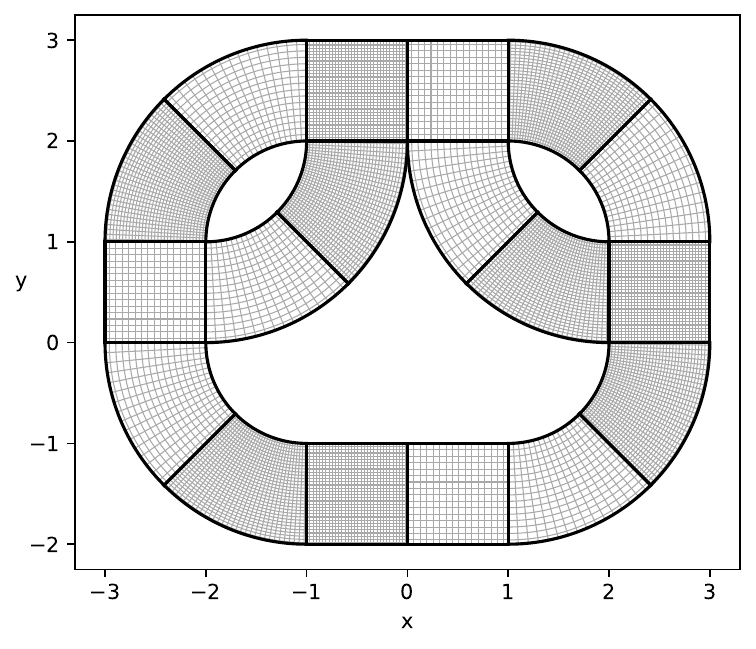}
        \caption{Refinement pattern of the pretzel-shaped
        domain.}
        \label{fig:th-maxwell-grid}
    \end{figure}

For the time-harmonic Maxwell problem of Appendix~\ref{subsec:thmaxwell}, the manufactured solution and source are taken from
\cite[Section~5.4]{guclu_broken_2023}. We compare uniform
discretizations with non-matching discretizations, where, starting from the uniformly refined grid, half of the patches are selected at
random and additionally refined by a factor of two. The resulting refinement
pattern is shown in
Figure~\ref{fig:th-maxwell-grid}.

\begin{figure}[!htb]
        \centering
        \includegraphics[width=0.45\columnwidth, trim={0 0cm 0 0}]{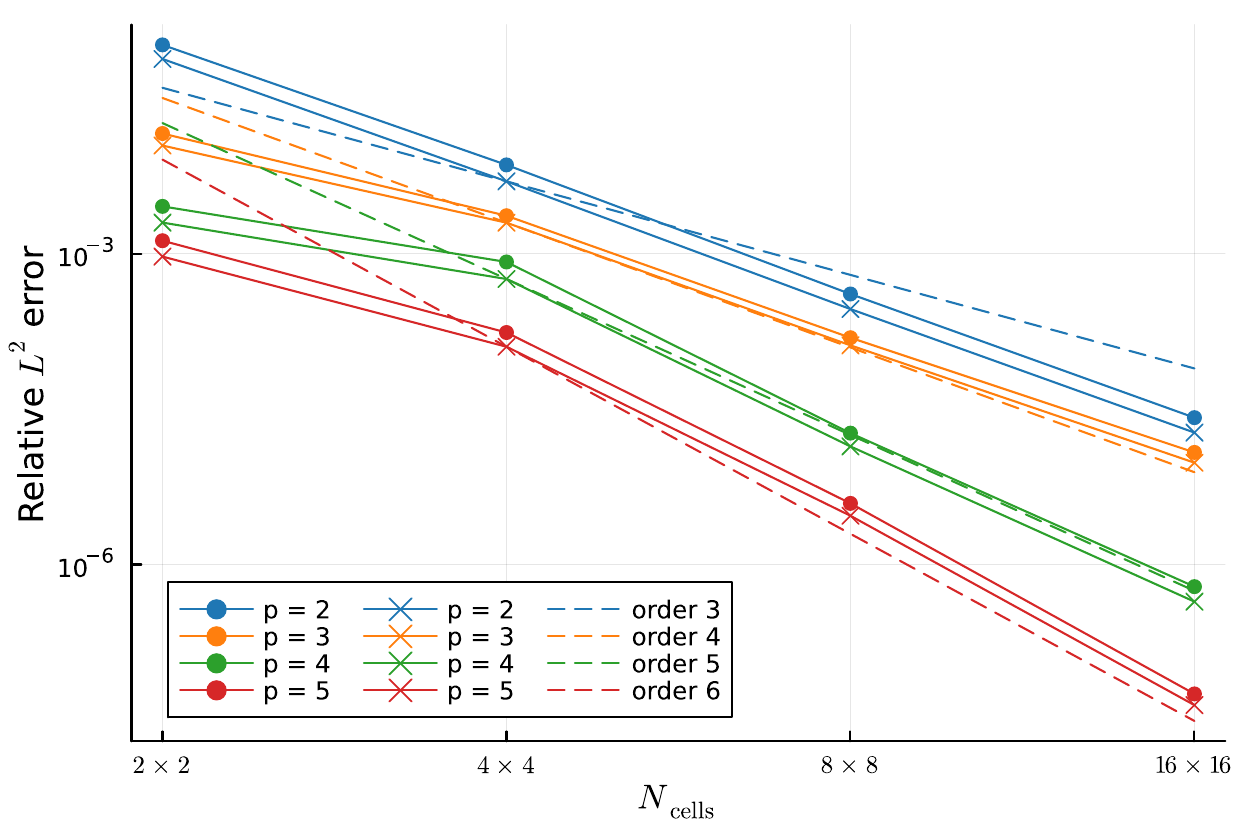}
        \hspace{1cm}
        \includegraphics[width=0.35\columnwidth, trim={0 0cm 0 0}]{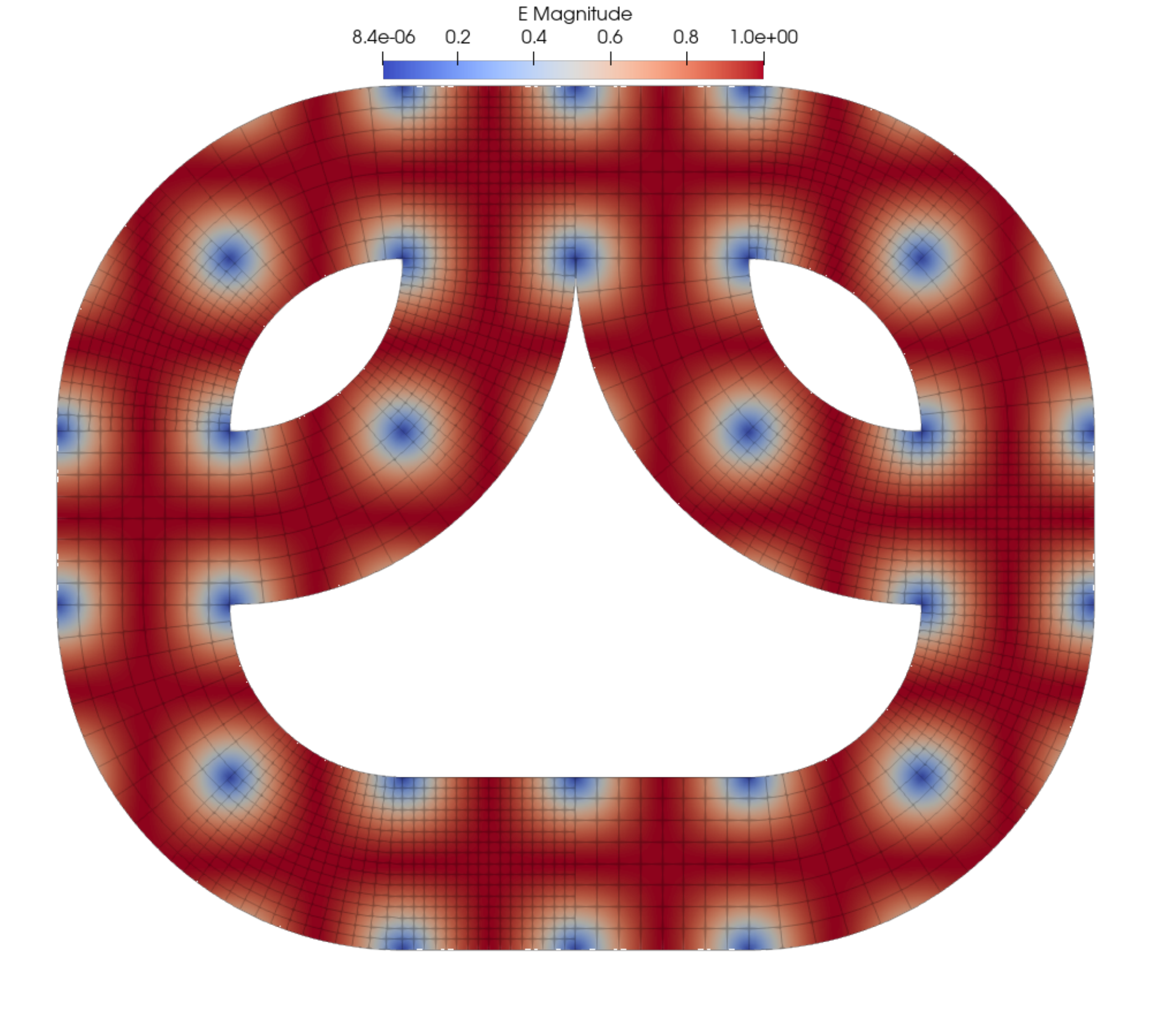}
        \caption{On the left, we plot the convergence curves for uniform (circle) and  refined (cross) 
        grids from Figure~\ref{fig:th-maxwell-grid}, 
        where the abscissa $\nc$ describes the unrefined number of cells on a patch.
        On the right, we show the amplitude of the numerical solution $\boldsymbol{E}_h$ 
       for spline degree $p=3$ on the refined domain in the bottom-left of Figure~\ref{fig:th-maxwell-grid}.}
        \label{fig:th-maxwell-conv}
    \end{figure}

\clearpage
    In Figure~\ref{fig:th-maxwell-conv}, the convergence rates for both cases are shown. They are qualitatively the same and coincide with the expected convergence rates of the method, i.e. they follow the reference orders $p+1$, so the patch-wise refinement does not reduce the asymptotic accuracy.
    The magnitude of the numerical
    solution on a non-matching grid is also shown in Figure~\ref{fig:th-maxwell-conv}.

    \subsection{\texorpdfstring{$\bcurl$-$\curl$}{curl-curl} eigenvalue problem} 
\begin{figure}[!ht]
    \centering
    \includegraphics[width=0.27\columnwidth]{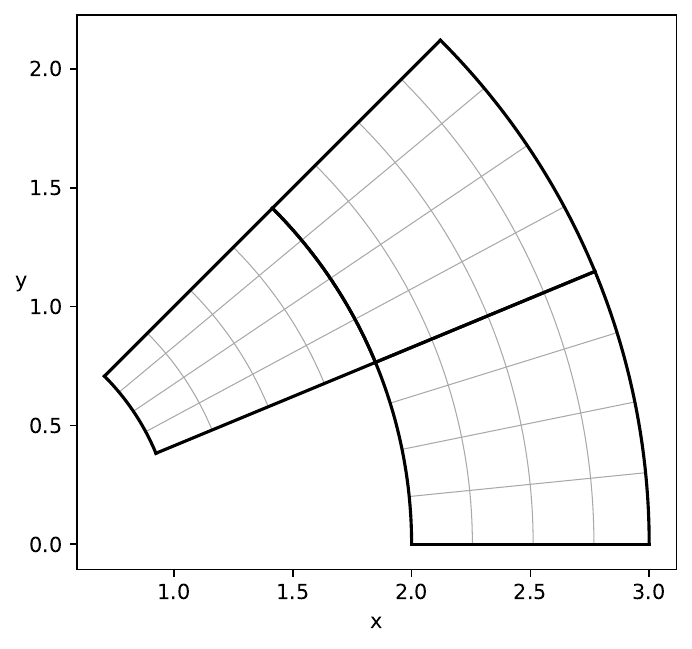}
    \includegraphics[width=0.27\columnwidth]{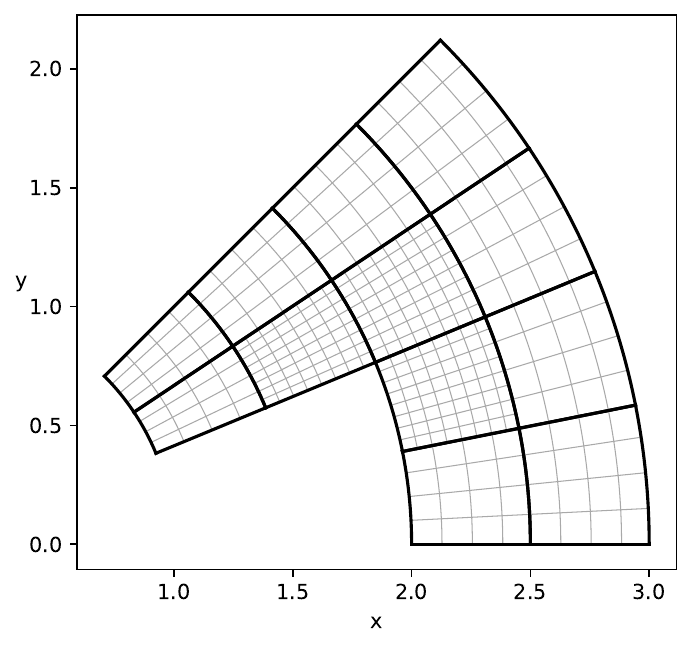}\\
    \includegraphics[width=0.27\columnwidth]{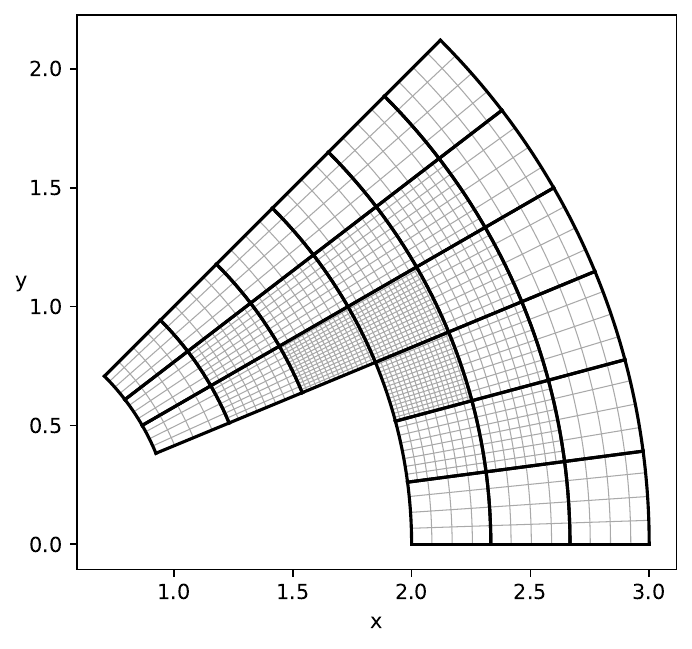}
    \includegraphics[width=0.27\columnwidth]{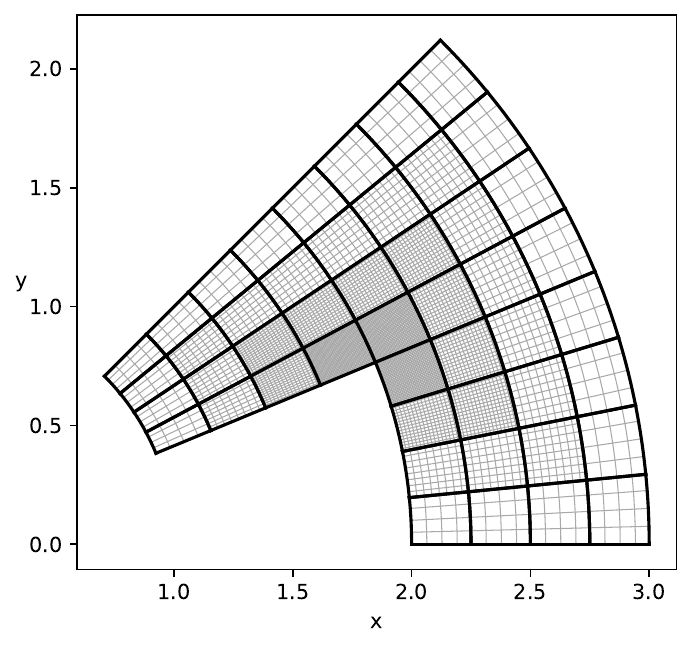}
    \caption{Refinement of the reentrant corner by adding patches of higher resolution.}
    \label{fig:cc-ev-grid}
\end{figure}

    The curved L-shaped domain is refined by layers of finer patches around the reentrant corner, 
    as shown in Figure~\ref{fig:cc-ev-grid}.

Solving the $\bcurl$-$\curl$ eigenvalue problem from Appendix~\ref{subsec:ccev} on this domain, 
we obtain the first five eigenvalues and compare them to a matching discretization
 with the same patch layout, where patches have been uniformly refined to obtain a comparable number of total degrees of freedom.

    \begin{figure}[!htb]
        \centering
        \includegraphics[width=0.45\columnwidth]{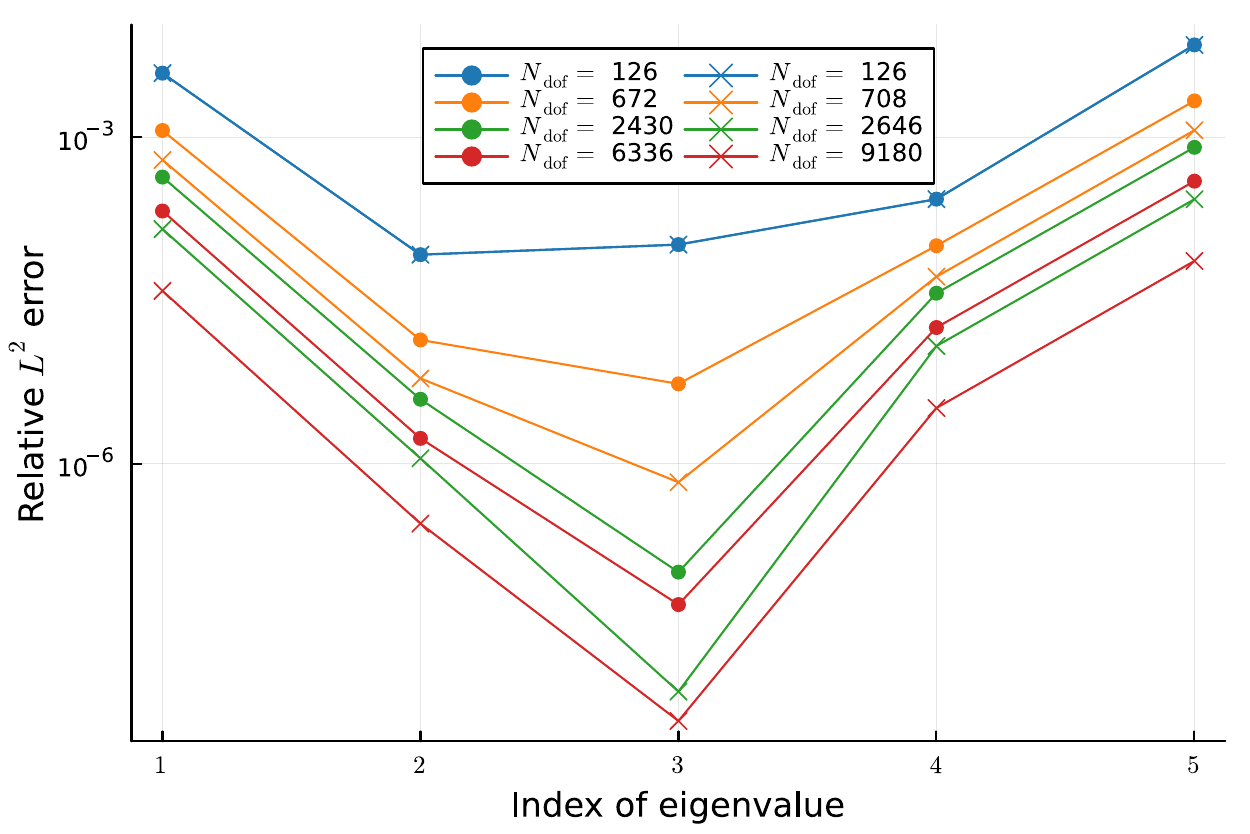}
        \hspace{1cm}
        \includegraphics[width=0.45\columnwidth]{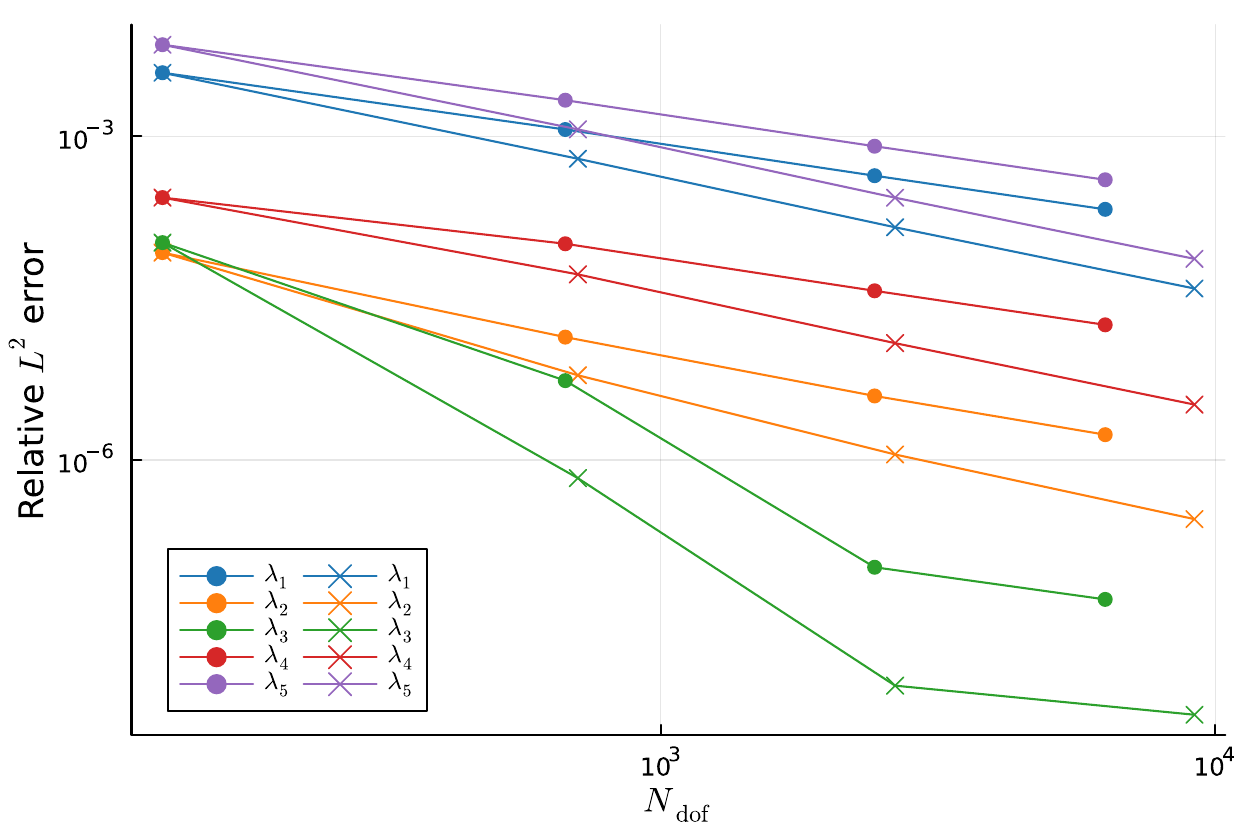}
        \caption{On the left, we plot the error of the first five eigenvalues for matching (circle) and non-matching (cross) domains
        computed with spline spaces of degree $p=3$.
        We use the refinement strategy of Figure~\ref{fig:cc-ev-grid} and compare the results to a matching discretization with a comparable
        number of DOF. 
        On the right, we show the corresponding convergence curves for the individual eigenvalues. 
        }
        \label{fig:cc-ev-conv}
       \end{figure}

    The results can be seen in 
    Figure~\ref{fig:cc-ev-conv}, where we plot the error for each eigenvalue for both 
    discretizations.
    Both discretizations are spectrally correct, and the non-matching refinement strategy shows a better approximation of the eigenvalues, which is expected due to the refinement of the reentrant corner.
     This case is particularly interesting, since the eigenfunctions have a singularity at the reentrant corner, which is resolved better by the non-matching refinement strategy.
    Additionally, in Figure~\ref{fig:cc-ev-conv}, we show the individual convergence curves for the first five eigenvalues, where we observe improved convergence rates for the non-matching refinement strategy. Surprisingly, the third eigenvalue, which is not associated with a singular eigenfunction, shows an improved rate of convergence as well. 

    \begin{figure}[!ht]
        \centering
        \includegraphics[width=0.35\columnwidth]{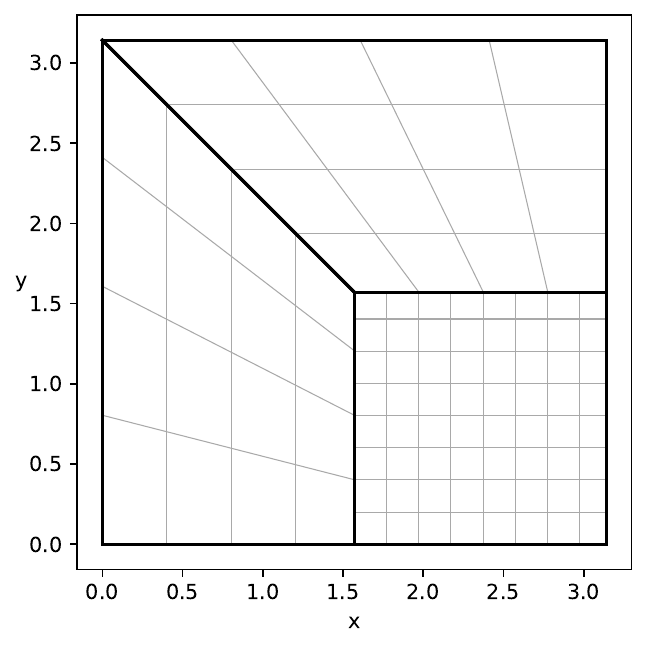}
        \hspace{1cm}
        \includegraphics[width=0.45\columnwidth]{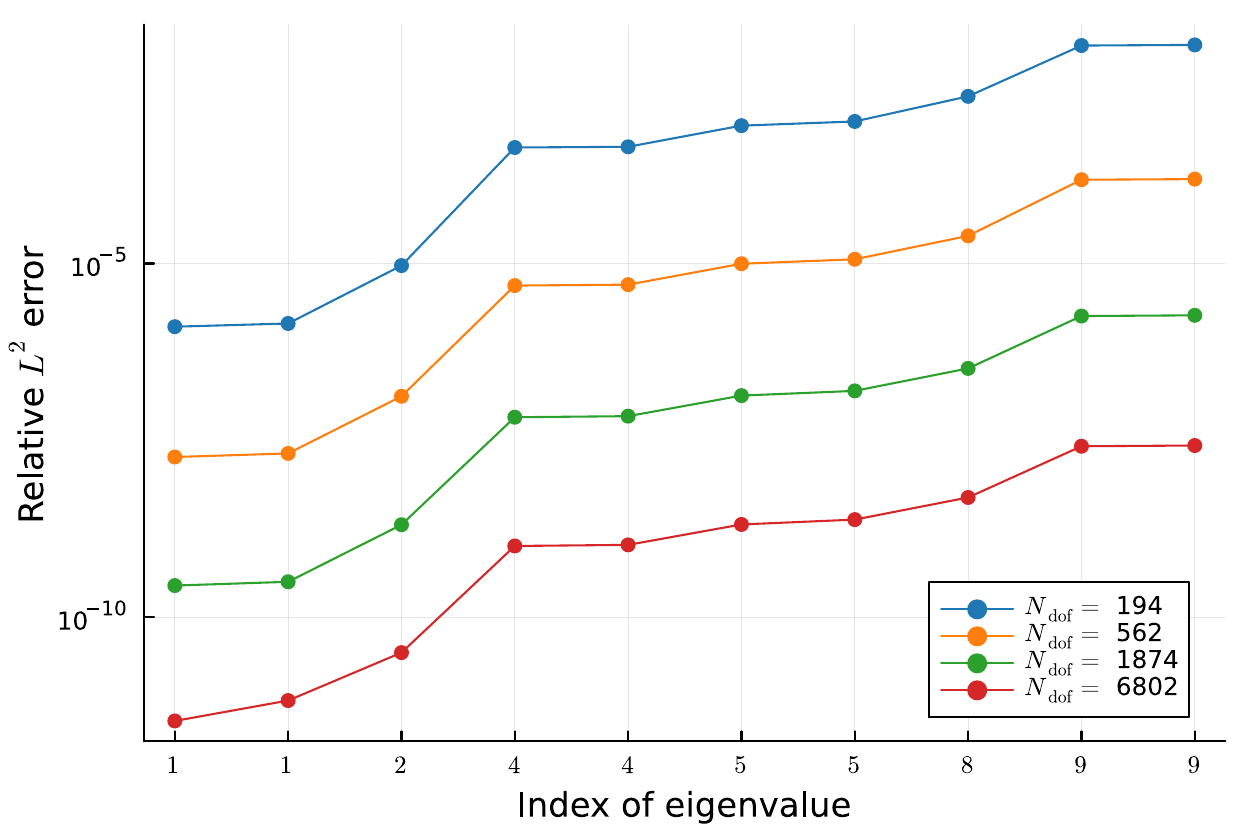}
        \caption{The $\bcurl$-$\curl$ eigenvalue problem on a three-patch domain, shown left, and the error of the first ten eigenvalues on the right, where we continuously refine the domain on the left by a factor of two.}
    \label{fig:three-patch-domain}
    \end{figure}

    In \cite{BCP},
    where we showed the existence of bounded commuting projections for these non-matching multipatch domains,
    we had limitations when it comes to certain grid configurations or refinements. For example, 
    the theory does not yet cover the case of domains with three patches around a vertex, but 
    we did not observe any issues in our numerical experiments.
    In particular in Figure~\ref{fig:three-patch-domain}, where we show such a three-patch domain with different resolutions and the errors of the first ten eigenvalues of the $\bcurl$-$\curl$ eigenvalue problem, we can observe that the method is still spectrally correct and converges as expected.

    \subsection{Time-domain Maxwell} \label{sec:td_maxwell_broken}
    Finally, consider the Maxwell system from Appendix~\ref{subsec:tdmaxwell} 
    on a square domain $\Omega = [0, \pi] \times [0, \pi]$ with absorbing boundary conditions.
    As an initial condition, we set:
    \begin{equation*}
        \bE(t = 0, \bx) = \frac{1}{\sigma^2}\begin{pmatrix}
        x_2  - \bar x_2\\ -( x_1 - \bar x_1)
        \end{pmatrix} e^\frac{- \|\bx - \bar \bx\|_2^2}{2 \sigma^2}, \quad \bB(t = 0, \bx) = \curl \bE, \quad \bJ(t, \bx) = 0,
    \end{equation*}
    where $\bar \bx = (0.5 \pi , 0.5 \pi )$ and $\sigma = 0.1$.
	For the time discretization, the semi-discretization \eqref{eq:tdmaxwell-semdisc} is cast into a simple leapfrog time stepping scheme for some time step $\Delta t > 0$, which reads
    \begin{equation} \label{eq:schemes_td_maxwell}
	\begin{aligned}
		\arrb^{n + \frac 12} &= \arrb^n - \frac{\Delta t}{2} \arrC \arrP^1 \arre^n, \\
		\arrM^1 \arre^{n+1} &= \arrM^1 \arre^n + \Delta t   (\arrC \arrP^1) ^T \arrM^2 \arrb^{n+\frac12} \\
		\arrb^{n+1} &= \arrb^{n+\frac12} - \frac{\Delta t}{2} \arrC  \arrP^1 \arre^{n+1}.
	\end{aligned}
\end{equation}

    \begin{figure}[!htb]
        \centering
        \includegraphics[width=0.3\columnwidth]{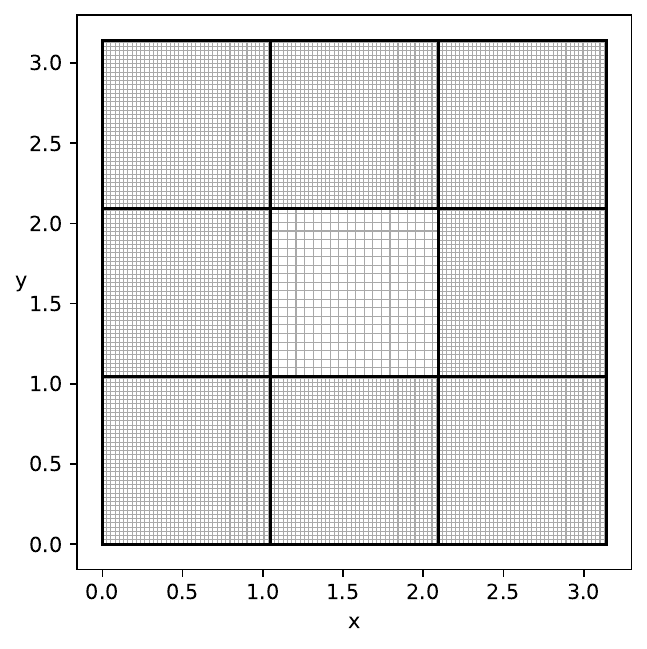}
        \includegraphics[width=0.3\columnwidth]{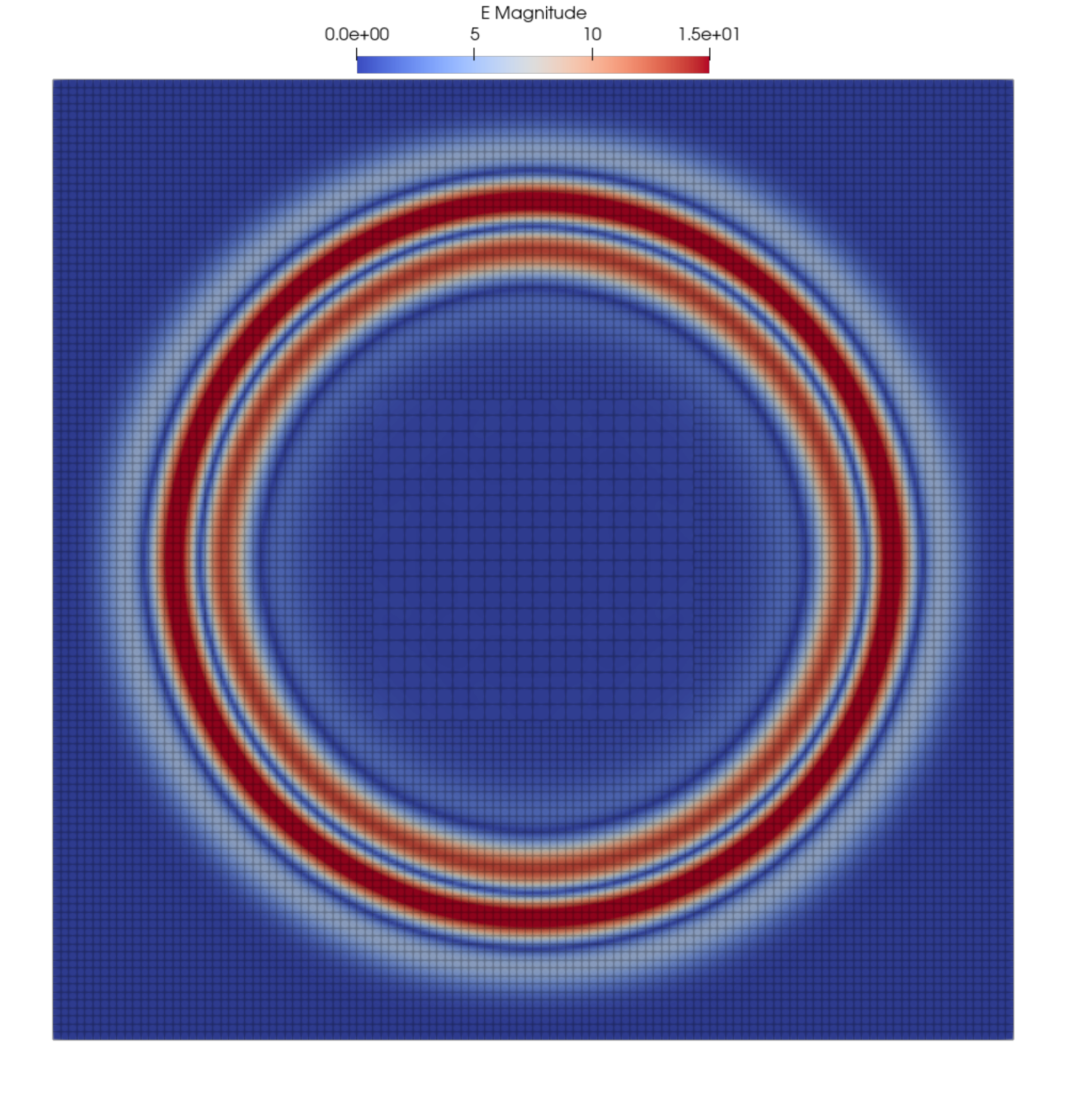}
        \includegraphics[width=0.3\columnwidth]{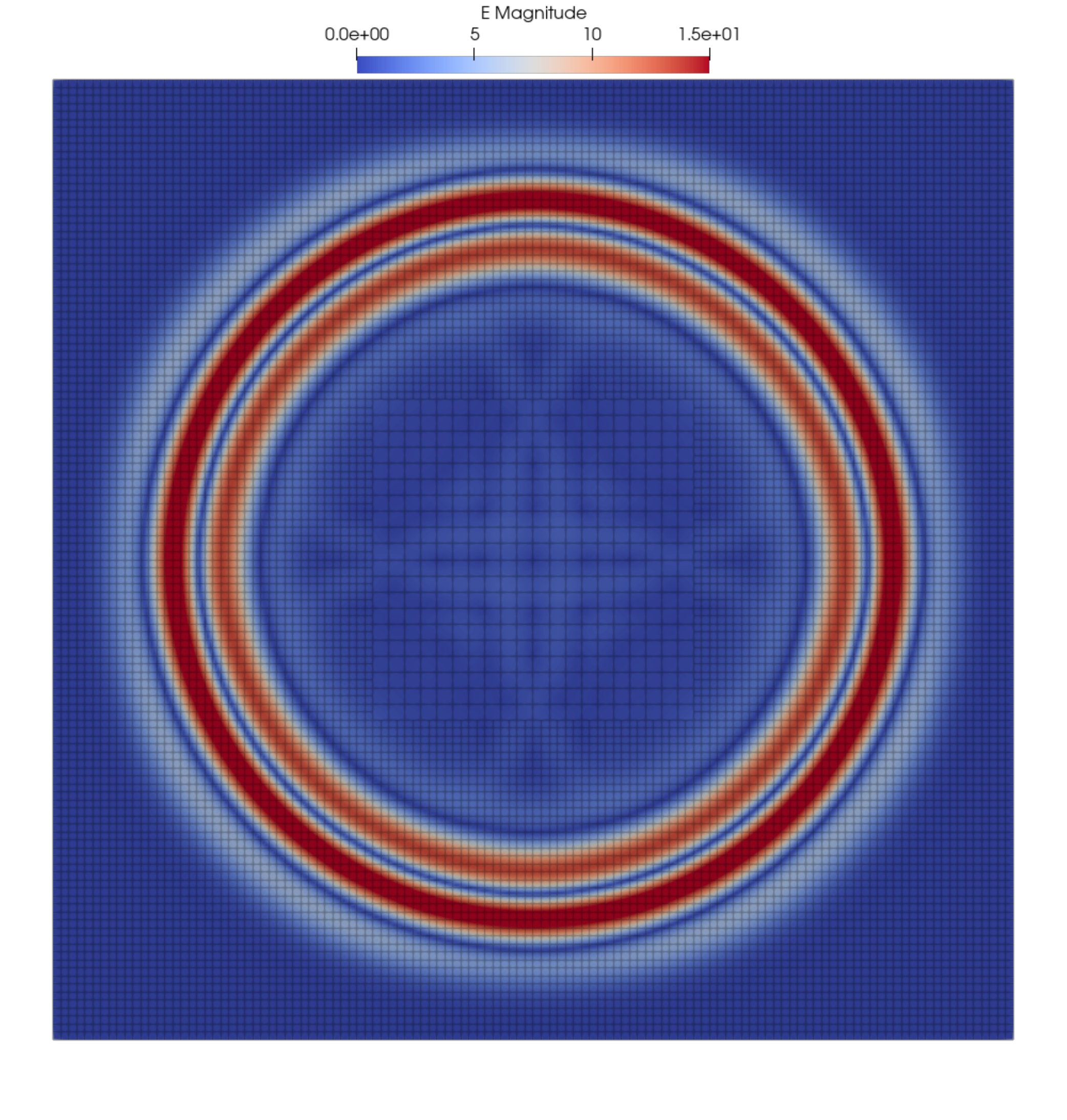} \\
        \includegraphics[width=0.3\columnwidth]{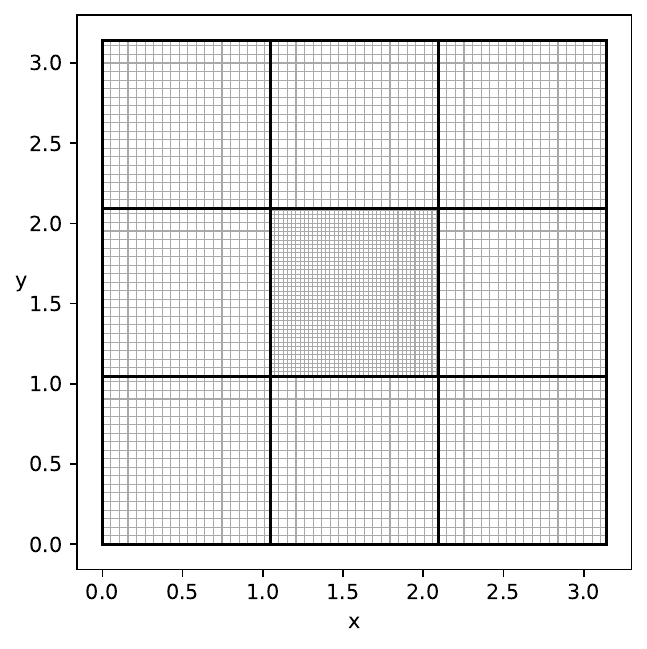}
        \includegraphics[width=0.3\columnwidth]{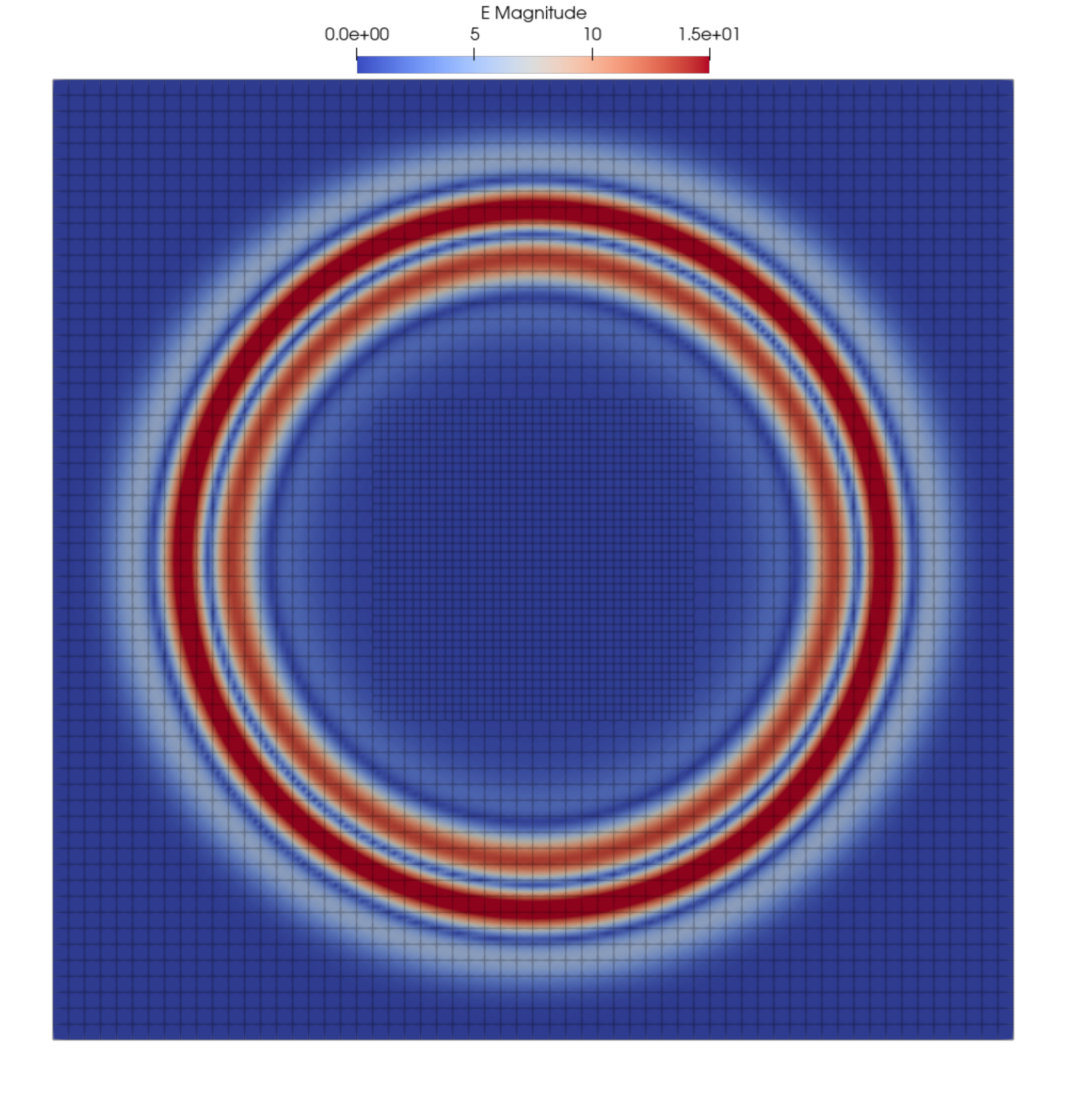}
        \includegraphics[width=0.3\columnwidth]{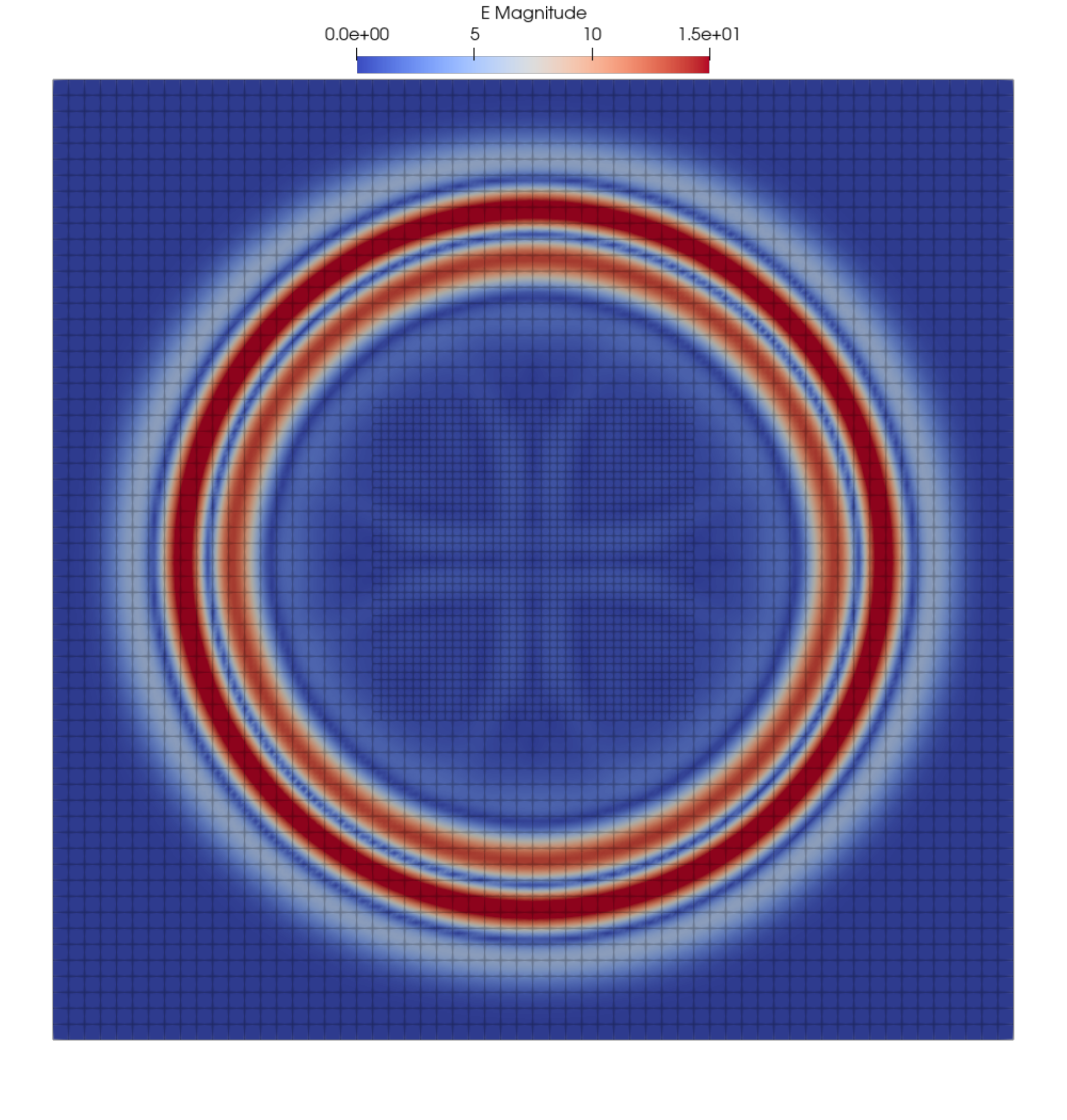} 
        \caption{The grid setup (left) and the amplitude of the electric field $\bE_h$ of the time-domain Maxwell problem after leaving the interior patch. 
        The discretization in the middle involves a conforming projection 
        that preserves polynomial moments up to the spline degree,
         whereas the discretization on the right does not.}   
             \label{fig:td-maxwell-comp}
    \end{figure}

    For our test cases, we consider two different patch configurations:
    a coarse patch enclosed with fine patches and a fine patch enclosed with coarse patches 
     as is shown in the left column of Figure~\ref{fig:td-maxwell-comp}.
     
     We use spline spaces of degree $p=3$ and set the order of preservation of polynomial moments $r=p$ to the spline degree, simulating for a total of $T=2$ time units with a time step $\Delta t = 0.006$.
    Observing the discrete solution over time, we see that as
    the Gaussian pulse leaves the interior patch,  we get
    spurious reflections if the conforming projection operators do not preserve polynomial moments, 
    as shown in Figure~\ref{fig:td-maxwell-comp}.
    These reflections are significantly reduced and virtually vanishing 
    if the conforming projection operators preserve polynomial moments.
    As we experimented with many grid configurations, this is a common behavior we observed, when 
    an edge is shared by patches of different resolution. 
    \bigskip 

    \begin{figure}[!htb]
        \centering
        \includegraphics[width=0.3\columnwidth]{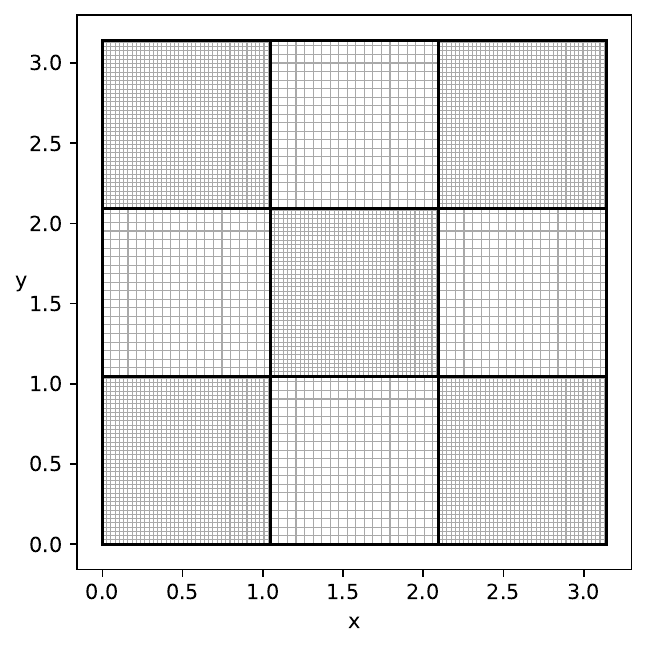}
        \includegraphics[width=0.3\columnwidth]{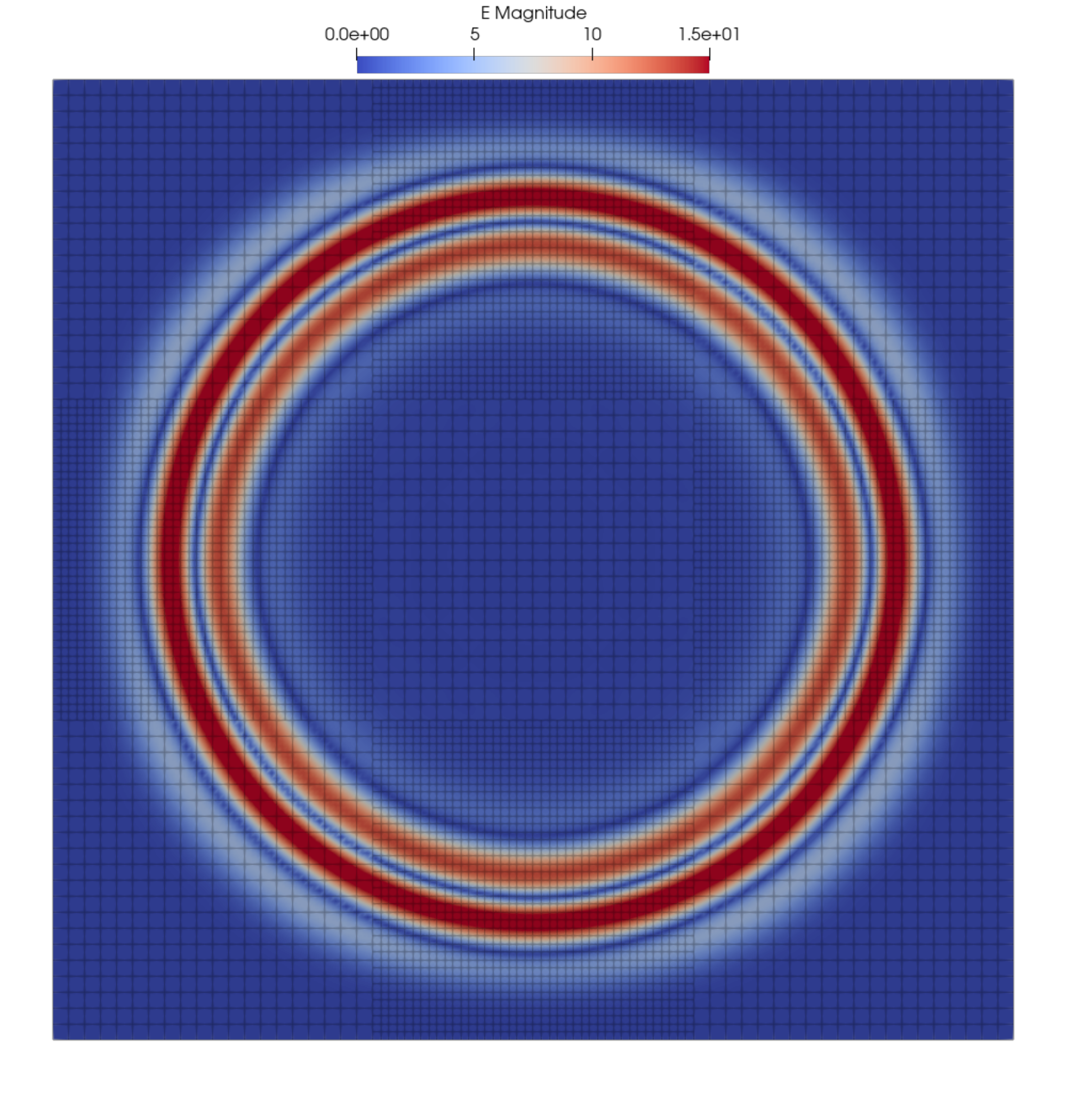}
        \includegraphics[width=0.3\columnwidth]{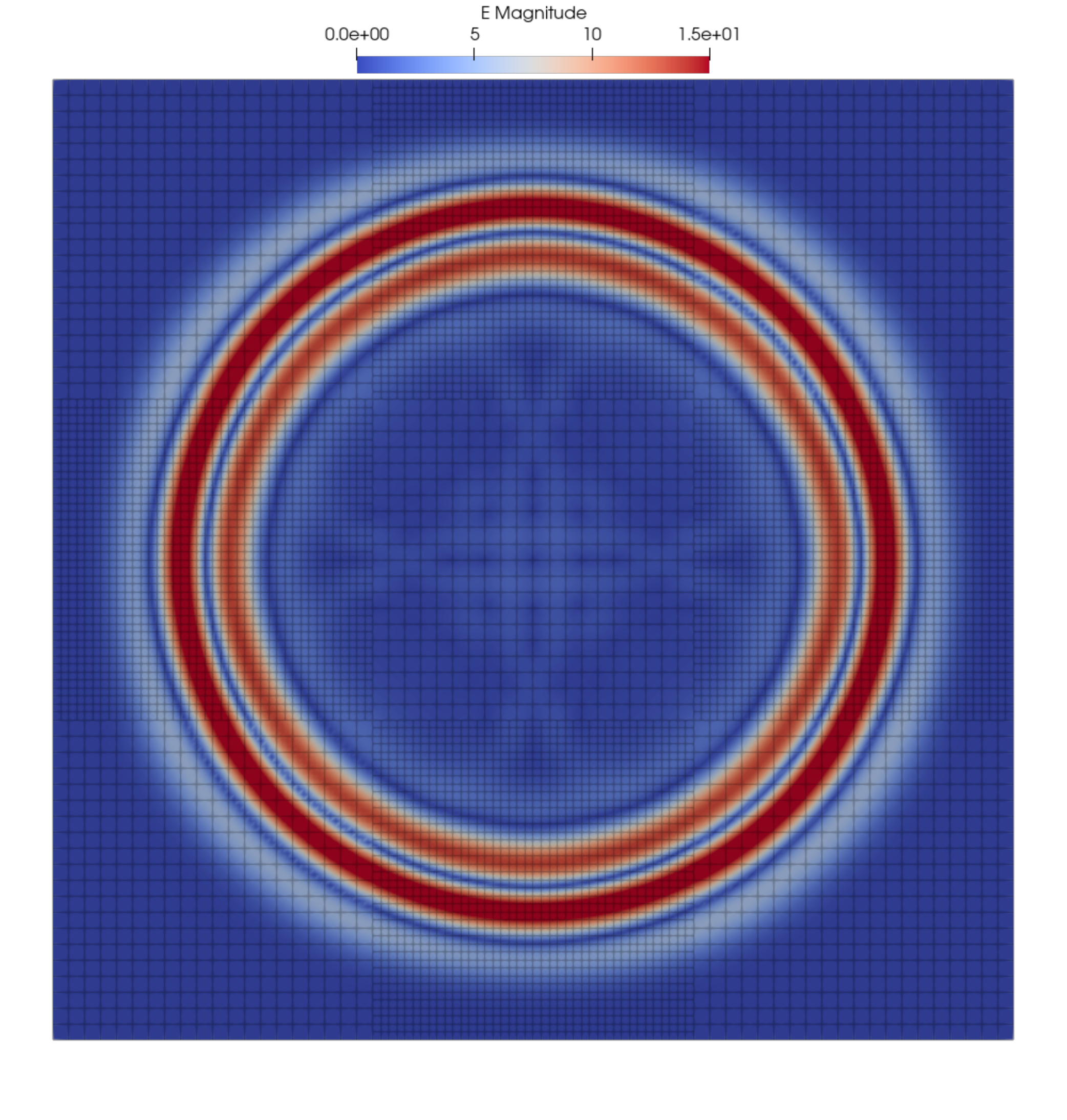}
        \caption{The checkerboard grid setup (left) and the amplitude of the electric field 
        $\bE_h$ of the time-domain Maxwell problem after leaving the interior patch. The 
        discretization in the middle involves a conforming projection 
        that preserves polynomial moments up to the spline degree,
         while the discretization on the right does not.}
        \label{fig:checkerboard-grid}
    \end{figure}

    Another grid setup that falls beyond the a priori stability analysis of \cite{BCP}
    is the case of refinements in a checkerboard pattern.
    In Figure~\ref{fig:checkerboard-grid},
    we show the electric field of the time-domain Maxwell problem on a grid with a checkerboard refinement pattern, 
    which has the same behavior as the two grid setups before.

\clearpage 

\section{Conclusion}
\label{sec:conclusion}
  This work extends the broken-FEEC approach of \cite{guclu_broken_2023} to non-matching interfaces, 
 which allows for a more flexible (through patch-local refinements) and efficient treatment of multipatch spaces. 
 The construction provides conforming projection operators for the 
 patch-wise multipatch spaces $V^0_\pw$ and $V^1_\pw$, with prescribed preservation of polynomial moments, and results in broken-FEEC de Rham sequences with accurate strong and weak derivatives.

In conclusion, the practical viability and efficiency of broken-FEEC algorithms using our novel conforming projection operators
in the numerical treatment of different PDE problems are demonstrated.
In particular, thanks to the preservation of polynomial moments, the reflections/oscillations in the context of electromagnetic wave equations are significantly reduced. 

The present construction is restricted to two dimensions, a common spline degree, and nested interface spaces obtained by refinement. Extending the approach to non-
nested trace spaces or varying degrees, and establishing uniform stability and error estimates for more general patch configurations, are natural directions for
further work. In three dimensions, conformity must be enforced simultaneously on patch faces, edges, and vertices. The correction operators must then preserve
bivariate polynomial moments on faces while remaining compatible with the face–edge–vertex incidence structure. Establishing the compatibility, locality, and
the required projection properties is the main challenge in such an extension.

\newpage 

\bibliographystyle{amsplain}
\bibliography{refs.bib}

\appendix
\section{Broken-FEEC discretizations of model PDEs}
\label{sec:disc-examples}

For the numerical verification of our conforming projection operators, 
we revisit the test problems and numerical schemes of the 
broken-FEEC framework from \cite[Section 3]{guclu_broken_2023}.

\subsection{Broken degrees of freedom}
Following \cite[Section 2.4]{guclu_broken_2023}, we assume that each local space $V^\ell_k$ is equipped with 
linear forms
\begin{equation} \label{eq:sigma}
	\sigma^{\ell, k}_\bi: V^\ell_k \rightarrow \mathbb{R}, \quad \bi \in \bcI^{\ell, k},
\end{equation}
which are characterized by duality to the basis functions:
\begin{equation}
	\sigma^{\ell, k}_\bi \left(  \Lambda^{\ell, k}_\bj \right) = \delta_{\bi,\bj}, \quad \bi, \bj \in \bcI^{\ell, k}.
\end{equation}
In the numerical examples of Section~\ref{sec:num_ex}, we will use B-spline basis functions so that $\sigma$ in \eqref{eq:sigma} are their corresponding coefficients.

This allows for the description of the differential operators in matrix form onto the coefficients, also called degrees of freedom (DOF), of the patch-local basis functions. For ease of notation, we will reduce/flatten the tuple of patch index $k \in \cK$ and multi-index $\bi \in \bcI^{\ell,k}$ to a single index 
\begin{equation} \label{eq:dof_ind}
	m : (k, \bi) \mapsto m \in \{1, \dots, N^\ell \coloneqq \dim(V^\ell_\pw) \},
\end{equation}
such that $\Lambda^\ell_m = \Lambda^{\ell, k}_\bi$.

On the patch-wise (broken) spaces, the differential operators $\bgrad$ and $\curl$ are represented by their patch-wise versions (which consist of summing every patch-local differential operator), whose operator matrices are of the form
\begin{align}
	\arrG_{i, j} &= \sigma_i^1(\bgrad_\pw \Lambda_j^0), \quad i = 1, \dots, N^1 j = 1, \dots, N^0,\\
	\arrC_{i, j} &= \sigma_i^2(\curl_\pw \Lambda_j^1),\quad i = 1, \dots, N^2, j = 1, \dots, N^1.
\end{align}
These matrices are of sizes $N^1 \times N^0$ and $N^2 \times N^1$ respectively with a patch-diagonal structure, that means the matrices are block-diagonal with blocks corresponding to the patch-local differential operator.

The conforming projection \eqref{intro:confproj} as an endomorphism on $V^0_\pw$ in matrix form reads
\begin{equation}
	\arrP^\ell_{i,j} = \sigma_i^\ell( P^\ell \Lambda^\ell_j), \quad i, j = 1, \dots, N^\ell.
\end{equation}
Since the conforming projection $P^\ell$ maps into the conforming space $V^\ell_h$, the matrix $\arrP^\ell$ is not block diagonal, but has to average/couple degrees of freedom on shared interfaces.

Introducing the coefficient spaces $\coeffs^\ell = \mathbb{R}^{N^\ell}$, this gives the commuting diagram:  
\begin{equation} \label{intro:comm_diag_disc}
    \begin{tikzpicture}[ampersand replacement=\&, baseline] 
        \matrix (m) [matrix of math nodes,row sep=3em,column sep=3em,minimum width=2em] {
            ~~ V^0_\pw ~ \bbb
            \& ~~ V^1_\pw ~ \bbb
                \& ~~ V^2_\pw ~ \bbb
    \\
	~~ \coeffs^0 ~ \bbb
            \& ~~ \coeffs^1 ~ \bbb
                \& ~~ \coeffs^2 ~ \bbb 
				\\
        };
        \path[-stealth]
        (m-1-1) edge node [above] {$\bgrad P^0$} (m-1-2)
        (m-1-2) edge node [above] {$\curl P^1$} (m-1-3)
        (m-1-1) edge node [right] {$\sigma^0$} (m-2-1)
        (m-1-2) edge node [right] {$\sigma^1$} (m-2-2)
        (m-1-3) edge node [right] {$\sigma^2$} (m-2-3)
		(m-2-1) edge node [above] {$\arrG \arrP^0$} (m-2-2)
        (m-2-2) edge node [above] {$\arrC \arrP^1$} (m-2-3)
        ;
        \end{tikzpicture}
\end{equation}

To write the coderivative matrices, we also need the mass matrices 
\begin{equation}
	\arrM^\ell_{i,j} = \sprod{\Lambda_i^\ell}{\Lambda_j^\ell}_{L^2(\Omega)}, \quad i, j = 1, \dots, N^\ell,
\end{equation}
which are block-diagonal as well.
Another matrix that we will use is the identity matrix $\arrI$, which is the identity operator on the coefficient space $\coeffs^\ell$ depending on the context.

\subsection{Weak differential operators} \label{subsec:wdiff}
Using the definitions in Section~\ref{sec:mom-pres}, we define the weak $\bcurl$ operator 
in the broken-FEEC framework as
\begin{equation}
	\tilde \arrC = (\arrM^1)^{-1} \left(\arrC \arrP^1\right)^T \arrM^2.
\end{equation}
The matrix form of the weak broken-FEEC $\Div$ operator can be derived analogously.

In Section~\ref{sec:weak-div-num-ex}, we verify its approximation properties 
when applying it to smooth test functions.

 \subsection{Poisson problem} \label{subsec:poisson}
Following \cite[Section 3.1]{guclu_broken_2023}, we consider the Poisson problem with homogeneous boundary conditions, i.e.
for given
$f \in L^2(\Omega)$, find $\phi \in V^0 = H^1_0(\Omega)$ such that
\begin{equation}
	-\Delta \phi = f.
\end{equation}
Within the broken-FEEC Framework, this equation is discretized as
\begin{equation}
	\left( \left(\arrG \arrP^0 \right)^T \arrM^1 \arrG \arrP^0 + \alpha \underbrace{ (\arrI  -  \arrP^0)^T \arrM^0  (\arrI -  \arrP^0)}_{\text{jump stabilization}} \right) \ \arr{\phi} =   (\arrP^0)^T \arrM^0 \arr{f},
\end{equation}
for some $\alpha > 0$, the coefficients $\arr{f} \in \coeffs^0$ of a projection of the source $f$
and the coefficients $\arr{\phi} \in \coeffs^0$ of the solution $\phi$. 
Notice that the right-hand-side is obtained by a filtered $L^2$-projection similar to the ones in Section~\ref{sec:mom-pres}.
Here we remind that the jump stabilization term is necessary to ensure the well-posedness of the discrete problem, and the precise form of the stabilization guarantees that for any $\alpha > 0$, the solution coincides with the standard conforming solution on the conforming subspace $V^0_h$ of $V^0_\pw$, see \cite[Sect~2.7]{guclu_broken_2023}.

\subsection{Time-harmonic Maxwell problem}\label{subsec:thmaxwell}
To study the time-harmonic Maxwell problem
with inhomogeneous boundary conditions, we again follow \cite[Section 3.3]{guclu_broken_2023}:
For $\omega \in \mathbb{R}$ and $\bJ \in L^2(\Omega)$, we solve
\begin{equation}
	- \omega^2 \bE + \bcurl \curl \bE  = \bJ,
\end{equation}
for  $\bE \in V^1 = H(\curl; \Omega)$.  
In our framework, that reads
\begin{equation}
	\left(- \omega^2 (\arrP^1)^T \arrM ^1 \arrP^1 + (\arrC \arrP^1 )^T \arrM^2 \arrC \arrP^1 + \alpha \underbrace{ (\arrI -  \arrP^1)^T \arrM^1  (\arrI -  \arrP^1)}_{\text{jump stabilization}}  \right)  \ \arre = (\arrP^1)^T \arrM^1 \arrj,
\end{equation}
for some $\alpha > 0$, with the discrete source coefficients $\arrj \in \coeffs^1$ and the solution coefficients $\arre \in \coeffs^1$.
Here as well, the jump stabilization term guarantees the well-posedness of the discrete problem, and that for any $\alpha > 0$, the solution coincides with the standard conforming solution on the conforming subspace $V^1_h$ of $V^1_\pw$, see \cite{guclu_broken_2023}.

\subsection{Curl-curl eigenvalue problem}\label{subsec:ccev}
    Following \cite[Section 3.4 and Section 5.5]{guclu_broken_2023}, we look at the spectral
     correctness of the $\bcurl$-$\curl$ eigenvalue problem:
     Find $\lambda > 0$ such that 
    \begin{equation} 
         \bcurl \curl \mathbf{E} = \lambda \mathbf{E}, 
        \end{equation}
    where $\mathbf{E} \in V^1 = H_0(\curl; \Omega)$.

    Within the broken-FEEC framework, this equation is discretized as
    \begin{equation*}
        \left(\arrC \arrP ^1 \right)^T \arrM ^2 \arrC \arrP^1 \ \arre = \lambda \left( \alpha \underbrace{(\arrI  -  \arrP^1)^T \arrM^1  (\arrI -  \arrP^1)}_{\text{jump stabilization}} +    (\arrP^1)^T \arrM^1 \arrP ^1 \right) \  \arre,
    \end{equation*}
    for $\alpha > 0$ and with the coefficients $\arre \in \coeffs^1$.
    Again the form of the stabilization guarantees that the broken and conforming eigenfunctions coincide for any $\alpha > 0$, see \cite{guclu_broken_2023}.

\subsection{Time-domain Maxwell} \label{subsec:tdmaxwell}
We consider the time-domain Maxwell problem on a square domain as described in \cite[Section 3.6 and Section 5.7]{guclu_broken_2023}: 
For  $\boldsymbol{J} \in L^2(\Omega)$, solve
	\begin{align*}
		\partial_t \boldsymbol{E} - \bcurl B &= - \boldsymbol{J},\\
		\partial_t B + \curl \boldsymbol{E} &= 0,
	\end{align*}
where  $\boldsymbol{E} \in H_0(\curl; \Omega)$ and $B \in L^2(\Omega)$. 
A semi-discrete scheme in the broken-FEEC framework then reads
\begin{align} \label{eq:tdmaxwell-semdisc}
	\arrM^1 \partial_t \arre  &= (\arrC \arrP^1) ^T \arrM^2 \arrb - (\arrP^1)^T \arrM^1 \arrj,\\
	\partial_t \arrb &= - \arrC \arrP^1 \arre,
\end{align}
with  coefficients $\arrb \in \coeffs^2$, $\arre \in \coeffs^1$ and $\arrj \in \coeffs^1$.
This is discretized in time using a leapfrog scheme. 
A similar approach can be used to solve the time-domain Helmholtz problem.

\newpage 

\section{Homogeneous boundary conditions}
\label{app:homogeneous}
In order to use our scheme for sequences with homogeneous boundary conditions, we can extend our construction. This boils down to do modifying the edge- and vertex-conforming projection operators for edges and vertices that lie on the domain boundary. 

First, we look at the $V^0$ case discussed in Section~\ref{sec:conf-proj}
and extend Definition~\ref{eq:Pv} to vertices on the boundary: 
Let $\vertex \in \vertices$ be a boundary vertex, then 
\begin{equation} 
    P_\vertex: \Lambda^k_\bj \mapsto \left\{ \begin{aligned} &\begin{aligned}
                            \sum_{m_1 = 1}^r \sum_{m_2 = 1}^r  \bgamma_{\mb}   \Lambda^{k}_{\mb}\\
                                \end{aligned} &&\text{ if } \bj = \bzero \text{ and } k \in \cK(\vertex), \\
    \end{aligned}\right.
\end{equation} 
where the point-value at the vertex is projected to zero, but we still have to account for the preservation of polynomial moments by interior basis functions close to the vertex. 

Continuing with the edge-based projections, we 
similarly extend Definition~\ref{eq:Pe} to edges on the boundary: 
Let $k$ be the patch adjacent to the boundary edge $\edge$

\begin{align}
    P_\edge : 
        \Lambda^{k}_\jb \mapsto \left\{ \begin{aligned}
                                    & \begin{aligned} & \cF^0 \left( \sum_{m = 1}^r  \gammaod_{m}  \lambda^{k}_{j_1} \otimes   \lambda^{k}_{m}  \right)
                                    \end{aligned}
                                                &&\text{ if } 0 \le j_1 \le n_k,~ j_2= 0, \\
                                    & \begin{aligned} &\cF^0 \left(  \sum_{m = 1}^r  \gammaod_{ m}  \mu_\vertex^- \otimes   \lambda^{k}_{m} \right)
                                    \end{aligned}
                                                &&\begin{aligned} 
                                                    \text{ if }& j_1 \in \{ 0, n_-\},~ j_2 = 0 \\
                                                    &\text{ for }
                                                     \vertex \text{ at } \jb 
                                                    \text{ and } \vertex \text{ shared by another patch}, 
                                                \end{aligned}\\ 
                                        &\begin{aligned} &\cF^0 \left(  \sum_{m_1 = 1}^r \sum_{m_2 = 1}^r  \bgamma_{\mb}   \Lambda^{k}_{\mb} \right)
                                                \end{aligned}
                                                            &&\text{ if }\jb = \bzero \text{ for } \vertex \text{ at } \jb \text{ and } \vertex \text{ is only on patch } k, \\
                                        \end{aligned}\right. 
\end{align}
where patch $k$ adjacent to the edge $\edge$ is treated as the coarse patch. Again, we project all point-values on the boundary to zero, but have to account for moment preservation similarly to before. 


This extension directly applies to the $V^1$ case from Section~\ref{sec:conf-proj-v1}. The edge-conforming projection in Definition~\ref{eq:Pe-v1} for a boundary edge $\edge$ in this case reads: 

\begin{align}  
    \hat P^1_\edge :
        \begin{aligned}
            &(\hat \bLambda^{1, -}_{1, \jb})_1 \mapsto  \begin{aligned}
                                    & \begin{aligned} & \sum_{m = 1}^r  \gammaod_{ m} \lambda^{1,k}_{j_1} \otimes \lambda^{k}_{m}            
                                    \end{aligned}
                                                &&\text{ if }  j_2= 0,\\
                                \end{aligned}
        \end{aligned}
\end{align}
with arguments similar to above. 

\newpage 
\section{Implementation of the conforming projection operators}
\label{app:implementation}
\subsection{Correction coefficients \texorpdfstring{$\gamma$}{gamma}}
\label{sec:impl_gamma}
As we discussed in Remark~\ref{rem:patch-indep}, we can solve the linear system for the correction coefficients
 $\ba$ and $\tilde \ba$ from Ansatz~\eqref{eq:ansatz-pv} independent of the patch. 
Similar to Section~\ref{sec:p_vertex}, we define the basis-polynomial duality-matrices and vertex-basis-polynomial vectors as 
\begin{align} \label{eq:pol_mat_const}
    \arrM = \left(\int_0^1 \lambda_{m} q_{j} \right)_{ 1 \le j \le r, 1 \le m \le r}, \quad \arrb &= \left( \int_0^1  \lambda_0 q_{j} \right)_{1 \le j \le r}.
\end{align}
For the polynomial basis on the interval $[0,1]$, we use Bernstein polynomials
\begin{equation} \label{eq:bernstein}
    q_j(\hat x) = \binom{r-1}{j-1} \hat x^{j-1}(1- \hat x)^{r-j}, \quad j = 1, \dots, r,
\end{equation} 
since the resulting duality-matrices have a good condition number in practice.
Thus, we define the one-dimensional correction coefficient 
\begin{equation}
    \gammaod =  \arrM^{-1} \arrb.
\end{equation}

\subsection{Matrix form of the extension- and restriction operators}
\label{sec:impl_erv0}
Since both spaces have the same degree and are nested, we use the knot insertion technique for B-splines, explained for example in \cite[Section 5.2]{NURBS},
to describe the extension operator $\scrE$ from \eqref{eq:extension} as the matrix $\arrE \in \RR^{n_+ +1\times n_- +1}$ acting on the degrees of freedom/coefficients in these basis representations, i.e. 
\begin{equation}
    \lambda^{-}_i = \sum_{j = 0}^{n_+} \arrE_{j, i} \lambda^{+}_j,
\end{equation}
for $i = 0, \dots, n_-$, where $n_\pm = \dim(\VV^0_\pm)-1$.
This form of the extension operator is used directly on the coefficients of the coarse space.

The restriction operator $\scrR$ in \eqref{eq:restr} consists of 
\begin{itemize}
\item The interior restriction operator $\tilde \scrR$ is defined 
by the $L^2$-projection in \eqref{eq:l2-r0}. We introduce the interior (mixed) mass matrices
\begin{align}
   \tilde \arrM^{-, +}_{i,j} &= \left( \langle \lambda^{-}_i, \lambda^{+}_j \rangle \right)_{i,j}, \quad i = 0, \dots, n_-, \ j = 1, \dots, n_+-1, \\
    \arrM^{-, -}_{i,j} &= \left( \langle \lambda^{-}_i, \lambda^{-}_j \rangle \right)_{i,j}, \quad i = 0, \dots, n_-, \ j = 0, \dots, n_-,
\end{align}
so that the projection has the matrix form  
\begin{equation}
    \tilde \arrR = (\arrM^{-, -})^{-1} \tilde \arrM^{-, +}.
\end{equation}
It is fully moment-preserving by definition of the $L^2$-projection. 

\item The coarse truncation $\tilde \cT$ in \eqref{eq:trunc-2}, with matrix form $\tilde \arrT  \in \RR^{n_--1 \times n_-+1}$ given by
        \begin{equation}
            \tilde \arrT_{i, j} = \left\{ \begin{aligned}
                &1 \quad &&\text{ for } 0< i,j<n_-, j = i, \\
                &\gammaod_i&& \text{ for } 0<i<r+1, j = 0, \\
                & \gammaod_{n_- - i} && \text{ for } n_--r-1<i<n_-, j = n_-, \\
                &0 \quad &&\text{ else.}
            \end{aligned} \right.
        \end{equation}
        The motivation for this matrix is given by looking at $\phi^- \in \VV^{0}_{-, \parallel}$ 
        and mapping $\tilde \cT \phi^- = \tilde \phi^- \in  \VV^{0}_{-, \parallel, 0}$ such that $(\tilde \phi^- )_0 = (\tilde \phi^- )_ {n_-} = 0$, 
        while preserving polynomial moments: 
        \begin{align}
            \int \tilde \phi^- q &= \int \left( \sum_{j = 1}^{n_- -1} \phi^-_j \lambda^-_j + \sum_{i = 1}^{r} \gammaod_i \phi^-_0 \lambda^-_i + \sum_{i = n_--r}^{n_- -1} \gammaod_{n_- -i} \phi^-_{n_-} \lambda^-_i \right) q = \int \phi^- q,
        \end{align}
        by the definition of $\gammaod$ in \eqref{eq:gamma1d}.

\item The truncation $\cT$ in \eqref{eq:trunc-1} in matrix form $\arrT  \in \RR^{n_+-1 \times n_+ +1}$ can be written as
        \begin{equation}
            (\arrT)_{i, j} = \left\{ \begin{aligned}
                &1 \quad &&\text{ for } j = i, 1 \le i \le n_+-1\\
                &-\arrE_{0, i} && \text{ for } j = 0 , 1 \le i \le n_+-1, \\
                &-\arrE_{n_+, i} && \text{ for } j = n_+, 1 \le i \le n_+-1, \\
                &0 \quad &&\text{ else.}
            \end{aligned} \right.
        \end{equation}
\end{itemize}

This leads to the matrix form of the moment-preserving interior restriction operator in \eqref{eq:int_rest} as
\begin{equation}
    \arrR_0 = \tilde \arrT \tilde \arrR \in \RR^{n_--1 \times n_+-1}.
\end{equation}
The full restriction matrix $\arrR$ is therefore
\begin{equation}
    \arrR =  \left(\begin{array}{cccc}
        1 & 0 & \cdots & 0 \\ 
        \hline
        \multicolumn{4}{c}{\arrR_0 \arrT} \\
         \hline 
        0 & \cdots & 0 & 1 \\ 
      \end{array}\right) \in \RR^{n_- +1\times n_+ +1}
\end{equation}

\begin{remark}
    Instead of using an $L^2$-projection for the restriction operator $\tilde \scrR$, 
    we could also choose a more local projection. 
\end{remark}

For the $V^1_\pw$ case, we can use the same approach as for $V^0_\pw$, since the nestedness still holds for the derived spaces. The extension operator is a knot insertion operator, whereas now the restriction operator is purely defined as the $L^2$-projection and thus already moment-preserving by definition.


\end{document}